\documentclass[12pt,oneside,english,american]{amsart}
\usepackage{libertine}
\usepackage[T1]{fontenc}
\usepackage[utf8]{inputenc}
\usepackage[active]{srcltx}
\usepackage{babel}
\usepackage{verbatim}
\usepackage{prettyref}
\usepackage{float}
\usepackage{mathrsfs}
\usepackage{mathtools}
\usepackage{varwidth}
\usepackage{amsbsy}
\usepackage{amstext}
\usepackage{amsthm}
\usepackage{amssymb}
\usepackage{graphicx}
\usepackage[a4paper]{geometry}
\usepackage{setspace}
\usepackage[authoryear]{natbib}
\usepackage[bookmarks=true,bookmarksnumbered=false,bookmarksopen=false,
 breaklinks=false,pdfborder={0 0 1},backref=false,colorlinks=false]
 {hyperref}
\hypersetup{pdftitle={Inference on eigenvectors of non-symmetric matrices},
 pdfauthor={Jerome R. Simons},
 hidelinks}

\makeatletter

\providecommand{\tabularnewline}{\\}
\newenvironment{cellvarwidth}[1][t]
    {\begin{varwidth}[#1]{\linewidth}}
    {\@finalstrut\@arstrutbox\end{varwidth}}

\floatstyle{ruled}
\newfloat{algorithm}{tbp}{loa}
\providecommand{\algorithmname}{Algorithm}
\floatname{algorithm}{\protect\algorithmname}

\numberwithin{equation}{section}
\numberwithin{figure}{section}
\numberwithin{table}{section}
\theoremstyle{definition}
 \newtheorem{example}{\protect\examplename}
\theoremstyle{plain}
\newtheorem{assumption}{\protect\assumptionname}
\theoremstyle{plain}
\newtheorem{thm}{\protect\theoremname}
\theoremstyle{plain}
\newtheorem{cor}{\protect\corollaryname}
\theoremstyle{plain}
\newtheorem{lem}{\protect\lemmaname}
\theoremstyle{remark}
\newtheorem{rem}{\protect\remarkname}

\newrefformat{eq}{\textup{(\ref{#1})}} 
\newrefformat{lem}{Lemma \ref{#1}}
\newrefformat{thm}{Theorem \ref{#1}}
\newrefformat{cha}{Chapter \ref{#1}}
\newrefformat{sec}{Section \ref{#1}}
\newrefformat{assu}{Assumption \ref{#1}}
\newrefformat{prop}{Proposition \ref{#1}}
\newrefformat{cor}{Corollary \ref{#1}}
\newrefformat{enu}{.\ref{#1}}
\newrefformat{exa}{Example \ref{#1}}
\newrefformat{subsec}{Section \ref{#1}}
\newrefformat{tab}{Table \ref{#1}}
\newrefformat{alg}{Algorithm \ref{#1}}
\newrefformat{fig}{Figure \ref{#1}}
\newrefformat{app}{Appendix \ref{#1}}
\usepackage{changepage}

\usepackage{svg}
\usepackage{tikz}
\usetikzlibrary{cd}
\tikzset{
  myls/.style={
    draw, text height=height("$d$")
    }
  }

\usepackage{csquotes}
\usepackage{ifpdf}
\usepackage{ragged2e}

\usepackage{booktabs}
\usepackage{hyperref}

\usepackage{algorithmic}

\makeatletter 
 
\@addtoreset{algorithm}{section} 
\makeatother

\usepackage{appendix}

\makeatother

\addto\captionsamerican{\renewcommand{\assumptionname}{Assumption}}
\addto\captionsamerican{\renewcommand{\corollaryname}{Corollary}}
\addto\captionsamerican{\renewcommand{\examplename}{Example}}
\addto\captionsamerican{\renewcommand{\lemmaname}{Lemma}}
\addto\captionsamerican{\renewcommand{\remarkname}{Remark}}
\addto\captionsamerican{\renewcommand{\theoremname}{Theorem}}
\addto\captionsenglish{\renewcommand{\algorithmname}{Algorithm}}
\addto\captionsenglish{\renewcommand{\assumptionname}{Assumption}}
\addto\captionsenglish{\renewcommand{\corollaryname}{Corollary}}
\addto\captionsenglish{\renewcommand{\examplename}{Example}}
\addto\captionsenglish{\renewcommand{\lemmaname}{Lemma}}
\addto\captionsenglish{\renewcommand{\remarkname}{Remark}}
\addto\captionsenglish{\renewcommand{\theoremname}{Theorem}}
\providecommand{\assumptionname}{Assumption}
\providecommand{\corollaryname}{Corollary}
\providecommand{\examplename}{Example}
\providecommand{\lemmaname}{Lemma}
\providecommand{\remarkname}{Remark}
\providecommand{\theoremname}{Theorem}

\begin{document}


\global\long\def\uwrite#1#2{\underset{#2}{\underbrace{#1}} }%
\global\long\def\blw#1{\ensuremath{\underline{#1}}}%
\global\long\def\abv#1{\ensuremath{\overline{#1}}}%
\global\long\def\vect#1{\mathbf{#1}}%


\global\long\def\smlseq#1{\{#1\} }%
\global\long\def\seq#1{\left\{  #1\right\}  }%
\global\long\def\smlsetof#1#2{\{#1\mid#2\} }%
\global\long\def\setof#1#2{\left\{  #1\mid#2\right\}  }%


\global\long\def\goesto{\ensuremath{\rightarrow}}%
\global\long\def\ngoesto{\ensuremath{\nrightarrow}}%
\global\long\def\uto{\ensuremath{\uparrow}}%
\global\long\def\dto{\ensuremath{\downarrow}}%
\global\long\def\uuto{\ensuremath{\upuparrows}}%
\global\long\def\ddto{\ensuremath{\downdownarrows}}%
\global\long\def\ulrto{\ensuremath{\nearrow}}%
\global\long\def\dlrto{\ensuremath{\searrow}}%


\global\long\def\setmap{\ensuremath{\rightarrow}}%
\global\long\def\elmap{\ensuremath{\mapsto}}%
\global\long\def\compose{\ensuremath{\circ}}%
\global\long\def\cont{C}%
\global\long\def\cadlag{D}%
\global\long\def\Ellp#1{\ensuremath{\mathcal{L}^{#1}}}%


\global\long\def\naturals{\ensuremath{\mathbb{N}}}%
\global\long\def\reals{\mathbb{R}}%
\global\long\def\complex{\mathbb{C}}%
\global\long\def\rationals{\mathbb{Q}}%
\global\long\def\integers{\mathbb{Z}}%


\global\long\def\abs#1{\ensuremath{\left|#1\right|}}%
\global\long\def\smlabs#1{\ensuremath{\lvert#1\rvert}}%
\global\long\def\bigabs#1{\ensuremath{\bigl|#1\bigr|}}%
\global\long\def\Bigabs#1{\ensuremath{\Bigl|#1\Bigr|}}%
\global\long\def\biggabs#1{\ensuremath{\biggl|#1\biggr|}}%
\global\long\def\norm#1{\ensuremath{\left\Vert #1\right\Vert }}%
\global\long\def\smlnorm#1{\ensuremath{\lVert#1\rVert}}%
\global\long\def\bignorm#1{\ensuremath{\bigl\|#1\bigr\|}}%
\global\long\def\Bignorm#1{\ensuremath{\Bigl\|#1\Bigr\|}}%
\global\long\def\biggnorm#1{\ensuremath{\biggl\|#1\biggr\|}}%


\global\long\def\Union{\ensuremath{\bigcup}}%
\global\long\def\Intsect{\ensuremath{\bigcap}}%
\global\long\def\union{\ensuremath{\cup}}%
\global\long\def\intsect{\ensuremath{\cap}}%
\global\long\def\pset{\ensuremath{\mathcal{P}}}%
\global\long\def\clsr#1{\ensuremath{\overline{#1}}}%
\global\long\def\symd{\ensuremath{\Delta}}%
\global\long\def\intr{\operatorname{int}}%
\global\long\def\cprod{\otimes}%
\global\long\def\Cprod{\bigotimes}%


\global\long\def\smlinprd#1#2{\ensuremath{\langle#1,#2\rangle}}%
\global\long\def\inprd#1#2{\ensuremath{\left\langle #1,#2\right\rangle }}%
\global\long\def\orthog{\ensuremath{\perp}}%
\global\long\def\dirsum{\ensuremath{\oplus}}%


\global\long\def\spn{\operatorname{sp}}%
\global\long\def\rank{\operatorname{rk}}%
\global\long\def\proj{\operatorname{proj}}%
\global\long\def\tr{\operatorname{tr}}%


\global\long\def\smpl{\ensuremath{\Omega}}%
\global\long\def\elsmp{\ensuremath{\omega}}%
\global\long\def\sigf#1{\mathcal{#1}}%
\global\long\def\sigfield{\ensuremath{\mathcal{F}}}%
\global\long\def\sigfieldg{\ensuremath{\mathcal{G}}}%
\global\long\def\flt#1{\mathcal{#1}}%
\global\long\def\filt{\mathcal{F}}%
\global\long\def\filtg{\mathcal{G}}%
\global\long\def\Borel{\ensuremath{\mathcal{B}}}%
\global\long\def\cyl{\ensuremath{\mathcal{C}}}%
\global\long\def\nulls{\ensuremath{\mathcal{N}}}%
\global\long\def\gauss{\mathfrak{g}}%
\global\long\def\leb{\mathfrak{m}}%

\global\long\def\prob{P}%
\global\long\def\Prob{\ensuremath{\mathbb{P}}}%
\global\long\def\Probs{\mathcal{P}}%
\global\long\def\PROBS{\mathcal{M}}%
\global\long\def\expect{\ensuremath{\mathbb{E}}}%
\global\long\def\probspc{\ensuremath{(\smpl,\filt,\Prob)}}%

\global\long\def\iid{\ensuremath{\textnormal{i.i.d.}}}%
\global\long\def\as{\ensuremath{\textnormal{a.s.}}}%
\global\long\def\asp{\ensuremath{\textnormal{a.s.p.}}}%
\global\long\def\io{\ensuremath{\ensuremath{\textnormal{i.o.}}}}%
\global\long\def\independent{\mathpalette{\independenT}{\perp}}%
 
\global\long\def\independenT#1#2{\mathrel{\rlap{$#1#2$}\mkern2mu {#1#2}}}%

\global\long\def\indep{\independent}%
\global\long\def\distrib{\ensuremath{\sim}}%
\global\long\def\distiid{\ensuremath{\sim_{\iid}}}%
\global\long\def\asydist{\ensuremath{\overset{a}{\distrib}}}%
\global\long\def\inprob{\ensuremath{\overset{p}{\goesto}}}%
\global\long\def\inprobu#1{\ensuremath{\overset{#1}{\goesto}}}%
\global\long\def\inas{\ensuremath{\overset{\as}{\goesto}}}%
\global\long\def\eqas{=_{\as}}%
\global\long\def\inLp#1{\ensuremath{\overset{\Ellp{#1}}{\goesto}}}%
\global\long\def\indist{\ensuremath{\overset{d}{\goesto}}}%
\global\long\def\eqdist{=_{d}}%
\global\long\def\wkc{\ensuremath{\rightsquigarrow}}%
\global\long\def\wkcu#1{\overset{#1}{\ensuremath{\rightsquigarrow}}}%
\global\long\def\plim{\operatorname*{plim}}%


\global\long\def\var{\operatorname{var}}%
\global\long\def\lrvar{\operatorname{lrvar}}%
\global\long\def\cov{\operatorname{cov}}%
\global\long\def\corr{\operatorname{corr}}%
\global\long\def\bias{\operatorname{bias}}%
\global\long\def\MSE{\operatorname{MSE}}%
\global\long\def\med{\operatorname{med}}%


\global\long\def\simple{\mathcal{R}}%
\global\long\def\sring{\mathcal{A}}%
\global\long\def\sproc{\mathcal{H}}%
\global\long\def\Wiener{\ensuremath{\mathbb{W}}}%
\global\long\def\sint{\bullet}%
\global\long\def\cv#1{\left\langle #1\right\rangle }%
\global\long\def\smlcv#1{\langle#1\rangle}%
\global\long\def\qv#1{\left[#1\right]}%
\global\long\def\smlqv#1{[#1]}%


\global\long\def\trans{\ensuremath{\prime}}%
\global\long\def\indic{\ensuremath{\mathbf{1}}}%
\global\long\def\Lagr{\mathcal{L}}%
\global\long\def\grad{\nabla}%
\global\long\def\pmin{\ensuremath{\wedge}}%
\global\long\def\Pmin{\ensuremath{\bigwedge}}%
\global\long\def\pmax{\ensuremath{\vee}}%
\global\long\def\Pmax{\ensuremath{\bigvee}}%
\global\long\def\sgn{\operatorname{sgn}}%
\global\long\def\argmin{\operatorname*{argmin}}%
\global\long\def\argmax{\operatorname*{argmax}}%
\global\long\def\Rp{\operatorname{Re}}%
\global\long\def\Ip{\operatorname{Im}}%
\global\long\def\deriv{\ensuremath{\mathrm{d}}}%
\global\long\def\diffnspc{\ensuremath{\deriv}}%
\global\long\def\diff{\ensuremath{\,\deriv}}%
\global\long\def\i{\ensuremath{\mathrm{i}}}%
\global\long\def\e{\mathrm{e}}%
\global\long\def\sep{,\ }%
\global\long\def\defeq{\coloneqq}%
\global\long\def\eqdef{\eqqcolon}%


\global\long\def\ris{\mathcal{S}}%

\global\long\def\rjsof{\mathcal{S}_{J}\left(M\right)}%

\global\long\def\risof{\mathcal{S}_{I}\left(M\right)}%

\global\long\def\schwein{\upsilon_{\perp}^{\trans}}%

\global\long\def\jacsvd{B_{\text{SVD}}}%

\global\long\def\jacsvdt{B_{\text{SVD}}^{\trans}}%

\global\long\def\eigs{\mathcal{L}}%
\global\long\def\spect{\mathcal{\lambda}}%
\global\long\def\trans{\mathsf{T}}%
\global\long\def\qcs{\textnormal{QCS}}%
\global\long\def\cs{\textnormal{CS}}%
\global\long\def\lu{{\scriptscriptstyle \textnormal{LU}}}%
\global\long\def\ur{{\scriptscriptstyle \textnormal{UR}}}%
\global\long\def\st{{\scriptscriptstyle \textnormal{ST}}}%
\global\long\def\sl{{\scriptscriptstyle \textnormal{SL}}}%
\global\long\def\ls{{\scriptscriptstyle \textnormal{LS}}}%
\global\long\def\lds{{\scriptscriptstyle \textnormal{L\ensuremath{\bullet}S}}}%
\global\long\def\sdl{{\scriptscriptstyle \textnormal{S\ensuremath{\bullet}L}}}%
\global\long\def\vek{\operatorname{vec}}%
\global\long\def\diag{\operatorname{diag}}%
\global\long\def\col{\operatorname{col}}%
\global\long\def\chol{\text{ch}}%
\global\long\def\proj{\mathcal{Q}}%
\global\long\def\smlfloor#1{\lfloor#1\rfloor}%
\global\long\def\sm{\mathfrak{sm}}%
\global\long\def\mg{\mathfrak{mg}}%
\global\long\def\M{\mathcal{M}}%
\global\long\def\N{\mathcal{N}}%
\global\long\def\wald{\mathcal{W}}%
\global\long\def\PHI{\boldsymbol{\Phi}}%
\global\long\def\like{\mathcal{\ell}}%
\global\long\def\likens{\like^{\ast}}%
\global\long\def\mn{\mathrm{MN}}%
\global\long\def\normdist{\mathrm{N}}%
\global\long\def\lr{\mathcal{LR}}%
\global\long\def\spc#1{\mathcal{#1}}%
\global\long\def\set#1{\mathscr{#1}}%
\global\long\def\err{\varepsilon}%
\global\long\def\radius{\rho}%
\global\long\def\xdet{\bar{x}}%
\global\long\def\ydet{\bar{y}}%
\global\long\def\zdet{\bar{z}}%
\global\long\def\errdet{\bar{\err}}%
\global\long\def\Zdet{\bar{Z}}%
\global\long\def\Ld{\set L_{\mathrm{d}}}%
\global\long\def\Ln{\set L_{\mathrm{n}}}%
\global\long\def\Ls{\set L_{\mathrm{s}}}%
\global\long\def\ci{\mathcal{C}}%
\global\long\def\cisimple{\ci_{\mathcal{S}}}%
\global\long\def\cibonf{\ci_{\textnormal{\emph{B}}}}%
\global\long\def\cibonfb{\ci_{\textnormal{M}}}%
\global\long\def\cibonfa{\ci_{\textnormal{A}}}%
\global\long\def\cijoh{\ci_{\textnormal{\ensuremath{J}}}}%
\global\long\def\largedec#1{\mathbf{#1}}%
\global\long\def\R{\largedec R}%
\global\long\def\G{\largedec G}%
\global\long\def\L{\largedec L}%
\global\long\def\y{\largedec y}%
\global\long\def\x{\largedec x}%
\global\long\def\BM{\mathrm{BM}}%
\global\long\def\ginv{\#}%
\global\long\def\irf{\mathrm{IRF}}%
\global\long\def\f{\varphi}%
\global\long\def\spcf{\mathcal{P}}%
\global\long\def\cone{\mathcal{K}}%
\global\long\def\snum{{\scriptscriptstyle \#}}%
\global\long\def\locest{f}%
\global\long\def\dist{\operatorname{dist}}%
\global\long\def\THETA{\boldsymbol{\Theta}}%
\global\long\def\GAMMA{\boldsymbol{\Gamma}}%
\global\long\def\PI{\boldsymbol{\Pi}}%
\global\long\def\Z{\mathcal{Z}}%
\global\long\def\VAR{\mathrm{VAR}}%
\global\long\def\MA{\mathrm{MA}}%
\global\long\def\geo{\mathrm{geo}}%
\global\long\def\alg{\mathrm{alg}}%
\global\long\def\eigs{\mathcal{L}}%
\global\long\def\lu{\mathrm{LU}}%
\global\long\def\sl{\mathrm{SL}}%
\global\long\def\st{\mathrm{ST}}%
\global\long\def\ls{\mathrm{LS}}%
\global\long\def\loc{\mathrm{LOC}}%
\global\long\def\fix{\mathrm{FIX}}%
\global\long\def\qcs{\mathrm{QCS}}%
\global\long\def\cs{\mathrm{CS}}%
\global\long\def\vek{\text{vec\ }}%
\global\long\def\ct{\mathrm{\dagger}}%
\global\long\def\Rspc{\mathrm{\mathcal{R}}}%
\global\long\def\simg{\mathcal{S}_{Q}}%
\global\long\def\simn{\mathcal{S}_{Q_{n}}}%
\global\long\def\simu{\mathcal{S}_{U}}%
\global\long\def\ex{\mathcal{X}}%

\global\long\def\ij{\text{IJ}}%

\href{https://arxiv.org/pdf/2303.18233.pdf}{Please click here for the arxiv version.}
\title{Hypothesis testing on invariant subspaces of non-diagonalizable matrices
with applications to network statistics}
\author{Jérôme R. Simons}
\date{8 October `25}
\begin{abstract}
We generalise the inference procedure for eigenvectors of symmetrizable
matrices of \citet{Tyler1981} to that of invariant and singular subspaces
of non-diagonalizable matrices. Wald tests for invariant vectors and
$t$-tests for their individual coefficients perform well in simulations,
despite the matrix being not symmetric. Using these results, it is
now possible to perform inference on network statistics that depend
on eigenvectors of non-symmetric adjacency matrices as they arise
in empirical applications from directed networks. Further, we find
that statisticians only need control over the first-order Davis-Kahan
bound to control convergence rates of invariant subspace estimators
to higher-orders. For general invariant subspaces, the minimal eigenvalue
separation dominates the first-order bound potentially slowing convergence
rates considerably. In an example, we find that accounting for uncertainty
in network estimates changes empirical conclusions about the ranking
of nodes' popularity.

\thanks{An earlier version of this paper was entitled \textquotedblleft Inference
on non-symmetric subspaces and network statistics\textquotedblright}\thanks{I would like to thank James A. Duffy, Steve Bond, Michael Leung, Richard
Samworth, David E. Tyler, Eric French, Alexei Onatskiy, Oliver Linton,
Richard Smith, Andrew Harvey, Patrick Allmis, Christian Ghiglino,
Carsten-Andreas Schulz, and Sam Gee for their comments and suggestions.
I am also grateful to the organizers and participants of Encounters
in Econometric Theory, seminars at Cambridge, Oxford, and the Summer
Meeting of the Econometric Society. I acknowledge funding from the
Keynes Fund at the Faculty of Economics, Cambridge.}
\end{abstract}

\maketitle

\section{Introduction}

This paper contributes hypothesis tests for both invariant subspace
and singular vectors of non-symmetric matrices that are not diagonalisable.
As an application of our theory, we specialise tests for a selection
of centrality and clustering statistics as they arise as functions
of network adjacency matrices. While network statistics are perhaps
most empirically relevant, our results are in the form of general
$t$- and Wald tests with the latter reducing to the procedure developed
in \citet{Tyler1981}, when the matrix has a real spectrum and is
diagonalizable.

The source of randomness in the network context is uncertainty about
the extent of weights or links where errors propagate to these statistics.
Allowing for non-symmetric adjacency matrices opens up many empirically
relevant applications. Directed, weighted networks for example arise
when weights depend on the flow direction between nodes. Specifically,
trade, input-output, and food chain networks trigger directed graphs
where direction matters. In this context, we assume a researcher has
a network adjacency matrix estimator at hand. For example, social
interaction models such as those described in \citet{depaula2023identifying,rothenhausler2015backshift,manresa}
treat adjacency matrix entries as estimands. Another variant is the
sampling of graphons that leads to noisy network matrices, developed
among others by \citet{10.1093/biomet/asac032,parise2023graphon}.
On the basis of such models, our results let us construct standard
errors for derived network statistics so that researchers can quantify
the propagated uncertainty.

In an application, we examine how confidence intervals for network
centralities arising in a simple network model provide a cautionary
tale about ranking nodes' popularity: reordering based on the upper
ends of the confidence intervals reorders the nodes' popularity in
one example but leaves the ordering undisturbed in another.

Besides, we also offer Monte Carlo evidence for the quality of the
distributional approximations. They perform well, but do depend on
the quality of the underlying matrix estimator. We also study the
performance of the $t$-test for a data-generating process that starts
with a random graph model, which experiences normally distributed
disturbances.

Beyond networks, there are many statistical applications that require
researchers to find eigenvectors of matrices estimated with error.
For example, companion matrices of vector auto-regressions are not
symmetric yet their spectrum carries information about the dynamics.
Eigenvectors associated with unit eigenvalues of these matrices identify
cointegrating relations, for which our inference methods are also
useful. Similarly, eigenvectors are used to estimate functional diversity
in ecology.

We calculate convergence rates of subspace-based estimators whenever
the convergence speed in the form of the Frobenius norm $\smlnorm{\hat{M}-M}_{\text{F}}$
is known. We also approximate higher-order bounds and learn that for
invariant subspaces, the eigenvalue separation dominates all higher-order
terms. To control convergence rates of eigenspaces, statisticians
only need control over the first-order bound as all higher-order bounds
are powers of the former. To first-order, we recover the version of
the \citet{davis1969some} bound found in \citet{demetrius}. These
results are helpful to strengthen consistency results to explicit
convergence rates.

Invariant subspaces of adjacency matrices also appear in latent space
graph models, where the latent space is either an invariant subspace
directly of the adjacency matrix or graph Laplacian or a hidden Euclidean
space that can be estimated via invariant subspaces as shown in \citet{ZhangXuZhu2022}.
Invariant subspaces are also used to approximate the latent spaces
in random dot product graph models \citep{10.1007/978-3-540-77004-6_11}
where \citet{XieXu2020} propose a method to estimate these spaces
and \citet{AthreyaEtAl2018} outline how to perform inference in random
dot product graphs.

An important distinction arises in the dimensionality of the $p\times p$
matrix $M$. In many scientific disciplines, network size $p$ is
modest and multiple observations of $M$ are available such as network
detection and modelling \citep{cattuto2010dynamics,krivitsky2014separable,prawesh2019small}
or health dynamics \citep{rothenberg1998social,cornwell2009network,christakis2010social},
or biology \citep{krause2009social,isella2011sociopatterns}. This
paper is predominantly about this case, which arises e.g. from the
measurement of many small networks, sampled over time. Importantly,
none of the Jacobians, the convergence results, and perturbation estimations
are sensitive to the size of $M$ whereas the hypothesis tests in
\prettyref{sec:hyp-tests} use fixed size CLTs and require estimation
of a $p^{2}\times p^{2}$ covariance matrix by standard methods.

Intuitively, symmetric matrices are convenient because small perturbations
to the entries correspond to small perturbations to the eigenvectors
and the eigenvector map is smooth. As one moves away from symmetric
matrices, eigenvalues acquire imaginary parts and come in complex
conjugate pairs. Therefore, eigenvalues may lie closely together though
differentiability still holds. Finally, if multiple eigenvalues correspond
to a single eigenvector, the eigenvector map is non-differentiable.
In this case, we can only distinguish groups of eigenvalues. Our methods
specialise in this latter case.

\section{Setup}

\label{sec:setup}

\subsection{General framework for subspace inference\label{subsec:General-framework}}

We define the framework for the invariant and singular subspace inference
problems. For a $p\times p$ matrix $M$, the columns of the $p\times q$
matrix $R$ span a right-invariant subspace of dimension $q$ iff
the relation

\begin{equation}
MR=R\Lambda\label{eq:inv-subs-rel}
\end{equation}
holds for some not necessarily full-rank $q\times q$ matrix $\Lambda$
with the analogous relationship $L^{\trans}M=\Lambda L^{\trans}$
and $R^{-1}=L^{\trans}$. In this case, the columns of $R$ span
a right-invariant subspace of $M$. For inference on eigenvectors,
a requirement by \citet{Tyler1981} is that there exists a positive
definite symmetric matrix $\Gamma$ such that $\Gamma M$ is symmetric,
which we relax. We illustrate this condition in
\begin{example}
\label{exa:symmetrisable}The matrix $M_{1}$ is symmetric in the
metric of $\Gamma$ for
\begin{align*}
\Gamma & =\begin{bmatrix}2 & 1\\
1 & 6
\end{bmatrix} &  & M_{1}=\begin{bmatrix}1 & 3\\
1 & 1
\end{bmatrix}
\end{align*}
whereas for parameters $\lambda,a$
\begin{align*}
M_{2}=\begin{bmatrix}\lambda & a\\
0 & \lambda
\end{bmatrix}
\end{align*}
is not symmetric in the metric of $\Gamma$. A $2\times2$ matrix
$\Gamma$ is not guaranteed to be positive-definite for $\Gamma M_{2}$
to be symmetric and $M_{2}$ is non-diagonalizable and only has a
single eigenvector $\left[\begin{smallmatrix}1 & 0\end{smallmatrix}\right]^{\trans}$.
Its single eigenvector however spans an invariant subspace. Therefore,
Tyler's procedure applies to $M_{1}$ but not to $M_{2}$.
\end{example}
Essentially, the requirement we relax stipulates that a matrix $M$
be symmetric with respect to some positive definite inner product
which implies diagonalizability and real eigenvalues.

We assume that we observe a sample of $n$ matrices $M_{1},\dots,M_{n}$,
whence we estimate the mean $\hat{M}_{n}$ so that the columns of
$\sqrt{n}(\hat{M}_{n}-M)$ converge weakly to a multivariate normal
distribution centered at zero. Denote the span of the columns of
interest in $R$ by $\risof$ with associated set of eigenvalues $\eigs_{I}$.
Then, our null hypothesis is
\begin{equation}
H_{0}:\upsilon\in\risof\label{eq:basic-hypothesis-random-matrix}
\end{equation}
against the one-sided alternative $H_{1}:\upsilon\notin\mathcal{S}_{I}(M)$.
The elements of $\eigs_{J}$ index all other directions of interest
so that if $M$ is diagonalizable, eigenvalues $\eigs$ equal $\eigs_{I}\union\eigs_{J}$.
To attain a statistic equalling zero under $H_{0}$, we focus on the
orthocomplement $\upsilon_{\perp}$, which is an element of the left-invariant
subspace $\rjsof$ formed by the columns of $(R)^{-\trans}$ so that
\begin{equation}
H_{0}^{\ast}\,:\,\upsilon_{\perp}\in\rjsof,\label{eq:perp-inference}
\end{equation}
which is equivalent to \prettyref{eq:basic-hypothesis-random-matrix}.
\prettyref{app:projections-and-gen-inverses} provides details. Therefore,
a researcher specifies up to $m\leq q$ columns of $\upsilon$ leading
to an orthocomplement with $p-m$ columns. It is, of course, possible
to only select a subset of the orthocomplement if $p$ is large as
long as the hypothesised vector of interest is orthogonal to $\upsilon_{\perp}$.

We also consider an extension to inference on the left-singular\footnote{An eigendecomposition of $M=R\Lambda^{\trans}$ implies that $R$
contains right eigenvectors, while for a singular value decomposition,
$M=U\Sigma V^{\trans}$, $U$ contains left singular vectors. Our
study focusses on those matrices appearing `on the left.'} subspace of $M$ spanned by the columns of $U$ in 
\[
M=U\Sigma V^{\trans}.
\]
Suppose a researcher observes a sequence of matrices $\{M_{t}\}_{t=1}^{n}$
that constitute noisy measurements from an underlying model 
\[
M_{t}=M+\varepsilon_{t},
\]
for an idiosyncratic error $\varepsilon_{t}$. The estimator $\hat{M}_{T}\defeq T^{-1}\sum_{t=1}^{T}M_{t}$
is consistent.

The hypothesis tests covered in the next section focus mainly on the
case where the size of $M$ is fixed although we make no formal restrictions
on it. For some of the network applications in \prettyref{sec:centralities},
the observations are rows and columns of $M$ which consequently grows
in size. To make our setup formal, we have
\begin{assumption}
\label{assu:general-ass}~Let $M\in\reals^{p\times p}$ be a general
matrix and let $\Omega\in\reals^{p^{2}\times p^{2}}$ denote a positive-definite
covariance matrix. Then,
\end{assumption}
\begin{enumerate}
\item \label{enu:tightness}If $p=n$, is growing, the estimator $\hat{M}$
is $1/r_{n}$ consistent (tight), i.e. $\text{\ensuremath{\smlnorm{\hat{M}_{n}-M_{n}}}}_{\text{F}}=O_{p}(r_{n})$
for the Frobenius norm.
\item \label{enu:asymptotically-normal-estimator}The model $M_{t}=M+\varepsilon_{t}$
generates data $\left\{ M_{t}\right\} _{t=1}^{T}$ where $\varepsilon_{t}$
are identical with general covariance matrix $\Omega=\expect\vek\varepsilon_{i}\,\vek\varepsilon_{j}^{\trans}$
for all $i,j=1,\dots,T$.
\item \label{enu:cov-mat-estimation}We can consistently estimate the covariance
matrix by some covariance estimator $\hat{\Omega}$ so that $\hat{\Omega}\inprob\Omega$.
\end{enumerate}
\prettyref{assu:general-ass}\prettyref{enu:tightness} implies
that $\hat{M}$ converges weakly and covers the setup where a network
researcher estimates network links based on a single sample of $M$.
Further, this assumption is the foundation for invariant subspaces
of large, random matrices. In those cases, special care has to be
taken. \citet[Example 2.3]{benaychgeorges2018lectureslocalsemicirclelaw}
show that eigenvectors of length $n$ tend to be scaled by $n^{-1/2}$.
Our focus in this study is simply how the convergence rate propagates
to estimators of invariant subspaces. \citet{cai2021network} provides
details.

\prettyref{assu:general-ass}\prettyref{enu:asymptotically-normal-estimator}
implies
\[
\sqrt{T}\,\vek\left(\hat{M}_{T}-M\right)\wkc N\left(0,\Omega\right),
\]
which admits a general covariance structure within $\Omega$ and imposes
a stationary error distribution over the index $t$ when $M$ has
fixed size $p\times p$. Finally, \prettyref{assu:general-ass}\prettyref{enu:cov-mat-estimation}
ensures that we can estimate the covariance matrix of the residuals.
The following assumption for invariant subspaces defines the set of
matrices we work with.
\begin{assumption}
\label{assu:invar-sub-ass}Let $\mathscr{M}\subset\reals^{p\times p}$
to be the set of matrices such that for every $M\in\mathscr{M}$,
\end{assumption}
\begin{enumerate}
\item \label{enu:semisimple}There exists an invariant subspace of interest
spanned by the columns of $R_{I}\in\reals^{p\times q}$ such that
$MR_{I}=R_{I}\Lambda_{I}$ for some matrix $\Lambda_{I}\in\mathbb{R}^{q\times q}$.
This requirement is equivalent to the existence of a non-singular
matrix $R\defeq\left[\begin{smallmatrix}R_{I} & R_{J}\end{smallmatrix}\right]\in\reals^{p\times p}$
with equivalent matrices $L=R^{-\trans}$ such that 
\begin{equation}
R^{-1}MR=\begin{bmatrix}\Lambda_{I} & \Lambda_{IJ}\\
0 & \Lambda_{J}
\end{bmatrix}=:\Lambda\label{eq:basic-rel}
\end{equation}
where $\Lambda_{I}\in\reals^{q\times q}$ and $\Lambda_{J}\in\reals^{r\times r}$
for $r=p-q$.
\item \label{enu:no-evs}The eigenvalues of $\Lambda_{I}$ and $\Lambda_{J}$
denoted by $\eigs_{I}$ and $\eigs_{J}$ obey $\eigs_{I}\intsect\eigs_{J}=\emptyset$
and $\lambda\in\eigs_{I}$ implies that $\lambda^{\ast}\in\eigs_{I}$,
i.e. $\eigs_{I}$ is closed under conjugation. Furthermore, $\abs{\lambda_{q}}>\abs{\lambda_{q+1}}$.
\item \label{enu:rankRlu}$\rank G^{\trans}R_{I}=q$ for a full-rank normalizing
matrix $G\in\reals^{q\times p}$. 
\end{enumerate}
\prettyref{assu:invar-sub-ass}\prettyref{enu:semisimple} defines
general invariant subspaces. The most important special case of those
subspaces are eigenspaces which obtain if $\Lambda_{IJ}=0$ and $\Lambda_{I}$
and $\Lambda_{J}$ are diagonal. An intermediate case is that of a
partially diagonalizable matrix, which we could achieve by requiring
$\Lambda_{I}$ to be diagonal with distinct eigenvalues, $\Lambda_{IJ}=0$,
and $\Lambda_{J}$ to be in Jordan normal form with blocks of arbitrary
size. The appeal of generic invariant subspaces is that they always
exist even in the most adverse circumstances. We have defined $\Lambda$
matrices to be real because even if eigenvalues appear as possibly
defective complex conjugates of another, it is always possible to
define the Jordan real form. For a conjugate pair, we can then find
an invariant subspace of twice the dimension of the multiplicity of
the eigenvalue associated with it.

Finally, \prettyref{assu:invar-sub-ass}\prettyref{enu:no-evs} ensures
that the map from $M$ to its invariant subspaces spanned by the columns
of $M$ is differentiable and that we can discriminate between vectors
of interest in sets $I$ and those in $J$. In particular, we do
not assume that we can discriminate among eigenvalues in $\eigs_{I}$
or $\eigs_{J}$.

An example of $\hat{M}$ that is consistent and has asymptotically
normal columns is the social interactions model \citep{depaula2023identifying,manresa}
for estimands $m_{ij}$ and $\gamma$, where the outcome
\[
y_{it}=\gamma\sum_{j\neq i}m_{ij}x_{jt}+\varepsilon_{it},
\]
for unit $i=1,\dots,p$ at time $t=1,\dots,T$ depends on the values
of individual-specific covariates $x_{it}$. Collecting estimands
$\hat{m}_{ij}$ into a matrix results in an estimated network adjacency
matrix whence all derived statistics such as centrality scores inherit
the uncertainty.\footnote{We consider statistics that can be derived from the spectrum of the
adjacency matrix but our methods extend straightforwardly to those
originating with graph Laplacians, too.} Furthermore, the estimation techniques in \citet{vliet} and \citet[Section 4.2]{rothenhausler2015backshift}
correspond to OLS estimation of $M$. Generally, OLS methods work
well if $T>>p$, which we explore in \prettyref{subsec:application}.

\subsection{Singular subspaces}

For singular subspaces, our setup changes ever so slightly. If $M$
is instead an $m\times l$ general, wide matrix so that $m<l$ and
of rank $q\leq\min\left(m,l\right)$, with $\Sigma_{I}\in\reals^{q\times q}$.
The columns of $V_{I}\in\reals^{l\times q}$ and $U_{I}\in\reals^{m\times q}$
span left- and right singular subspaces of dimension $q$ iff
\begin{equation}
M=U_{I}\Sigma_{I}V_{I}^{\trans}+U_{J}\Sigma_{J}V_{J}^{\trans}\label{eq:svd}
\end{equation}
holds for $U_{J}\in\reals^{m\times\left(l-q\right)}$, $V_{J}\in\reals^{\left(l-q\right)\times l}$,
and $\Sigma_{J}\in\reals^{\left(l-q\right)\times\left(l-q\right)}$.
It is also instructive to consider that we could implement inference
on $U_{I}$ via the procedure in \citet{Tyler1981} applied to $MM^{\trans}$.
While using eigenvectors of $MM^{\trans}$ may make no difference
asymptotically, it may be less efficient because $\var\vek MM^{\trans}\geq\var\vek M$.
In analogy to \prettyref{assu:invar-sub-ass}, we have for the singular
subspace decomposition
\begin{assumption}
\label{assu:svd-ass}Define $\mathscr{S}$ such that for every $M\in\mathscr{S}$,
$M$ is an $m\times l$ real matrix with $m\leq l$, and of rank $q\leq m$.
Then,
\begin{enumerate}
\item \label{enu:svd-split}The singular subspaces split according to $M=U_{I}\Sigma_{I}V_{I}^{\trans}+U_{J}\Sigma_{J}V_{J}^{\trans}$
where $\rank M=q$ implies that $\Sigma_{I}\in\reals^{q\times q}$
is a diagonal, square matrix and $\Sigma_{J}=0$.
\item For a singular subspace of dimension $q$, we have $U_{I}\in\reals^{m\times q}$
where $q\leq F$.
\end{enumerate}
\end{assumption}
\prettyref{assu:svd-ass} defines a general setting that also covers
graph models, where $M$ may be sparse so that $q<<\min\left(m,l\right)$.

\section{Hypothesis tests}

\label{sec:hyp-tests}

We present Wald and $t$-tests for inference on basis vectors of invariant
and singular subspaces. Proofs of results appear in \prettyref{app:results-hyp-tests}.

To test $H_{0}$ in \eqref{eq:basic-hypothesis-random-matrix}, we
construct the orthocomplement to the candidate vectors $\upsilon_{\perp}$
so that under the null hypothesis, $\upsilon_{\perp}^{\trans}R_{I}L_{I}^{\trans}=0$.
Let $m$ denote the column dimension of $\upsilon_{\perp}$ which
is chosen by the researcher based on hypothesised vectors. Normally,
for a hypothesised eigenspace of dimension $q$, we have $m=p-q$
columns in $\upsilon_{\perp}$ although we could add additional columns
so that $m$ could exceed $p-q$. 

Letting $\hat{\Omega}_{W}^{+}$ denote a generalized inverse of $\hat{\Omega}_{W}$,
defined in \prettyref{eq:gen-inv-eig} and \prettyref{eq:gen-inv-svd},
for $\hat{\Omega}_{W}:=\hat{B}\hat{\Omega}\hat{B}$, we have
\begin{equation}
\hat{W}_{n}\left(\upsilon_{\perp}\right)\defeq T\vek\left(\upsilon_{\perp}^{\trans}\hat{R}_{I}\hat{L}_{I}^{\trans}\right)^{\trans}\hat{\Omega}_{W}^{+}\vek\left(\upsilon_{\perp}^{\trans}\hat{R}_{I}\hat{L}_{I}^{\trans}\right)\label{eq:waldstat-1}
\end{equation}
and 
\begin{equation}
B\defeq\left(L_{I}\otimes\upsilon_{\perp}^{\trans}R_{J}\right)\left\{ \left(\Lambda_{I}^{\trans}\otimes I_{J}\right)-\left(I_{I}\otimes\Lambda_{J}\right)\right\} ^{-1}\left(R_{I}^{\trans}\otimes L_{J}^{\trans}\right).\label{eq:jacobian}
\end{equation}
into which we can insert sample analogues. The Wald test is based
on a first-order approximation of how the estimated subspace $R_{I}$
changes in response to perturbations in the matrix estimate $\hat{M}$,
with the Jacobian $B$ capturing this response. In our next result,
we focus on inference but \prettyref{sec:jacobians-ho-perturbations}
lays out Jacobians and further analytic results in more detail. For
eigenspaces, $\Lambda_{IJ}=0$ which does not affect the test statistic
but the speed with which test statistics converge as we see in \prettyref{subsec:higher-order-dk}.
Therefore, the approximations to test statistics are robust to $\Lambda_{IJ}\neq0$
and converge fastest for eigenspaces.
\begin{thm}
\label{thm:waldstatdistrib}Suppose \prettyref{assu:general-ass}
holds. Then, $\hat{W}_{n}\left(\upsilon_{\perp}\right)\wkc\chi_{qm}^{2}.$
\end{thm}
\noindent The proof appears in \prettyref{app:distributions}.

Analogously to \eqref{eq:basic-hypothesis-random-matrix}, we write
hypotheses for singular vectors as $H_{0}^{\ast}:\upsilon_{\perp}^{\trans}U_{I}=0$
for a specified $m\times(m-h)$ matrix $\upsilon_{\perp}$ corresponding
to some $\upsilon_{0}$. For the covariance, define
\begin{align}
\jacsvd= & \left(\Sigma_{I}^{-1\trans}V_{I}^{\trans}\otimes\schwein U_{J}U_{J}^{\trans}\right)\label{eq:svd-jacobian}\\
 & +\left(I_{F}\otimes\schwein U_{I}\right)\left[\left(\indic^{\trans}\otimes\vek D_{\text{d}}\right)\cdot\left(\left(\Sigma_{I}^{\trans}V_{I}^{\trans}\otimes U_{I}^{\trans}\right)+\left(U_{I}^{\trans}\otimes\Sigma_{I}V_{I}^{\trans}\right)\right)\right].\nonumber 
\end{align}
In analogy to inference on invariant subspaces above, $\upsilon_{\perp,0}^{\trans}\in\reals^{m\times\left(m-h\right)}$
where $h\leq F$ is the size of the hypothesized singular subspace.

We define the covariance matrix $\hat{\Omega}_{S}:=\widehat{\jacsvd}\hat{\Omega}\widehat{\jacsvd}^{\trans}$
and denote by $\hat{\Omega}_{S}^{+}$ a generalized inverse of $\hat{\Omega}_{S}$.
The statistic is then
\begin{equation}
\hat{W}_{\text{SVD},n}\left(\upsilon_{\perp}\right)\defeq T\vek\left(\upsilon_{\perp}^{\trans}\hat{U}_{I}\right)^{\trans}\hat{\Omega}_{S}^{+}\vek\left(\upsilon_{\perp}^{\trans}\hat{U}_{I}\right)\label{eq:waldstat-svd}
\end{equation}
for which we obtain
\begin{thm}
\label{thm:svd-wald-test}Suppose \prettyref{assu:svd-ass} holds,
then $\hat{W}_{\text{SVD},n}\left(\upsilon_{\perp}\right)\wkc\chi_{q\left(m-h\right)}^{2}$.
\end{thm}
\prettyref{thm:svd-wald-test} is directly analogous to \prettyref{thm:waldstatdistrib}
and the proof is identical and available in \prettyref{app:results-hyp-tests}.

\noindent Although straightforward to compute, a Wald test statistic
is only required when a researcher specifies a full vector or matrix
hypothesis. We imagine that this scenario is most likely of interest
when testing whether localized basis vectors belong to an invariant
or singular subspace of $M$. Empirical applications may require inference
on individual coefficients of invariant or singular vectors, in which
case a scalar $t$-test is useful.

\noindent To make invariant subspace coordinates unique, ensuring
any associated estimates are consistent across experiments, we normalise
$\upsilon_{\perp}=:\begin{smallmatrix}[\upsilon_{1} & \upsilon_{2}]\end{smallmatrix}$
where $\upsilon_{1}\in\reals^{q\times q}$ and $\upsilon_{2}\in\reals^{r\times q}$.
Then letting $D^{\trans}\defeq\upsilon_{2}\upsilon_{1}^{-1}\in\reals^{r\times q}$,
we have

\begin{equation}
\upsilon_{\perp}\defeq\begin{bmatrix}I_{q}\\
-D^{\trans}
\end{bmatrix}\label{eq:eigenvector-normalization}
\end{equation}

\noindent If an eigenvector represents a node's centrality score,
this normalization corresponds to choosing a numeraire node with unit
centrality so that coordinates are comparable across samples. However,
this normalization may cause outliers making it potentially unstable.
In our simulations, we did not find any evidence of any problems.

Analogous to \prettyref{eq:eigenvector-normalization}, we can derive
a closed-form expression of $\hat{D}_{n}$ appearing in $R_{I}$.
Let $R_{I,1}\in\reals^{r\times q}$ and $R_{I,2}\in\reals^{q\times q}$
with $\rank R_{I,2}=q$ so that
\begin{equation}
\begin{bmatrix}R_{I,1}\\
R_{I,2}
\end{bmatrix}\defeq R_{I}.\label{eq:partitionforahat-1}
\end{equation}
Consequently, 
\begin{equation}
\hat{D}_{I,n}^{\trans}\defeq\hat{R}_{I,1,n}\hat{R}_{I,2,n}^{-1}.\label{eq:ahat-1}
\end{equation}
so the normalized vector is $\begin{smallmatrix}[-D & I_{q}]^{\trans}\end{smallmatrix}$.
Observe that \prettyref{eq:ahat-1} defines a unique estimator.

\noindent For inference on individual entries of $D$, define $d_{ij}\defeq e_{i}^{\trans}De_{j}$
for $e_{i}\in\reals^{q\times1}$ and $e_{j}\in\reals^{r\times1}$
where the vectors $e_{i}$ have unit entries at $i$ and zero elsewhere.
Let $\hat{R}_{I,2}^{-1}$ be the empirical analogue of $R_{I,2}^{-1}$
defined in \prettyref{eq:partitionforahat-1} and define
\begin{equation}
B_{ij}=\left(e_{i}^{\trans}R_{2,I}^{\trans-1}\otimes e_{j}^{\trans}\upsilon_{\perp}^{\trans}R_{J}\right)\left[\left(\Lambda_{I}^{\trans}\otimes I_{J}\right)-\left(I_{I}\otimes\Lambda_{J}\right)\right]^{-1}\left(R_{I}^{\trans}\otimes L_{J}^{\trans}\right),\label{eq:jac-coeff}
\end{equation}
in analogy to \eqref{eq:jacobian}. We write the scalar variance
of $\sqrt{T}(\hat{d}_{ij}-d)$ as $\sigma_{ij}^{2}\defeq B_{ij}\Omega B_{ij}^{\trans}$
and define the $t$-test statistic for inference on the $ij$th coefficient
of $D$
\begin{equation}
t_{ij,n}\left(d_{0}\right)\defeq\frac{\hat{d}_{n,ij}-d_{0}}{\sqrt{\hat{\sigma}_{ij}^{2}/T}}.\label{eq:tstat-1-1}
\end{equation}
We construct the estimator of the variance, $\hat{\sigma}_{ij}$ by
replacing $B_{ij}$ with $\hat{B}_{ij}$ which contains sample analogues
$\hat{R}_{i}$ and $\hat{L}_{i}$ for $i\in\left\{ I,J\right\} $.
The distribution is given in
\begin{thm}
\label{thm:tdist}Suppose \prettyref{assu:general-ass} holds, then
$t_{ij,n}\left(d_{0}\right)\wkc N\left(0,1\right).$
\end{thm}
\noindent Eigenvector-based centralities are based on the absolute
values $|\hat{d}_{ij}|$. We define for this result the folded normal
distribution with cumulative distribution function 
\begin{equation}
F_{G}\left(x;d_{ij},\sigma_{ij}\right)\defeq\Phi\left(\frac{x-d_{ij}}{\sigma_{ij}}\right)-\Phi\left(\frac{-x-d_{ij}}{\sigma_{ij}}\right)\label{eq:cdf-folded}
\end{equation}
where $\Phi\left(\frac{x-\mu}{\sigma}\right)$ denotes the normal
cumulative distribution function with mean $\mu$ and variance $\sigma^{2}$.
For the asymptotic distribution of the absolute values of the normalized
eigenvector entries, we obtain 
\begin{cor}[Folded normal distribution]
\label{cor:folded-normal}Suppose \prettyref{assu:general-ass} holds,
then $\sqrt{n}|\hat{d}_{ij}|\wkc G$ where $G$ is a random variable
that has a folded normal distribution with c.d.f. \textup{given in
\eqref{eq:cdf-folded}.}
\end{cor}
\prettyref{subsec:Network-centralities-as} gives an application of
this result.

\subsection{Estimators for basis vectors\label{sec:estimators}}

Our next result characterizes the distribution of $\hat{D}_{I,n}^{\trans}$.
Define 
\begin{equation}
B_{D^{\trans}}\defeq\left(R_{I,2}^{\trans-1}\otimes\upsilon_{\perp}^{\trans}R_{J}\right)\left\{ \left(\Lambda_{I}^{\trans}\otimes I_{J}\right)-\left(I_{I}\otimes\Lambda_{J}\right)\right\} ^{-1}\text{\ensuremath{\left(R_{I}^{\trans}\otimes L_{J}^{\trans}\right)}}.\label{eq:d-jac}
\end{equation}
Then, we obtain
\begin{thm}[Distribution of basis vectors]
\label{thm:estdist}Let \prettyref{eq:ahat-1} define $\hat{D}_{n}$.
Then, \prettyref{assu:general-ass} implies $n^{1/2}\vek\{\hat{D}_{n}^{\trans}-D^{\trans}\}\wkc N(0,B_{D^{\trans}}\Omega B_{D^{\trans}}^{\trans}).$
\end{thm}
\noindent The above results rely on consistent covariance matrix estimators
which we discuss alongside other proofs in \prettyref{app:distributions}.

To make precise statements about smoothness of invariant and singular
subspace maps, we define $\psi\left(M;\upsilon\right)$ that represents
invariant subspace statistics and has informative null distributions
on spaces spanned by a subset of the columns of $R$ or $U$. Let
$\psi\,:\,\set M\rightarrow\reals^{r\times q}$ for $M\in\mathscr{M}$
in \prettyref{assu:invar-sub-ass} by
\begin{equation}
\psi\left(M,\upsilon_{\perp}^{\trans}\right)\coloneqq\upsilon_{\perp}^{\trans}R_{I}L_{I}^{\trans}\left(M\right)\label{eq:inv-subspace-map}
\end{equation}
and analogously $\Psi(M,\upsilon_{\perp}^{\trans})\defeq\upsilon_{\perp}^{\trans}U_{I}$. 

We define the map $\Psi:\mathscr{S\rightarrow\reals}^{r\times q}$
via 

\begin{equation}
\Psi\left(M,\upsilon_{\perp}^{\trans}\right)\defeq\schwein U_{I}\left(M\right),\label{eq:sin-subspace-map}
\end{equation}
so that $\Psi\left(M,\upsilon_{0,\perp}^{\trans}\right)=0$ under
the null hypothesis.\footnote{Without loss of generality, we restrict ourselves to left-singular
vectors as inference on $\col V$ follows from $M^{\trans}$ instead
of $M$.}

The zero-level set $\psi(M,\upsilon_{0,\perp}^{\trans})=0$ defines
non-rejection regions for $\upsilon_{\perp}$ and $\upsilon$. We
have set estimators

\begin{equation}
\widehat{\upsilon}_{\perp,n}\defeq\left\{ \upsilon_{\perp}\in\reals^{p\times m}\,:\,\hat{\Phi}\left(\schwein\right)=0\right\} \label{eq:betahat-1}
\end{equation}
where 
\[
\hat{\Phi}\left(\schwein\right)\defeq\begin{cases}
\psi\left(\hat{M}_{n},\schwein\right) & \text{invariant subspace,}\\
\Psi\left(\hat{M}_{n},\schwein\right) & \text{singular subspace.}
\end{cases}
\]
Estimators of $\psi\left(M\right)$ obtain from evaluating $\psi$
at the estimate $\hat{M}$, so that 
\[
\psi\left(\hat{M}_{n}\right)=\upsilon_{\perp}^{\trans}\hat{P}_{I}\left(\hat{M}_{n}\right).
\]
Using \prettyref{lem:jacobian}\prettyref{enu:invariant}, we obtain
the distribution of the sample analogue of $\psi(M)$ in
\begin{lem}
\label{lem:proj-distrib}Suppose \prettyref{assu:general-ass} holds
for invariant and additionally \prettyref{assu:svd-ass} holds for
singular subspaces. Let $B$ be as in \prettyref{eq:jacobian} or
\prettyref{eq:svd-jacobian} and let $\upsilon_{\perp}$ be as in
\prettyref{eq:per-def}. Then,\textup{
\[
\sqrt{T}\,\vek\hat{\Phi}\left(\schwein\right)\wkc N\left(0,B\Omega B^{\trans}\right).
\]
}
\end{lem}
For general eigenvectors, i.e. if the set of eigenvalues $\eigs_{I}$
is not closed und\ensuremath{\le}er conjugation, $\sqrt{n}\,\vek\psi(\hat{M}_{n})$
converges to a multivariate complex normal distribution, with covariance
$B\Omega B^{\trans}$ and relation $B\Omega B^{\prime}$. For now,
we shall continue to assume that $\eigs_{I}$ is closed under conjugation
or, in the case of $\abs{\eigs_{I}}=1$, that the Perron-Frobenius
theorem applies although our results also accommodate these cases.

\section{Jacobians and smoothness\label{sec:jacobians-ho-perturbations}}

\subsection{\label{subsec:expansion-psi}Invariant subspace map $\psi$}

In the background, the hypothesis tests in \prettyref{sec:hyp-tests}
relied on a first-order Taylor expansion argument to find standard
errors. To make our arguments about higher-order perturbations to
invariant subspaces precise, we use the map $\psi$ that assigns such
spaces to non-symmetric matrices. \prettyref{assu:invar-sub-ass}\prettyref{enu:no-evs}
ensures that eigenvalues are semi-simple within their respective groups,
which allows differentiability of the invariant subspaces spanned
by columns of $R$ belonging to $\eigs_{I}$. This setup allows eigenvalues
to coincide within $\eigs_{I}$ so that our setup covers matrices
whose Jordan blocks have non-unit size, i.e. are not simple. Generally,
such eigenvalues imply that the corresponding eigenvectors are not
differentiable because small perturbations in the underlying matrix
have unpredictable consequences for the corresponding eigenvectors.
Our baseline setup covers exactly this case while some results in
\prettyref{sec:centralities} reduce to the case of simple eigenvalues.

In this section, we establish that $\psi$ and its singular subspace
companion $\Psi$ are infinitely differentiable: as their arguments
change slightly, the maps respond in a controlled and predictable
way enabling expansions of arbitrary order. Smoothness is therefore
not just a technical regularity: infinite differentiability enables
refined inference, including higher-order likelihood and Edgeworth
expansions. Our main novelty lies in achieving smoothness and deriving
Jacobians for generic invariant subspaces without assuming symmetry
or diagonalisability, extending the classical results. Unsurprisingly,
smoothness applies to singular subspace maps, too. Furthermore, we
characterise higher-order bounds of invariant subspaces. For estimations
of bounds of singular subspace maps beyond second order, we recommend
adapting invariant subspace expressions\footnote{Invariant subspaces of $MM^{\trans}$ are singular subspaces of $M$.}
because of the tractability of the perturbation theory for invariant
subspaces. 

To summarize our results succinctly, we define the resolvent map $\mathscr{P}\rightarrow\reals^{p\times p}$,
for
\begin{equation}
X\mapsto\left(X-\alpha I\right)^{-1}\label{eq:res}
\end{equation}
which forms the basis of all spectral-based statistics. The main result
is
\begin{thm}
\label{thm:resolvent}The resolvent from $\mathscr{P}\rightarrow\reals^{p\times p}$
in \prettyref{eq:res} defines a smooth map.
\end{thm}
This result guarantees that network statistics based on \prettyref{eq:res},
or, equivalently, those based on spectral data of $X$, have existent
higher-order expansions and hence inherit a certain stability.\footnote{If $X$ is less than full rank, then a generalized inverse can be
used without altering the result.} Furthermore, all statistics based on \eqref{eq:res} are consistently
estimable. Because singular vectors of $M$ are eigenvectors of $MM^{\trans}$,
\prettyref{thm:resolvent} applies to $\Psi$ as well. A proof of
\prettyref{thm:resolvent} appears in \prettyref{app:smooth-maps-perturb-expansions}.

Jacobians appear in
\begin{lem}[Jacobians]
\label{lem:jacobian}~
\end{lem}
\begin{enumerate}
\item \label{enu:invariant}Let \prettyref{eq:inv-subspace-map} define
$\psi$. Then, the Jacobian matrix with respect to $M$ in the direction
of $E\in\reals^{p\times p}$ of $\psi$ is $B\in\reals^{p^{2}\times rq}$
such that$\diff\vek\psi(M,E)=B\vek E$ for $B$ in \eqref{eq:jacobian}.
\item \label{enu:singular}Let \prettyref{eq:sin-subspace-map} define $\Psi$.
Then, the Jacobian matrix with respect to $M$ in the direction of
$E\in\reals^{m\times l}$ of $\Psi$ is $\jacsvd\in\reals^{F\left(m-h\right)\times ml}$
such that $\diff\vek\Psi\left(M;E\right)=\jacsvd\vek E$ for $B_{\text{SVD}}$
in \eqref{eq:svd-jacobian}.
\end{enumerate}
\noindent\prettyref{app:Additional-details-on-invar} contains a
proof. Importantly, \eqref{eq:jacobian} reduces to $C_{w}$ in \citet[4.3]{Tyler1981}
when $M$ is symmetrizable. In this case, the columns of $R_{I}$
and $L_{I}$ satisfy $r_{i}l_{i}^{\trans}\equiv r_{i}r_{i}^{\trans}\Gamma$
where $\Gamma$ is a real positive definite symmetric matrix such
that $M\Gamma$ is symmetric. If such $\Gamma$ is found, the distinction
between left and right eigenvectors disappears as $\Gamma r_{i}$
takes the place of $l_{i}$ for $i\in I$.

\subsection{Extension of Davis-Kahan bound to higher-order perturbations}

\label{subsec:higher-order-dk}

The Jacobian in \eqref{eq:jacobian} reduces to the bound in \citet{demetrius},
which asserts that small perturbations produce changes bounded by
the ratio of the perturbation and the size of the spectral gap. For
a $2\times2$ matrix, the first-order term
\begin{align}
\norm{\diff\psi\left(M,E\right)}_{\text{F}} & \leq\frac{\norm E_{\text{F}}}{\lambda_{1}-\lambda_{2}}\label{eq:davis-kahan-disguise}
\end{align}
where $\text{d}\psi$ is the infinitesimal response to perturbation
$E$. This bound differs from the original found in \citet{davis1969some}
by the fact that only population eigenvalues appear in the denominator.
The relation \eqref{eq:davis-kahan-disguise} inspires the question
whether we can use perturbative arguments to derive higher-order bounds.
The answer is affirmative and to simplify the exposition, we define
the size of the perturbation
\begin{align*}
a & \defeq\max\left\{ \norm{\hat{M}-M}_{\text{F}},\norm{L_{j}^{\trans}\left(M-\hat{M}\right)R_{i}}\right\} \\
 & =\norm{\hat{M}-M}_{\text{F}}
\end{align*}
for $i,j\in\left\{ I,J\right\} $ where the equality follows from
the proof of \prettyref{thm:stoch-bound-prop} and the discussion
in \prettyref{app:inequality-block-diagonal}. We define the eigenvalue
separation as 
\[
s\defeq\norm{\left\{ \left(\Lambda_{I}^{\trans}\otimes I_{J}\right)-\left(I_{I}\otimes\Lambda_{J}\right)\right\} ^{-1}}_{\text{F}}.
\]
Similarly, we have the distance 
\[
\alpha_{12}\defeq\norm{\left(L^{\trans}MR\right)_{1:p-q,q+1:p}}_{\text{F}}
\]
which measures how far $L^{\trans}MR$ is from being block-diagonal.

For eigenvectors, higher-order bounds are simply the first-order bound
raised to a higher power so that control over $as$ implies control
over all orders. For general invariant subspaces $\alpha_{12}\neq0$
which amplifies higher-order effects, revealing a qualitative difference
in robustness.

Denote by $h(s^{x},a^{y})$ a polynomial with leading powers $x$
and $y$ in $s$ and $a$, respectively. Then, the first-order bound
$\gamma_{1}$ is given by \eqref{eq:davis-kahan-disguise}. While
first-order error is known to scale with $a/s$, it it is unclear
how higher-order terms behave for eigenspaces or invariant subspaces.
For general invariant subspaces, $\alpha_{12}>0$ and the eigenvalue
separation $s$ has a progressively stronger effect on higher-order
error terms than the estimation precision. This relationship contrasts
sharply with the eigenspace case where $\alpha_{12}=0$, where higher-order
bounds are simply powers of the first-order term.
\begin{thm}[Higher-order Davis-Kahan bounds]
\label{thm:higher-order-dk}The first order bound in \eqref{eq:davis-kahan-disguise}
extends to higher orders according to
\begin{align*}
\gamma_{2} & \approx h\left(s^{2},a^{2}\right)+\indic\left\{ \alpha_{12}>0\right\} \alpha_{12}s\gamma_{2}\approx h\left(s^{2},a^{2}\right)+\alpha_{12}\indic\left\{ \alpha_{12}>0\right\} h\left(s^{3},a^{2}\right)\\
\gamma_{3} & \approx h\left(s^{3},a^{3}\right)+\indic\left\{ \alpha_{12}>0\right\} \alpha_{12}^{2}\gamma_{3}\approx h\left(s^{3},a^{3}\right)+\alpha_{12}^{2}\indic\left\{ \alpha_{12}>0\right\} s^{2}h\left(s^{3},a^{3}\right)\\
\gamma_{4} & \approx h\left(s^{4},a^{4}\right)+\indic\left\{ \alpha_{12}>0\right\} \alpha_{12}^{3}h\left(s^{7},a^{4}\right)\\
 & \vdots\\
\gamma_{k} & \approx h\left(s^{k},a^{k}\right)+\indic\left\{ \alpha_{12}>0\right\} \alpha_{12}^{k-1}h\left(s^{2k-1},a^{k}\right).
\end{align*}
For eigenspaces $\alpha_{12}=0$, in which case the higher order bounds
are homogeneous polynomials in $as$.
\end{thm}
The approximations to higher-order terms only list the highest power
of $\alpha_{12}$ but a term linear in $\alpha_{12}$ nevertheless
appears and is a polynomial with leading powers $k$ in $s$ and $a$,
($h(s^{k},a^{k})$). The first two equalities follow from $\gamma_{2}=O(s^{2}a^{2})$
multiplying $s$and $\gamma_{3}=O(s^{3}a^{3})$ multiplying $s^{2}$.
The lesson of \prettyref{thm:higher-order-dk} is that invariant subspace
estimation converges slowly if $\alpha_{12}>1$. Moreover, a quickly
converging estimator with $a\goesto0$ cannot easily balance out a
larger reciprocal eigenvalue separation $s$ because its higher-order
contributions rise faster in $s$ than they do in $a$. Specifically,
its powers dominate with $2k-1$ relative to $k$ for $a$ so that
statisticians estimating invariant subspaces are advised to check
$\alpha_{12}$ whenever Davis-Kahan arguments are applied. A proof
is given in \prettyref{app:perturb-to-deriv}.

\section{Network statistics\label{sec:centralities}}

In this section, we study the scenario where a researcher has an estimator
$\hat{M}$ of $M$ that fulfills $\smlnorm{\hat{M}-M}_{\text{F}}=O_{p}\left(r_{n}\right)$
for some sequence $r_{n}$. Using $r_{n}$ as an input, we find the
stochastic order of $\smlnorm{\hat{\psi}-\psi}$ for subspace-based
network statistics $\psi$ as well as the node-wise and overall clustering
coefficients. Moreover, we apply the inference methods of \prettyref{sec:hyp-tests}
to the problem of finding standard errors for network centrality measures,
all of which are based on invariant subspace decompositions of the
network adjacency matrix. Proofs appear in \prettyref{app:results-network-centralities}.

\subsection{\label{subsec:Network-centralities-as}Convergence rates of network
statistics}

For our results to apply, the leading eigenvalue of the adjacency
matrix has to be distinct from all the others. In this case, \prettyref{assu:invar-sub-ass}\prettyref{enu:no-evs}
becomes the singleton $\eigs_{I}=\{\lambda_{\text{max}}\}$. Both
this assumption as well as the Perron-Frobenius theorem in the case
of adjacency matrices ensure that this condition is met. The application
of the Perron-Frobenius theorem here is crucial because it guarantees
that the top eigenvalue is simple so that the single top eigenvector
map is differentiable.

\subsection*{Convergence rate of $\psi$ as a network statistic for matrices with
potentially growing size}

The map $X\mapsto(X-\alpha I)^{-1}$ forms the basis of a large
class of network centrality measures, which we term \emph{resolvent-based
centrality scores }and in \prettyref{thm:resolvent}, we establish
that it is infinitely differentiable. The immediate consequence is
that all network statistics admit expansions to arbitrary order. For
$g\in\{n,T\}$, recall \prettyref{assu:general-ass}\prettyref{enu:tightness},
which states that
\begin{equation}
\text{\ensuremath{\norm{\hat{M}-M}}}_{\text{F}}=O_{p}\left(r_{g}\right)\label{eq:stoch-bd}
\end{equation}
and suppose that a researcher has shown that \prettyref{eq:stoch-bd}
holds for some estimator. The question that naturally arises is to
what extent 
\[
\norm{\psi\left(\hat{M}\right)-\psi\left(M\right)}_{\text{F}}
\]
is then bounded in probability. Naturally, $\psi(\hat{M})\inprob\psi(M)$
as $r_{g}\goesto0$ follows from \prettyref{thm:resolvent} as $\alpha\goesto\lambda_{\text{max}}$.\footnote{As $\alpha\goesto\lambda_{\text{max}}$, the resolvent approaches
the top eigenvector. See the proof of \prettyref{lem:equivalence-centr}
in Appendix \ref{subsec:network-centrality-proofs} for an explanation.} But we wish to quantify the exact convergence rate. Denote by $\psi_{1}$
a statistic based on the top eigenvector of $M$. Then, we have
\begin{thm}
\label{thm:stoch-bound-prop}Suppose the stochastic bound \eqref{eq:stoch-bd}
(\prettyref{assu:general-ass}\prettyref{enu:tightness}) holds and
that the Perron-Frobenius theorem applies to $M$. Then, the following
apply:
\begin{enumerate}
\item Statistics based on the principal component of $M$, $\psi_{1}$,
obey 
\begin{equation}
\norm{\psi_{1}\left(\hat{M}\right)-\psi_{1}\left(M\right)}_{\text{F}}=O_{p}\left(\frac{r_{g}}{\lambda_{1}-\lambda_{2}}\right).\label{eq:principal-comp-bound}
\end{equation}
\item If the statistic encompasses more than one eigenvector (invariant
vector), let $1\leq r\leq s\leq\rank A$ and $I=\left\{ r,r+1,\dots,s\right\} $
with ordered eigenvalues $\lambda_{r}<\lambda_{r+1}<\dots<\lambda_{s}$
where $\spn\psi\left(\hat{M}\right)=\spn\left\{ R_{I}\right\} +o_{p}\left(1\right)$.
In this case, the bound changes to 
\[
O_{p}\left(\frac{r_{g}}{\min\left\{ \left(\lambda_{r-1}-\lambda_{r}\right),\left(\lambda_{s}-\lambda_{s+1}\right)\right\} }\right).
\]
\end{enumerate}
\end{thm}

\noindent The bound attained in \prettyref{thm:stoch-bound-prop}
resembles the Davis-Kahan theory \citep{davis1969some,demetrius}
except that we derived it from perturbative arguments.
\begin{rem}
The result in \prettyref{thm:stoch-bound-prop} is valid for $M$
with fixed size $p\times p$ or $n\times n$ where $n$ is understood
to be possibly divergent.
\end{rem}
\noindent For eigenvector-based statistics, we can set $\alpha_{12}=0$
in \prettyref{thm:higher-order-dk} and combine the result with \prettyref{thm:stoch-bound-prop}
to a bound to $k$th order which reads 
\begin{equation}
\sum_{i=1}^{k}O_{p}\left(\frac{r_{g}}{\min\left\{ \left(\lambda_{r-1}-\lambda_{r}\right),\left(\lambda_{s}-\lambda_{s+1}\right)\right\} }\right)^{i}.\label{eq:higher-order-bound}
\end{equation}
An important case arises if the Perron-Frobenius theorem does not
apply or if eigenvalues are not simple, so that individual eigenvectors
are no longer differentiable. In this case, we have to focus the analysis
on generic invariant subspaces, where $\alpha_{12}\neq0$ in \prettyref{thm:higher-order-dk}.
It is possible to extend \prettyref{thm:stoch-bound-prop} to these
cases.

\subsection*{Clustering coefficient}

Follow \citet{jackson2008social} and define the node-wise clustering
coefficient $\text{cl}_{i}$ via 
\begin{equation}
\text{cl}_{i}\defeq\frac{\sum_{j\neq i,k\neq j,k\neq i}m_{ij}m_{ik}m_{jk}}{\sum_{j\neq i,k\neq j,k\neq i}m_{ij}m_{ik}},\label{eq:clustering-coeff}
\end{equation}
where the sums run over all indices but $i$. If we are dealing with
the overall clustering coefficient we have instead
\[
\text{cl}\defeq\frac{\sum_{i,j\neq i,k\neq j,k\neq i}m_{ij}m_{ik}m_{jk}}{\sum_{i,j\neq i,k\neq j,k\neq i}m_{ij}m_{ik}}.
\]
Estimators obtain from replacing $m_{ij}$ with elements of $\hat{M}$. 
\begin{thm}[Convergence rate of clustering coefficients]
\label{thm:convergence-rate-clustering-coeff}Suppose the stochastic
bound \eqref{eq:stoch-bd} (\prettyref{assu:general-ass}\prettyref{enu:tightness})
holds. For the node-wise and individual clustering coefficients, $\smlabs{\hat{\text{\text{cl}}}_{i}-\text{cl}_{i}}=O_{p}\left(r_{n}\right)$
and $\smlabs{\hat{\text{\text{cl}}}-\text{cl}}=O_{p}\left(r_{n}\right).$
\end{thm}
\prettyref{thm:convergence-rate-clustering-coeff} shows that the
node-wise and overall clustering coefficients converge at the same
rate as the adjacency matrix estimates.

\subsection{Network centrality measures}

In the following, we shall use standard errors for centrality scores
to construct $t$-tests. Note that these all centrality scores satisfy
the convergence rate given in \prettyref{thm:stoch-bound-prop}. Recall
that standardisation by $g\in\{n,T\}$ depends on whether $M$ has
growing or fixed size where \citet[Example 2.3]{benaychgeorges2018lectureslocalsemicirclelaw}
provide details on required normalization.

\subsection*{Eigenvector centrality}

This centrality score $c_{i}$ assigns popularity to node $i$ as
the sum of the popularity of its neighbors. Intuitively speaking,
if someone is connected to very popular nodes, one is themselves very
popular so that $c_{i}\coloneqq\frac{1}{\lambda}\sum_{j\in N_{i}}c_{j}.$
The sum runs over all $j$ that are connected to node $i$, denoted
by $N_{i}$ and $\lambda$ is a normalizing constant. For the full
vector of scores $c$, we make the substitution $\sum_{j\in N_{i}}c_{j}=\sum_{j=1}^{N}m_{ij}c_{j}$
and obtain the eigenvector problem $Mc=\lambda c$. Because $c$ is
a basis vector belonging to the eigenspace of $M$ associated with
the largest eigenvalue of $M$, we can apply the results of \prettyref{sec:estimators}
directly to this network statistic. Importantly, this measure is not
invariant to normalizations. We thus suggest applying the normalization
in \prettyref{eq:eigenvector-normalization} which defines centralities
uniquely by denoting one node as the reference and use the test in
\prettyref{thm:tdist} for inference on individual coefficients. Formally,
we have
\begin{thm}
\label{thm:eigenvec-smooth}The eigenvector centrality $c\left(M\right)$
 is a smooth function of $M$. Therefore, $\hat{M}_{n}\inprob M$
implies $\hat{c}\inprob c$ and for an adjacency matrix with distribution
$N\left(\vek M,\Omega\right)$, we have 
\[
\frac{\hat{c}_{i}-c_{i}}{\sqrt{B_{1j}\Omega B_{1j}^{\trans}/g}}\wkc N\left(0,1\right)
\]
for the Jacobian of $B_{1j}$ defined in \eqref{eq:jac-coeff}.
\end{thm}
The above result is a direct consequence of \prettyref{thm:estdist}
and allows one to construct standard errors. In applied work, we wish
to test hypotheses of the form $H_{0}:c_{i}=c_{0}$ against the alternative
that $c_{i}\neq c_{0}$. We use the absolute value version of the
coefficient-wise $t$-test based on \prettyref{cor:folded-normal}
to construct one-sided confidence intervals and define the one-sided
confidence interval with level $\alpha$ for $\abs{c_{i}}$\footnote{This centrality measure is usually implemented as the abs. value of
the coefficients, see e.g. \href{https://github.com/JuliaGraphs/Graphs.jl/blob/master/src/centrality/eigenvector.jl}{https://github.com/JuliaGraphs/Graphs.jl/blob/master/src/centrality/eigenvector.jl.}.} via$\mathcal{C}_{i}\defeq\left[\abs{\hat{c}_{i}},x_{0}\right],$where
the upper limit $x_{0}>0$ solves $F_{G}\left(x_{0};\abs{c_{i}},\sigma_{i}\right)=1-\alpha$
for $x_{0}>0$. Expression \prettyref{eq:cdf-folded} defines $F_{G}$
with consistent point estimates $\abs{\hat{c}_{i}}$ and $\hat{\sigma}_{i}$
substituted for $\abs{s_{i}}$ and $\sigma_{i}$ for $i=1,\dots,p-1$.
Consequently, $\forall c\in\mathcal{C}$, we cannot reject $H_{0}\,:\,\abs{c_{i}}=c$.

\subsection*{PageRank, Katz, and Diffusion centralities}

A related family of centrality measures is conceptually very similar
to the eigenvector centrality but performs better for dense graphs.
Formally, the PageRank centrality \citep{brin1998anatomy} $c_{P}$
is given by 
\begin{align}
c_{P} & =(I-\alpha M)^{-1}\beta\label{eq:pagerank-def}
\end{align}
where $\alpha\in\reals$ and $\beta\in\reals^{p}$ are constants.
The damping factor $\alpha$ should be chosen to lie between $0$
and $1/\lambda_{\text{max}}\left(M\right)$. We can easily derive
\prettyref{eq:pagerank-def} from the eigenvector centrality by adding
a small amount of ``centrality'' to each node in the eigenvector
relation in the form of $\beta$ to obtain $c=\alpha Mc+\beta$ for
some choice of $\alpha$ and solving for $c$. To motivate \prettyref{eq:pagerank-def}
further, we imagine a walk on a graph for $S$ periods with the entries
of $c$ providing a measure of which nodes would be visited most often
in the iterative scheme $c^{\left(t\right)}=\alpha Mc^{\left(t-1\right)}+\beta$
where $c^{\left(t\right)}$ converges to $c_{P}$ for any initial
choice $c^{\left(0\right)}$ as long as $\alpha$ is chosen appropriately.
\citet{banerjee2013diffusion} introduce this idea formally via the
diffusion centrality as
\begin{equation}
c_{D}\defeq\left(\sum_{s=0}^{S}\alpha^{s}M^{s}\right)\indic,\label{eq:diffusion-def}
\end{equation}
which we can interpret as a ``finite'' walk PageRank centrality
with exogenous centrality $\beta=\indic$. Its limiting case is the
Katz centrality defined via

\begin{align}
c_{K} & \defeq(I-\alpha M)^{-1}\indic,\label{eq:katz-def}
\end{align}
which we obtain by setting $\beta=\indic$ in the PageRank centrality
\eqref{eq:pagerank-def}. We summarize the relationship between eigenvector,
Katz, PageRank, and Diffusion centralities in
\begin{lem}
\label{lem:equivalence-centr}~
\begin{enumerate}
\item \label{enu:diff-katz-equi}As the chains in the diffusion centrality
become infinitely long, it converges to the Katz centrality, i.e.
as $S\goesto\infty$, $c_{D}\goesto c_{K}$
\item \label{enu:eigenvec-equiv}As $\alpha\goesto1/\lambda_{1}$, Katz
and PageRank centralities converge to the eigenvector centrality,
i.e. $c_{K}\goesto c$ and $c_{P}\goesto c$.
\end{enumerate}
\end{lem}
For inference on PageRank and Katz centralities, we derive the Jacobian
$\diff c_{P,i}=B_{P,i}\diff\vek M$ as 
\begin{equation}
B_{P,i}=\left(-\alpha\beta^{\trans}\left(I-\alpha M^{\trans}\right)^{-1}\right)\otimes\left(e_{i}^{\trans}\left(I-\alpha M\right)^{-1}\right),\label{eq:jac-pr}
\end{equation}
while for the diffusion centrality we have
\begin{equation}
B_{D,i}=\sum_{s=0}^{S}\alpha^{-s}\sum_{j=1}^{s}\left(\indic^{\trans}M^{\trans}\right)^{s-j}\otimes e_{i}^{\trans}M^{j-1}.\label{eq:jac-diff}
\end{equation}
The following result allows conducting inference on eigenvector-based
centrality measures.
\begin{thm}
\label{thm:pagerank-katz-smooth}For an adjacency matrix estimator
that satisfies $\sqrt{n}\,\vek(\hat{M}_{n}-M)\wkc N\left(0,\Omega\right)$,
we have the following results.
\begin{enumerate}
\item \label{enu:pr-katz}The PageRank and Katz centralities in \prettyref{eq:pagerank-def}
and \prettyref{eq:katz-def} are smooth functions of $M$. Furthermore,
$\hat{M}_{n}\inprob M$ implies $\hat{c}_{P,i}\inprob c_{P,i}$ and
for an adjacency matrix with distribution $N\left(\vek M,\Omega\right)$,
we have the $t$-statistic
\[
\sqrt{g}\frac{\left(\hat{c}_{P,i}-c_{P,i}\right)}{\sigma_{P,i}}\wkc N\left(0,1\right)
\]
 for \textup{$\sigma_{P,i}=B_{P,i}\Omega B_{P,i}^{\trans}$} defined
in \prettyref{eq:jac-pr}. Similarly, setting $\beta=\indic$ in $\sigma_{P,i}$
lets us obtain the same result for the Katz centrality.
\item \label{enu:diffusion}The diffusion centrality in \prettyref{eq:diffusion-def}
is a smooth function of $M$. Furthermore, $\hat{M}_{n}\inprob M$
implies $\hat{c}_{D,i}\inprob c_{D,i}$. Further, we have 
\[
\sqrt{g}\frac{\left(\hat{c}_{D,i}-c_{D,i}\right)}{\sigma_{D,i}}\wkc N\left(0,1\right)
\]
for $\sigma_{D,i}=B_{D,i}\Omega B_{D,i}^{\trans}$, defined in \prettyref{eq:jac-diff}.
\end{enumerate}
\end{thm}
From a computational point of view, the evaluation of the Jacobian
$B_{P,i}$ can be memory-intensive and sometimes researchers may prefer
to avoid inversion altogether. It is straightforward to construct
a variance estimate based on the Moore-Penrose inverse instead. See
\citet[Ch. 8, Thm. 5]{magnus2019matrix}.

\subsection*{Degree centrality}

Lastly, we examine degree centrality, defined as $c_{N}\defeq M\indic$
for a vector of ones $\indic$ of length $p$. It is straightforward
to compute the distribution of the degree centrality via
\begin{thm}
\label{thm:deg-centrality}The degree centrality satisfies $c_{N}=\left(\indic^{\trans}\otimes I_{p}\right)\vek M$
and $\hat{M}_{n}\inprob M$ implies that $\hat{c}_{N}\inprob c_{N}$.
For an adjacency matrix estimator that satisfies $\sqrt{n}\,\vek\left(\hat{M}_{n}-M\right)\wkc N\left(0,\Omega\right)$,
we have for the $i$th coefficient 
\[
\frac{\hat{c}_{N,i}-c_{N,i}}{\sqrt{\sigma_{N,i}^{2}/g}}\wkc N\left(0,1\right),
\]
 for $\sigma_{N,i}^{2}=\left(\indic\otimes e_{i}\right)^{\trans}\Omega\left(\indic\otimes e_{i}\right)$.
\end{thm}

\section{Applications and Simulations\label{sec:Simulation-study}}

Our results cover directed networks with potential self-loops, which
imply non-symmetric adjacency matrices. For the simplest case, the
entries $M_{ij}$ are equal to one if a connection exists between
$i$ and $j$ and zero otherwise. For directed graphs, $M$ is not
symmetric and in the weighted case, its entries are proportional to
the intensity of the connection. The applications demonstrate the
utility of our methods for eigenvector centralities but equally apply
to related statistics used, e.g. by clustering algorithms.

\subsection{\label{subsec:application}Network models of international trade
and input-output}

\begin{figure}
\caption{\label{fig:one-sided-intervals}Confidence intervals ($95\%$) for
eigenvector centralities of (a) trade network (one-sided) and (b)
input-output network (two-sided).}

\centering

\begin{tabular}{cc}
\includegraphics[width=7cm,height=6cm,keepaspectratio]{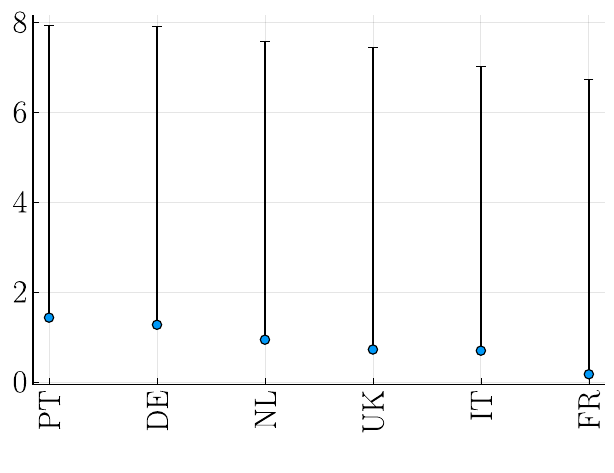} & \includegraphics[width=7cm,height=6cm,keepaspectratio]{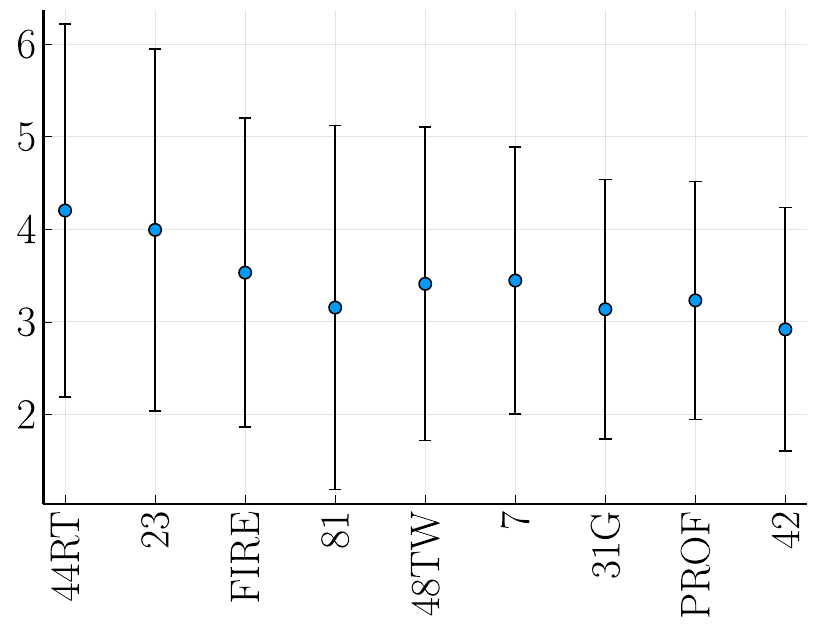}\tabularnewline
(a) Centralities for EU trade network (logarithms). & (b) Centralities for input-output network.\tabularnewline
\end{tabular}

\justifying

\noindent{\footnotesize\emph{Notes}}{\footnotesize : The plots show
estimated centrality scores based on weighted, directed trade networks
estimated with error. The left panel shows how accounting for uncertainty
in the estimated trade links provides upper bounds for the centrality
scores, though leaves the importance ordering based on the point estimates
intact. In the right panel (b), we present estimated input-output
links based on the sectors of the American economy, where we ordered
the nodes by the upper end of the confidence interval. We see that
the point estimates vary substantially and could have led to misreporting
of node importance.}{\footnotesize\par}
\end{figure}

Suppose that an empirical researcher wishes to quantify the uncertainty
in the associated centrality estimates using the results of \prettyref{sec:centralities},
in particular \prettyref{thm:eigenvec-smooth}. To demonstrate versatility,
we follow computational convention applying \prettyref{cor:folded-normal}
to $\smlabs{d_{i}}$ implying one-sided intervals for one application,
and rely on the PF theorem for our second, leading to two-sided intervals.
\prettyref{fig:one-sided-intervals} shows plots of intervals obtains
for both examples for centralities of (a) trade and (b) input-output
networks.
\prettyref{tab:sector-codes} links the abbreviations in \prettyref{fig:one-sided-intervals}
(b) to the relevant sector. We see from panel (a) that accounting
for uncertainty in the centrality score estimates leaves the importance
ranking among the countries intact, whereas it reorders it for the
sectors of the US economy when considering the upper bounds.
\begin{table}
\caption{\label{tab:sector-codes}Sector codes for the input-output network
model.}
\centering
\begin{tabular}{l|l}
\toprule   Code   & Description                                                       \\ 
\hline   23 &  Construction \\  31G &                                                     Manufacturing \\   42 &                                                   Wholesale trade \\ 44RT &                                                      Retail trade \\ 48TW &                                    Transportation and warehousing \\ FIRE &              Finance, insurance, real estate, rental, and leasing \\ PROF &                                Professional and business services \\    7 & Arts, entertainment, recreation, accommodation, and food services \\   81 &                                 Other services, except government \\    G &                                                        Government \\ \bottomrule 
\end{tabular}
\justifying

\noindent{\scriptsize\emph{Notes}}{\scriptsize : This table shows
the sectors used in estimating the input-output network adjacency
matrix. Data originate with the Bureau of Economic Analysis.}{\scriptsize\par}
\end{table}

\begin{figure}
\caption{\label{fig:estimated-networks}Estimated networks from (a) trade data
and (b) input-output of sectors of the US economy. Arrow thickness
indicates trade volume while node size indicates estimated centrality
score.}

\centering%
\begin{tabular}{c}
\resizebox{11.5cm}{!}{\includegraphics[width=9cm,height=9cm,keepaspectratio]{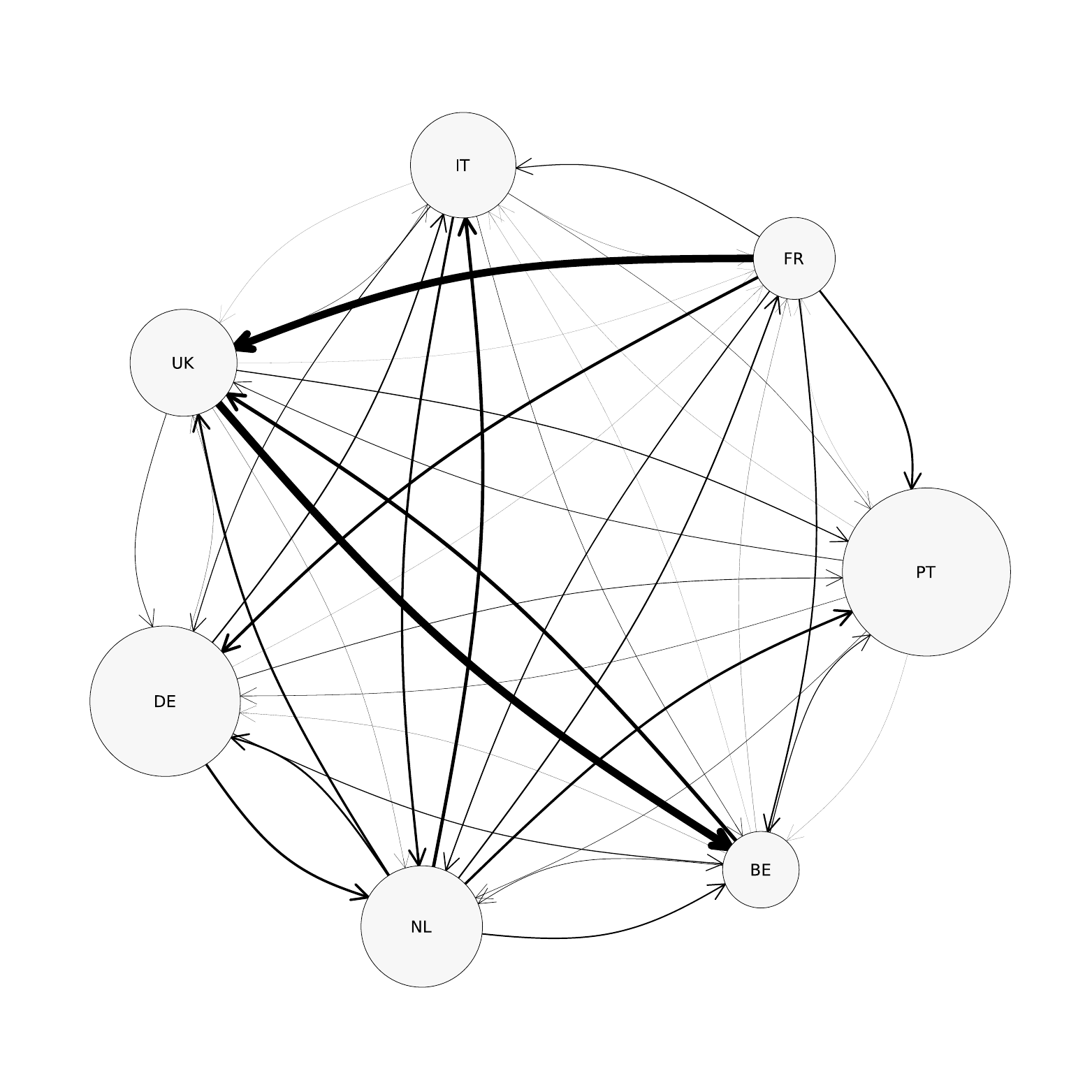}}\tabularnewline
(a) Digraph of estimated EU pharmaceutical trade network. Date range
is from '04 to '18, EUROSTAT.\tabularnewline
\begin{cellvarwidth}[t]
\centering
\resizebox{11.5cm}{!}{

\includegraphics[width=9cm,height=9cm,keepaspectratio]{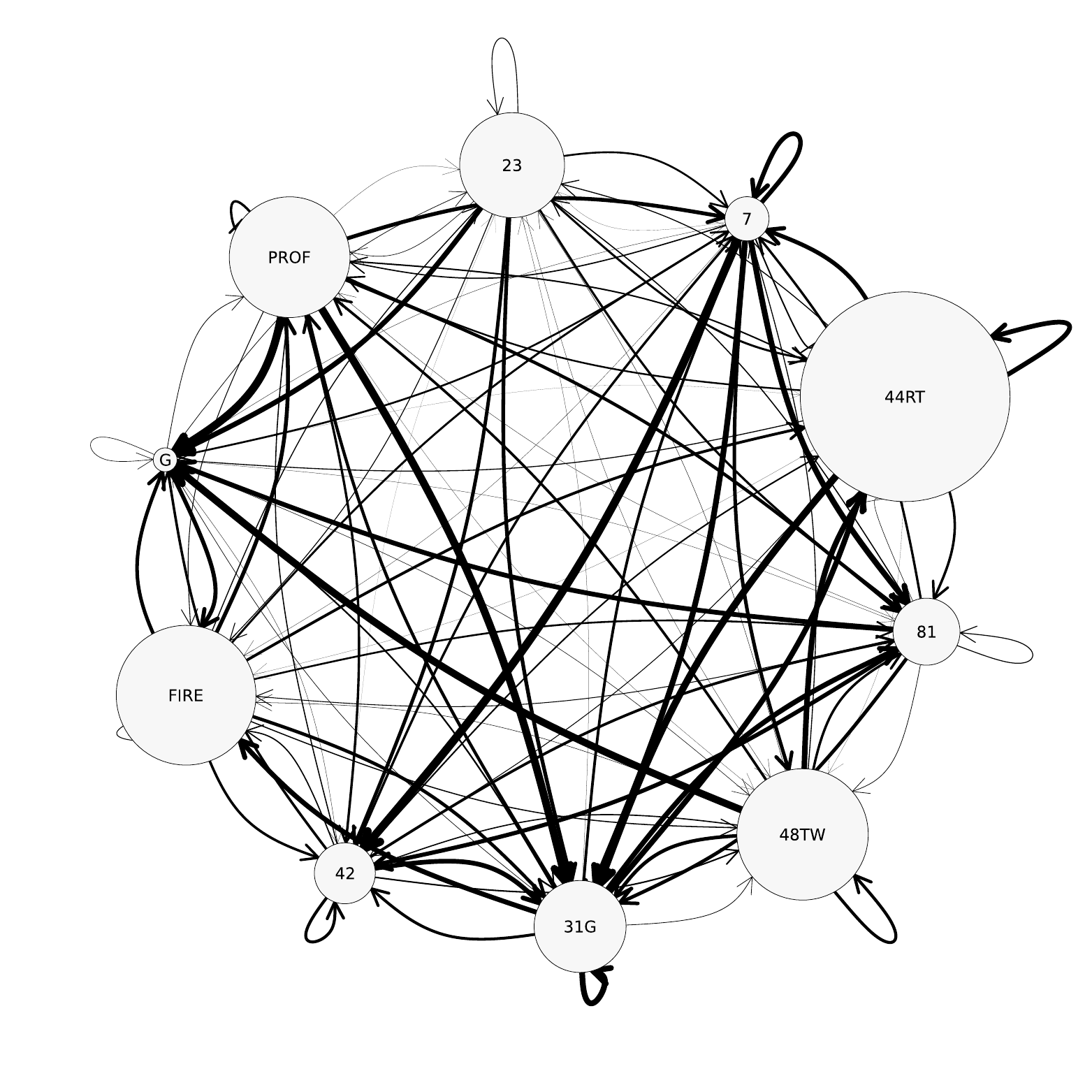}}
\end{cellvarwidth}\tabularnewline
(b) Digraph of estimated sectoral trade of US economy. Date range
is from '17 to '22, BEA.\tabularnewline
\end{tabular}

\justifying

\noindent{\footnotesize\emph{Notes: }}{\footnotesize The plots show
estimated networks where the arrow indicates the direction of trade.}{\footnotesize\par}
\end{figure}

\prettyref{fig:estimated-networks} visualizes the estimated networks,
where arrow thickness denotes the weight and node size corresponds
to centrality. Panel (a) displays the EU trade network and gives some
intuition why Portugal is such a central node, despite having lower
weights than the others. What it lacks in trade volume, it compensates
with diversity of connections and being connected to other very central
nodes, as is intended for eigenvector centrality. In panel (b), we
see the estimated network of sectors of the US economy. There, trade
tends to be fairly balanced shown by equally thick arrows. Pointwise,
the retail sector (44RT) is the largest, but it is estimated with
quite a bit of noise. Arts, food, and entertainment (7) could be much
more central than retail despite its lower point estimate. The lesson
here is clear: it is likely that arts, food, and entertainment is
more volatile than retail, so that perhaps a time-dependent graph
model may be better able to capture the dependencies and reduce the
uncertainty thus quantified.

\subsection{Digraph model}

\begin{figure}
\caption{\label{fig:eigenvector-centrality}Example graph measured with noise
and quality of the asymptotic approximation for inference on eigenvector
centralities. $1000$ MC repetitions were used for a sample size of
500. Q-Q plots are theoretical ($y$) vs. empirical ($x$).}

\centering

\begin{tabular}{cc}
\includegraphics[height=6cm]{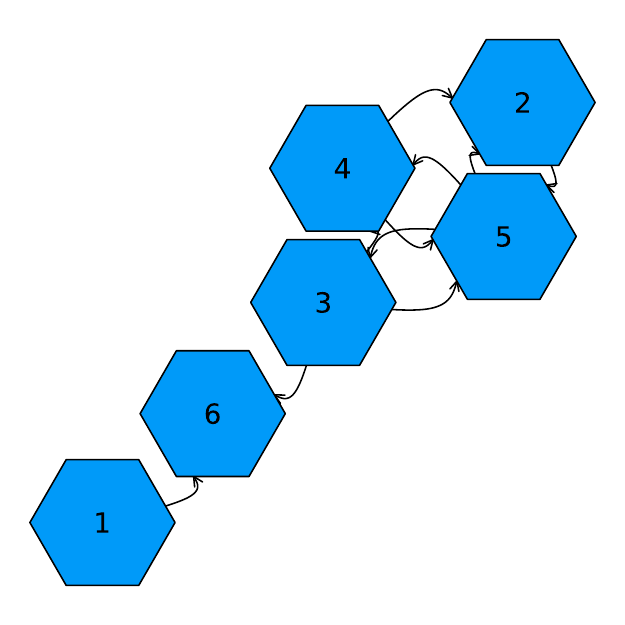} & \includegraphics[height=6cm]{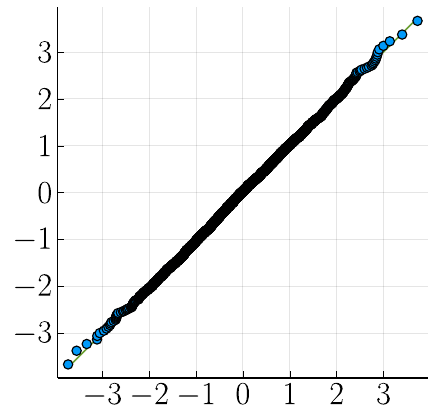}\tabularnewline
(a) Digraph with adjacency matrix $M$. & (b) Quality of approximation.\tabularnewline
\end{tabular}

\justifying

\noindent{\footnotesize Notes: The left panel shows the graph that
gives rise to adjacency matrix $M$. It has six nodes and ten connections.
On the right panel, we see the Q-Q plot of the $t$-statistic for
inference on a single centrality score on the $x$-axis and the theoretical
quantiles on the $y$-axis.}{\footnotesize\par}
\end{figure}

We studied the finite sample performance of the $t$-test based on
a random, binary directed graph with six nodes and ten connections
displayed in panel (a) of \prettyref{fig:eigenvector-centrality}.
The Q-Q plot for the $t$-test statistic in panel (b) shows that the
approximation performs very well. Our method is therefore well-suited
for binary digraphs, too.

\subsection{\label{subsec:Monte-Carlo-Evidence}Generic invariant subspaces}

\begin{figure}
\caption{\label{fig:quality-approx-invar-sub}Quality of asymptotic approximation.
Left column: Wald test, right column: t-test. $5000$ MC repetitions
were used for a sample size of $100$. : Q-Q plots are theoretical
($y$) vs. empirical ($x$).}

\centering

\begin{tabular}{ccc}
\includegraphics[width=8cm,height=6cm,keepaspectratio]{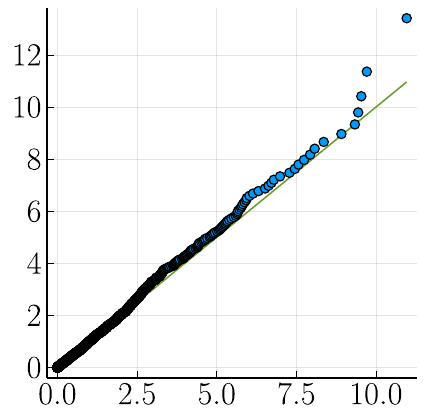} &  & \includegraphics[width=8cm,height=6cm,keepaspectratio]{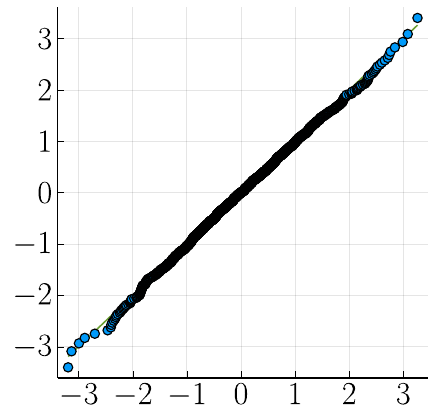}\tabularnewline
(a) Wald test quantiles. &  & (b) $t$-test quantiles.\tabularnewline
\includegraphics[width=8cm,height=6cm,keepaspectratio]{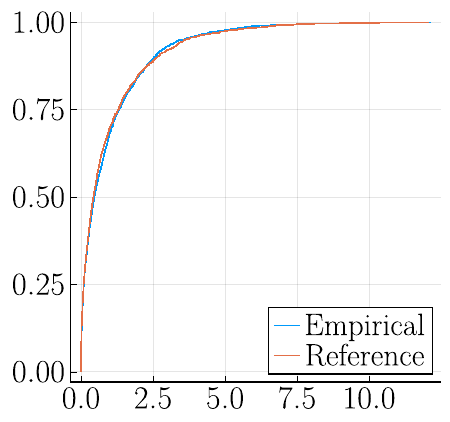} &  & \includegraphics[width=8cm,height=6cm,keepaspectratio]{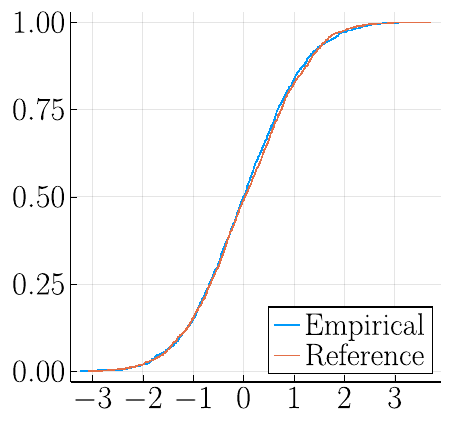}\tabularnewline
(c) CDF comparisons of Wald test statistic. &  & (d) CDF comparisons of $t$-test statistic.\tabularnewline
\end{tabular}

\justifying

\noindent{\footnotesize\emph{Notes}}{\footnotesize : The figures show
Q-Q plots and CDFs for the Wald and $t$-test statistics to verify
the accuracy of \prettyref{thm:waldstatdistrib} and \prettyref{thm:tdist}.
The Q-Q plots visualise the outliers in the tails of the Wald tests
which are less pronounced for the $t$-test. However, even for the
Wald-test, the number of outliers is only moderate in light of the
$5,000$ MC repetitions. The bottom panels, (c) and (d), show empirical
and reference CDFs which show close tracking across the entire support.} 
\end{figure}

To check how well our methods worked for inference on generic invariant
subspaces, we used a dense matrix with $p=2$ whose stacked columns
have simplified covariance matrix $\Omega=I_{p}\otimes\Omega_{M}$.
Full details on the data-generating process appear in \prettyref{alg:wald-and-t}
in \prettyref{app:results-network-centralities}. The online supplement
to this paper contains further DGPs as well as examples on how to
include heteroskedasticity and autocorrelation-robust covariance matrix
estimators. We ran further simulations to study the performance of
the $t$- and Wald tests using Q-Q plots as well as the empirical
cumulative distribution functions compared with their theoretical
counterparts. These visual aids demonstrate the quality of the asymptotic
approximations found in \prettyref{sec:hyp-tests}.\footnote{Complementary to the results presented in this section, we refer the
reader to the extra material hosted at \href{https://github.com/jsimons8/networkmodelssubspaces}{https://github.com/jsimons8/networkmodelssubspaces.}} The panels in \prettyref{fig:quality-approx-invar-sub} show that
both the Wald and $t$-test statistics perform well in simulation
exercises. Left panels let us judge the approximation made in \prettyref{thm:waldstatdistrib}
while right ones display the quality of the $t$-test approximation
of \prettyref{thm:tdist}. The overall performance is very good with
only few outliers. Histograms overlain with densities displayed similar
results and appear in the online supplement. The bottom panels show
the same pattern where the cumulative distribution functions track
their empirical counterparts well.

\subsection{Singular subspace inference: Monte Carlo evidence}

We considered a data-generating process using a covariance matrix
$\Omega=I_{m}\otimes\Omega_{W}$ for a positive-definite $\Omega_{W}\in\reals^{m\times m}$
for $m=3$. \prettyref{alg:wald-and-t-svd} in \prettyref{app:results-network-centralities}
details the steps of how we generate samples for the $t$- and Wald
statistics.
\begin{figure}
\caption{\label{fig:quality-approx-svd}Quality of asymptotic approximation
for Wald test-based inference on singular vectors.  $2000$ MC repetitions
were used for a sample size of $500$. : Q-Q plots are theoretical
($y$) vs. empirical ($x$).}
 \centering

\begin{tabular}{cc}
\includegraphics[width=8cm,height=6cm,keepaspectratio]{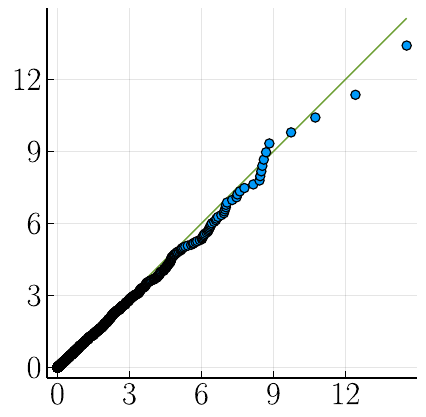} & \includegraphics[width=8cm,height=6cm,keepaspectratio]{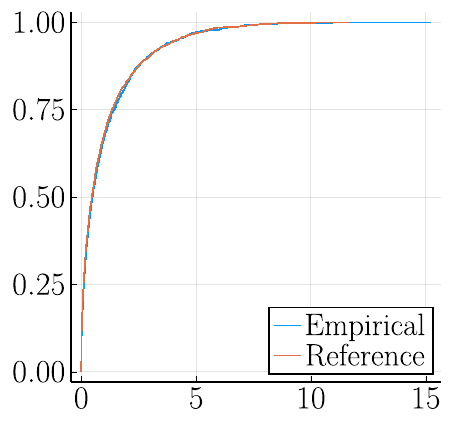}\tabularnewline
(a) Wald test quantiles. & (b) CDF comparisons of Wald test statistic.\tabularnewline
\end{tabular}

\justifying

\noindent{\footnotesize\emph{Notes}}{\footnotesize : The plots show
performance of the Wald test, both in terms of quantiles and the CDF.
In panel (a), we see only a few outliers and generally good agreement
between empirical and reference distributions throughout the support
in panel (b).}{\footnotesize\par}
\end{figure}

The Q-Q plot in panel (a) of \prettyref{fig:quality-approx-svd} shows
that the asymptotic approximation for the SVD performs well with few
outliers. Similarly, panel (b) shows the case for $q=m=1$ and that
the empirical CDF tracks the implied $\chi_{1}^{2}$ benchmark well.
These results are encouraging that our first-order expansion of the
SVD map delivers a good approximation for conducting inference on
singular subspace vectors. Details on the DGP appear in \prettyref{alg:wald-and-t-svd}.

\section{Conclusion\label{sec:Conclusion}}

This paper has extended the inferential theory of \citet{Tyler1981}
to cover non-diagonalizable matrices and applied the results to network
statistics. In addition to the Wald test for full vector hypotheses,
a $t$-test is practically useful because it allows inference on individual
coefficients. The method of smooth eigenvector estimation in \citet{10.1093/biomet/asad018}
presents a useful extension to possibly avoid introducing outliers
through the proposed normalization and to take advantage of infinite
differentiability.

Regarding the underlying perturbation theory, there are strong resemblances
between \citet{Sun1991} and \citet{Kato} although the former also
discusses the extensions to Gateaux derivatives. We are hopeful that
the present exposition can aid researchers in similar settings and
present a way to find Jacobians of maps that are related to invariant
or singular subspaces. We refer readers new to the literature on perturbation
theory to \citet{2greenbaum:2019} who offer a pedagogic and detailed
treatment of the subject. We are hopeful that our arguments are easily
adaptable for statistics that depend on matrix subspaces in a more
general way both in the graph domain and others.

The leading application of our results is to invariant subspaces of
estimated network adjacency matrices, which inherit the measurement
error of the estimated network. Similarly, our results on singular
vectors may be applied to network clustering algorithms to allow quantifying
the uncertainty in low-dimensional representations of networks. Generally,
invariant subspaces find applications in the analysis of not only
adjacency matrices but also graph Laplacians where they enable algorithms
for spectral clustering, community detection, or the finding of mixing
rates for random walks on graphs. In this vein, using covariate-level
information as suggested in \citet{BinkiewiczVogelsteinRohe2017}
could be useful to constrain invariant subspaces and potentially tighten
confidence intervals.

Network centralities are often used to identify interventions. However,
point estimates may be misleading if the uncertainty in the network
link identification is large. Therefore, we advocate reporting confidence
bands. 

The higher-order Davis-Kahan bounds reveal important distinctions
between general invariant subspace and eigenspace perturbations. In
the case of the former, the eigengap dominates the perturbation estimations
while for the latter, higher-order bounds simply consist of the first-order
Davis-Kahan bound raised to a higher power. The consequence is that,
for eigenspaces, statisticians only need control over the Davis-Kahan
bound while for invariant subspaces, the eigenvalue gap is more important
than the estimation precision.

We have also shown in \prettyref{thm:higher-order-dk} and \prettyref{thm:stoch-bound-prop}
that the eigengap governs the convergence speed of subspace-based
network statistics and explicitly calculated convergence rates of
clustering coefficients whenever the network is estimated using a
single large matrix. A useful extension would be to apply methods
similar to those in \citet{10.1093/biomet/asac032} where a considered
network statistic is the maximum eigenvalue of the adjacency matrix
corresponding to a random eigengap.

The inference results assume fixed matrix size but the convergence
rate and perturbation bound calculations are agnostic about matrix
size. Therefore, an extension of this study is to consider our results
for large, random matrices in more detail, which we leave for future
work.

\bibliographystyle{chicago}
\bibliography{eigenvector-inference,biometrika-one-jasa-references,biometrika_spectral_graph,empirical-applications}
\appendix
\appendixpage

\subsection*{Organization}

This appendix is organized as follows: \prettyref{app:projections-and-gen-inverses}
reviews projections as well as the relationship between perturbation
expansions and derivatives and shadows \prettyref{sec:setup}. \prettyref{app:distributions}
shadows \prettyref{sec:hyp-tests}, \nameref{sec:hyp-tests}. \prettyref{app:Jacobians}
derives expressions for the Jacobians and shadows \prettyref{sec:jacobians-ho-perturbations}.
\prettyref{app:using-ift} derives results of smoothness and infinite
differentiability of invariant and singular subspace maps. Finally,
\prettyref{app:results-network-centralities} shadows \prettyref{sec:centralities},
\nameref{sec:centralities}.

\section{Linear algebra fundamentals and further details for \prettyref{sec:setup}:
\nameref{sec:setup}\label{app:projections-and-gen-inverses}}

In this section, we follow \citet{Tyler1981} and briefly review the
theory of projections onto subspaces as well as other fundamentals
of linear algebra.

Regarding \prettyref{exa:symmetrisable}, we illustrate the requirement
by \citet{Tyler1981} that we relax. Define the inner product $\inprd uv_{\Gamma}\defeq u^{\trans}\Gamma v$
so that $M$ satisfies 
\[
\inprd{Mu}v_{\Gamma}=\inprd u{Mv}_{\Gamma}
\]
or written out $M^{\trans}\Gamma=\Gamma M$. As $\Gamma$ is positive-definite,
it has a symmetric square root. Define $B\defeq\Gamma^{\frac{1}{2}}M\Gamma^{-\frac{1}{2}}$
and find 
\begin{align*}
B^{\trans} & =\Gamma^{-\frac{1}{2}}M^{\trans}\Gamma^{\frac{1}{2}}\\
 & =\Gamma^{-\frac{1}{2}}M^{\trans}\Gamma\Gamma^{-1}\Gamma^{\frac{1}{2}}\\
 & =\Gamma^{-\frac{1}{2}}M^{\trans}\Gamma\Gamma^{-\frac{1}{2}}\\
 & =\Gamma^{-\frac{1}{2}}\Gamma M\Gamma^{-\frac{1}{2}}\\
 & =\Gamma^{\frac{1}{2}}M\Gamma^{-\frac{1}{2}}
\end{align*}
where the last line equals $B$. Hence, $M$ is similar to a symmetric
matrix and therefore has real eigenvalues.

To introduce inference on invariant subspaces, we will use the case
of diagonalizable $M$ and introduce eigenvalues (or singular values)
where necessary. Associating subspace vectors with a set of such values,
we assume in either case that the elements of $\eigs$ split according
to $\eigs:=\eigs_{I}\union\eigs_{J}$.

Throughout and without loss of generality, we assume that $\rank\upsilon=r\leq\abs{\eigs_{I}},$
i.e. the number of linearly independent columns in $\upsilon$ is
less than or equal to the dimension of the subspace spanned by the
vectors associated with elements in $\eigs_{I}$. If $\rank\upsilon>\abs{\eigs_{I}}$,
i.e. the number of linearly independent columns in $\upsilon$ exceeds
the number of elements in $\eigs_{I}$, we define the alternative
null hypothesis 
\[
H_{0}^{\ast}\,:\,\risof\in\spn\,\upsilon,
\]
which can be tested by constructing the $p\times p-q$ matrix $\upsilon_{\perp}$
such that 
\begin{equation}
\upsilon_{\perp}^{\trans}\upsilon=0.\label{eq:per-def}
\end{equation}

The case of $\rank\upsilon\leq\abs{\eigs_{I}}$ is intuitive because
the columns of $\upsilon$ span a lower-dimensional subspace than
$\risof$ while $\rank\upsilon>\abs{\eigs_{I}}$ is effectively an
overdetermination of the subspace of interest so that the estimand
$\risof$, under the null hypothesis, is wholly contained in the candidate
space $\spn\,\upsilon$. For eigenvectors, \citet{Tyler1981} implements
this scenario by inverting the hypothesis \prettyref{eq:basic-hypothesis-random-matrix}
via orthocomplementation to achieve \prettyref{eq:perp-inference}
which reduces to $H_{0}$ again. The orthocomplement, $\upsilon_{\perp}$,
is always easy to compute from $\upsilon$ using, e.g. an LU factorization
or the Gram-Schmidt algorithm. Moreover, working with $H_{0}^{\ast}$
and hence $\upsilon_{\perp}$ implies normalization-invariant test
statistics and offers numerical advantages as Monte Carlo experiments
showed. We shall therefore write our arguments using $\upsilon_{\perp}$
is equivalent to inference on $\upsilon$. From a practical perspective,
a researcher can test for a maximum of $q$ linearly independent vectors
that may lie in $\risof$, which corresponds to $p-q$ independent
elements in $\upsilon_{\perp}$. Because we prefer to write hypotheses
about $\upsilon_{\perp}$, researchers choose maximally $m$ columns
for $\upsilon$ with $m\leq q$ which translates to $p-m$ columns
for $\upsilon_{\perp}$.

Take the setup in \prettyref{assu:invar-sub-ass} and assume further
that the eigenvalues of $\Lambda_{I}$ are semi-simple so that each
of them corresponds to a linearly independent column of $R_{I}$.
To make $M$ diagonalizable, set $\Lambda_{IJ}=0$. Define the eigenprojection
onto the subspace spanned by eigenvectors associated with roots in
$\eigs_{I}$ as
\[
P_{I}:=\sum_{\lambda\in\mathcal{L}_{I}}\lambda r_{\lambda}l_{\lambda}^{\trans},
\]
where $r_{\lambda}$ and $l_{\lambda}$ belong to root $\lambda$.
Now, consider
\begin{example}
\label{exa:projectors} Let 
\begin{equation}
M=\begin{bmatrix}0.8 & 0.5\\
0 & 0.4
\end{bmatrix}.
\end{equation}
Then, the right eigenvectors are $r_{1}^{\trans}=\left[\begin{smallmatrix}1 & 0\end{smallmatrix}\right]$
and $r_{2}^{\trans}=\left[\begin{smallmatrix}1 & -0.8\end{smallmatrix}\right]$
with eigenvalues are $\lambda_{1}=0.8$ and $\lambda_{2}=0.4$. Its
left eigenvectors are $l_{1}^{\trans}=\left[\begin{smallmatrix}1 & 1.25\end{smallmatrix}\right]$
and $l_{2}^{\trans}=\left[\begin{smallmatrix}0 & 1.6\end{smallmatrix}\right]$.
The associated eigenprojections are
\begin{equation}
P_{\lambda_{i}}=r_{i}l_{i}^{\trans}\label{eq:outer_proj}
\end{equation}
for $i\in\left\{ 1,2\right\} $ the pair of matrices
\begin{align*}
P_{\lambda_{1}} & =\begin{bmatrix}1 & 1.25\\
0 & 0
\end{bmatrix} &  & P_{\lambda_{2}}=\begin{bmatrix}0 & -1.25\\
0 & 1
\end{bmatrix}
\end{align*}
We see that $P_{\lambda_{i}}$ are rank-1 projection matrices for
which $P_{\lambda_{i}}P_{\lambda_{j}}=P_{\lambda_{i}}\delta_{ij}$
and if $\upsilon$ belongs to an eigenspace associated with root $\lambda_{i}$,
we have $\upsilon^{\trans}P_{\lambda_{i}}=\upsilon^{\trans}$ and
$P_{\lambda_{1}}+P_{\lambda_{2}}=I$. Finally, under the maintained
null hypothesis that $\upsilon\in\risof$ for $I=\left\{ 1\right\} $,
\[
P_{\lambda_{1}}\upsilon_{\perp}=0.
\]
\end{example}
Based on eigenprojections, it is easy to construct generalized inverses
$A^{+}$ such that $AA^{+}A=A$. For a matrix of rank $s$, let the
eigenprojection of $A$ be as in \prettyref{exa:projectors}. Then,
a generalized inverse
\begin{equation}
A^{+}=\sum_{i=1}^{s}\lambda_{i}^{-1}P_{\lambda_{i}}.\label{eq:gen-inv-eig}
\end{equation}
If we wish to relax the diagonalizability assumption, we may easily
construct a generalized inverse via a singular value decomposition,
which exists independently of diagonalizability. Proceed via
\begin{equation}
A^{+}=\sum_{i=1}^{s}\iota_{i}^{-1}v_{i}u_{i}^{\trans},\label{eq:gen-inv-svd}
\end{equation}
where $\left\{ \iota_{i}\right\} _{i=1}^{s}$ collects the non-zero
singular values and $u_{i}$ and $v_{i}$ are left- and right-singular
vectors of $A$. For the purpose of inverting estimated covariance
matrices, either \prettyref{eq:gen-inv-eig} or \prettyref{eq:gen-inv-svd}
would be suitable.

\subsection{Covariance matrices\label{app:Covariance-matrices}}

Regarding estimation of associated covariance matrices belonging to
$\hat{M}_{T}\defeq T^{-1}\sum_{t=1}^{T}M_{t}$ is, 
\[
\hat{\Omega}\defeq T^{-1}\sum_{t=1}^{n}\left(\vek\hat{M}_{t}-\vek M_{t}\right)\left(\vek\hat{M}_{t}-\vek M_{t}\right)^{\trans},
\]
which we can make more efficient if we know that $\Omega=I_{p}\otimes\Omega_{M}$
by enforcing this structure during estimation. If the index $t$ refers
to time, we can also accommodate autocorrelation in the residuals,
see \citet{b76ccb64-7fa5-32f9-a4fe-72298146be7d} and the empirical
application in \prettyref{subsec:application}. To be able to construct
Wald tests, we need to consistently estimate $\Omega^{-1}$. However,
neither the rank nor a generalized inverse of $\Omega$ are continuous.
To deal with this issue, we recommend the procedure outlined in \citet{LUTKEPOHL1997315}
and invoke \citet[Lemma 2.1]{doi:10.1080/01621459.2023.2202435} to
establish that $\hat{\Omega}^{+}\inprob\Omega^{-1}$.

An important special case of covariance arises when the columns of
the error matrix are homoskedastic, which we define for some suitable
$\Omega_{M}$ and $\varepsilon_{t}\defeq\ensuremath{\left[\textbf{a}_{1}\,\textbf{a}_{2}\,\dots\,\textbf{a}_{p}\right]_{t}}$
via 
\begin{equation}
\mathbb{E}\textbf{a}_{i}\textbf{a}_{j}^{\trans}=\Omega_{M}\delta_{ij}\label{eq:homoskedastic-cols}
\end{equation}
so that $\Omega=I_{p}\otimes\Omega_{M}$. None of our results are
sensitive to these restrictions and, if equipped with a suitable covariance
matrix estimator, one could allow for autocorrelated errors and heteroskedasticity
across time, as well, an application of which we consider in \prettyref{subsec:application}.

\section{Results and proofs for \prettyref{sec:hyp-tests}: \nameref{sec:hyp-tests}}

\label{app:results-hyp-tests}\label{app:distributions}

Before proving \prettyref{lem:proj-distrib}, we establish how the
first derivative \prettyref{eq:firstderivative-1} allows assembly
of a first-order expansion, which we pursue in the following result.
Furthermore, we show how the ingredients of Wald statistics behave
individually afterwards.

Define 

\begin{equation}
\begin{bmatrix}\Delta_{I}\\
\Delta_{\text{JI}} & \Delta_{J}
\end{bmatrix}\defeq\begin{bmatrix}L_{I}^{\trans}ER_{I}\\
L_{J}^{\trans}ER_{I} & L_{J}^{\trans}ER_{J}
\end{bmatrix},\label{eq:argument-1-1-1}
\end{equation}
which corresponds to the block-diagonalization of the perturbing matrix
$E=\tilde{M}-M$.

\subsection{Consistency}

To establish consistency, we have
\begin{thm}[Consistency of invariant subspaces]
\label{thm:consistency}Under \prettyref{assu:general-ass}, the
following hold:
\begin{enumerate}
\item \label{enu:estim}$\col\hat{M}_{n}\inprob\col M$.
\item \label{enu:.right-cons}$\col\hat{R}_{n}\inprob\col R$.
\item \label{enu:left-cons}$\col\hat{L}_{n}\inprob\col L$.
\item \label{enu:consistent-cov-mat-estimator}$\hat{\Omega}_{n}\inprob\Omega$.
\item \label{enu:root-cons}$\hat{\Lambda}\inprob\Lambda$ up to reordering
of eigenvalues within blocks $I$ and $J$.
\item \textup{\label{enu:B-cons}$\hat{B}\inprob B$.}
\item \label{enu:.cov-mat}$\hat{B}_{j}\hat{\Omega}\hat{B}_{j}^{\trans}\inprob B_{j}\Omega B_{j}$
for $j\in\left\{ \emptyset,D^{\trans},W\right\} $ denoting the generic,
basis vector, and Wald test Jacobians.
\item \label{enu:estimator}$\hat{D}_{n}\inprob D$.
\end{enumerate}
\end{thm}
Analogously to \prettyref{thm:consistency}, we have
\begin{cor}[Consistency of singular subspaces]
\label{cor:consistency-svd}Under \prettyref{assu:svd-ass}, the
following hold:
\begin{enumerate}
\item \label{enu:B-svd-cons}$\hat{B}_{\text{SVD}}\inprob\jacsvd$.
\item \label{enu:right-sg-cons}$\col\hat{U}_{I}\inprob\col U_{I}$.
\end{enumerate}
\end{cor}
The above results guarantee consistent estimation of the Jacobian
associated with the singular subspace map and right-singular subspace.
Note that \prettyref{thm:consistency}\prettyref{enu:estim}-\ref{enu:left-cons}
refer to the fact that we only identify column spans of these matrices.
For \prettyref{thm:consistency}\prettyref{enu:consistent-cov-mat-estimator},
we need to identify the matrix explicitly as $\Omega$ must be symmetric.
If only the column span of $\Omega$ was identified, we would obtain
invariance under scalar multiplication of its columns which would
alter the row spans. The main utility of \prettyref{thm:consistency}
is the application of Slutsky's theorem to establish distributional
results of our test statistics.

\subsection{Distributions of invariant subspace estimators}

Define the short-hand $P\defeq RL^{\trans}$.
\begin{lem}
\label{lem:expansion-and-convergence}~
\begin{enumerate}
\item \label{enu:expansion}The standardized term admits the expansion centered
around the true projection matrix\textup{
\[
\sqrt{n}\upsilon_{\perp}^{\trans}\left(\hat{P}_{n,I}-P_{I}\right)=\sqrt{n}\upsilon_{\perp}^{\trans}\hat{R}_{n,J}\boldsymbol{S}^{-1}\left(\Delta_{n,\text{IJ}}\right)\hat{L}_{n,I}^{\trans}+\sqrt{n}r\left(\hat{M}_{n}-M\right).
\]
}
\item \label{enu:remainder-term}The remainder term obeys $r\left(\hat{M}_{n}-M\right)=O\left(n^{-1}\right).$
\item \label{enu:order}Let $\hat{\Delta}_{n,i}$ with $i\in\left\{ I,\text{IJ},J\right\} $
be as in \prettyref{eq:argument-1-1}. Then
\[
\begin{bmatrix}\hat{\Delta}_{n,I}\\
\hat{\Delta}_{n,\text{IJ}} & \hat{\Delta}_{n,J}
\end{bmatrix}=\begin{bmatrix}O_{p}\left(n^{-1/2}\right)\\
O_{p}\left(n^{-1/2}\right) & O_{p}\left(n^{-1/2}\right)
\end{bmatrix}.
\]
\end{enumerate}
\end{lem}
\prettyref{lem:expansion-and-convergence}\prettyref{enu:expansion}
lets us derive an estimator of the covariance matrix while \prettyref{lem:expansion-and-convergence}\prettyref{enu:remainder-term}
ensures that the approximation error is $O\left(n^{-1/2}\right)$
and therefore converges in probability to zero. Our next result ensures
the desired convergence. Define $B_{\text{IJ}}^{\trans}\defeq R_{I}\otimes L_{J}$.
Then, we have
\begin{lem}
\label{lem:deriv-conv}Let $\Delta_{n,\text{IJ}}$ be as in \eqref{eq:argument-1-1-1}.
Then, for $\upsilon_{\perp}\in\reals^{p\times\left(p-m\right)}$ ,
\begin{enumerate}
\item $\sqrt{n}\vek\Delta_{n,\text{IJ}}\wkc N\left(0,B_{\text{IJ}}\Omega B_{\text{IJ}}^{\trans}\right).$
\item \label{enu:proje}$\sqrt{n}\vek\left\{ \upsilon_{\perp}^{\trans}\hat{R}_{n,J}\boldsymbol{S}^{-1}\left(\Delta_{n,\text{IJ}}\right)\hat{L}_{n,I}^{\trans}\right\} \wkc N\left(0,\Omega_{W}\right).$
\item \label{enu:deg}$\rank\Omega_{W}=qm$.
\end{enumerate}
\end{lem}
\begin{proof}[Proof of \prettyref{lem:deriv-conv}.]
~ 
\begin{enumerate}
\item By \prettyref{lem:derivative-1}, 
\[
\vek\hat{\Delta}_{n,\text{IJ}}=\left(R_{I}^{\trans}\otimes L_{J}^{\trans}\right)\vek\left(\hat{M}_{n}-M\right).
\]
The result follows by the continuous mapping theorem and \prettyref{assu:general-ass}.
\item By the preceding argument, and the definition of $\boldsymbol{S}^{-1}$
in \prettyref{lem:derivative-1}, we see that
\[
\upsilon_{\perp}^{\trans}R_{n,J}\boldsymbol{S}^{-1}\left(\hat{\Delta}_{n,\text{IJ}}\right)L_{n,I}^{\trans}
\]
is a linear and continuous transformation of $\hat{\Delta}_{n,\text{IJ}}$.
Vectorizing $\boldsymbol{S}^{-1}\left(\hat{\Delta}_{n,\text{IJ}}\right)$
and taking the expectation of the outer product provides the result.
\item To establish the rank of the covariance matrix, we shall take the
constituent elements apart. For $\Omega_{W}=B_{W}\Omega B_{W}$ where
we reproduce \eqref{eq:jacobian} for convenience as
\[
B_{W}=\left(L_{I}\otimes\upsilon_{\perp}^{\trans}R_{J}\right)\left\{ \left(\Lambda_{I}^{\trans}\otimes I_{J}\right)-\left(I_{I}\otimes\Lambda_{J}\right)\right\} ^{-1}\left(R_{I}^{\trans}\otimes L_{J}^{\trans}\right).
\]

First, we establish that $\rank\upsilon_{\perp}^{\trans}R_{J}=m$.
Let $r\defeq p-q$. By \prettyref{assu:general-ass}\prettyref{enu:semisimple},
$R_{J}\in\reals^{p\times r}$ has $r$ linearly independent columns
so that $\rank R_{J}=r$. By the properties of ranks as linearly independent
columns we obtain the inequality $\rank\upsilon_{\perp}^{\trans}R_{J}\leq\min\left(p-m,p-q\right)\leq p-m$
where the last inequality follows from $m\leq q$. By Sylvester's
rank inequality, we obtain the lower bound $m+r-p\leq\rank\upsilon_{\perp}^{\trans}R_{J}\leq p-m$,
which implies $m-q\leq\rank\upsilon_{\perp}^{\trans}R_{J}\leq p-m$.
To see that the lower bound is slack and the upper bound is tight,
we induct on $q$: put $q=1$ and assume without loss of generality
that $m=r$. Then, $r-1\leq\rank\upsilon_{\perp}^{\trans}R_{J}\leq q$,
so premultiplication of $\upsilon_{\perp}^{\trans}$ causes the rank
of $R_{J}$ to drop by one, which implies that one column of $\upsilon_{\perp}^{\trans}$
is perpendicular to one column of $R_{J}$. However, by \prettyref{assu:general-ass}\prettyref{enu:no-evs},
we have that no columns parallel to any in $R_{I}$ must reside in
$R_{J}$, which means that $\rank\upsilon_{\perp}^{\trans}R_{J}=r$.
Therefore, per the induction hypothesis the claim holds for $q=1$.
Now, we show that it holds for $q+1$. Suppose that $\rank\upsilon_{\perp}^{\trans}R_{J}=r-q-1$,
which implies that we have again one column of $\upsilon_{\perp}^{\trans}$
that is perpendicular to a column of $R_{J}$, leading to the same
contradiction as before. Therefore, for any $q$ and $q+1$ we must
have that the lower bound has to be slack while the upper bound binds,
thus establishing the claim for $m=r$. As there was nothing special
about this choice of $m$, the claim also holds for all $m\leq r$.

Next, we use the property of the Kronecker product $\rank L_{I}\otimes\upsilon_{\perp}^{\trans}R_{J}=\rank L_{I}\rank\upsilon_{\perp}^{\trans}R_{J}=qm$.
The Jacobian $\left\{ \left(\Lambda_{I}^{\trans}\otimes I_{J}\right)-\left(I_{I}\otimes\Lambda_{J}\right)\right\} ^{-1}$
is a $qr\times qr$ square matrix of full rank. Then, $\rank\left(R_{I}^{\trans}\otimes L_{J}^{\trans}\right)=qr$
so that $\rank\left\{ \left(\Lambda_{I}^{\trans}\otimes I_{J}\right)-\left(I_{I}\otimes\Lambda_{J}\right)\right\} ^{-1}\left(R_{I}^{\trans}\otimes L_{J}^{\trans}\right)=qr$.
Putting $A_{1}:=\left(R_{I}^{\trans}\otimes L_{J}^{\trans}\right)$
and $A_{2}:=\left\{ \left(\Lambda_{I}^{\trans}\otimes I_{J}\right)-\left(I_{I}\otimes\Lambda_{J}\right)\right\} ^{-1}\left(R_{I}^{\trans}\otimes L_{J}^{\trans}\right)$
we have $qm\pm qr\leq\rank AB\leq\min\left(qm,qr\right)=qm$ where
the first relation follows from Sylvester's rank inequality and the
last equality follows from $m\leq r$. Hence, $qm\leq\rank A_{1}A_{2}\leq qm$
or $\rank A_{1}A_{2}=\rank B_{W}=qm$. The same argument then establishes
that $\rank B_{W}\Omega B_{W}^{\trans}=qm$: in detail, $qm\pm p^{2}\leq\rank B_{W}\Omega^{1/2}\leq\min\left(qm,p^{2}\right)=qm$.
Then, $qm\pm p^{2}\leq\rank B_{W}\Omega B_{W}^{\trans}\leq\min\left(qm,p^{2}\right)=qm$
and the claim follows.
\end{enumerate}
\end{proof}
Considering \prettyref{lem:proj-distrib}, we have the
\begin{proof}[Proof of \prettyref{lem:proj-distrib}.]
\prettyref{lem:proj-distrib} follows from application of \prettyref{lem:deriv-conv}\prettyref{enu:proje}
to \prettyref{lem:expansion-and-convergence}\prettyref{enu:expansion}.
For singular subspaces, we apply \prettyref{lem:jacobian}\prettyref{enu:singular}
and the delta method to $\sqrt{n}\,\vek\upsilon_{\perp}^{\trans}\hat{U}_{I}\left(\hat{M}_{n}\right)$
to obtain the result, where $\jacsvd$ is the Jacobian of the map
$\Psi$ defined in \eqref{eq:sin-subspace-map}.
\end{proof}
Recall \prettyref{thm:estdist} that asserts that \prettyref{assu:general-ass}
implies 
\[
n^{1/2}\vek\{\hat{D}_{n}^{\trans}-D^{\trans}\}\wkc N(0,B_{D^{\trans}}\Omega B_{D^{\trans}}^{\trans}),
\]
where \prettyref{eq:ahat-1} defines $\hat{D}_{n}$.
\begin{proof}[Proof of \prettyref{thm:estdist}.]
The result follows from applying \prettyref{lem:deriv-conv} to the
expansion in \prettyref{cor:expansion-and-convergence-estimator}
by substituting $R_{n,2,I}^{-1}$ for $L_{n,I}^{\trans}$.
\end{proof}
\begin{proof}[Proof of \prettyref{thm:consistency}.]
\prettyref{thm:consistency}\prettyref{enu:estim} follows from a
law of large numbers for $\hat{M}_{n}=M+n^{-1}\sum_{i=1}^{n}\varepsilon_{i}$
and $\expect\varepsilon_{i}=0$ where convergence in column span is
implied by pointwise convergence of $\hat{M}$. Similarly, \prettyref{thm:consistency}\prettyref{enu:consistent-cov-mat-estimator}
follows from a standard argument: let $\hat{\varepsilon}_{i\left(j\right)}$
denote the $j$th column of the $i$th residual matrix $\hat{\varepsilon}_{i}\in\reals^{p\times p}$.
Then, $\hat{\Omega}=\frac{1}{np}\sum_{j=1}^{p}\sum_{i=1}^{n}\hat{\varepsilon}_{i\left(j\right)}\hat{\varepsilon}_{i\left(j\right)}^{\trans}\inprob\Omega$
where we have exploited the homoskedasticity across columns of $\hat{M}-M$.
Then, \prettyref{lem:neudeck-deriv} implies that perturbation expansions
of $\psi\left(M\right)$ in terms of basis vectors of right eigenvectors
are equal to their map as defined implicitly by \prettyref{eq:inv-subspace-def}.
Therefore, $\psi\left(M\right)$ is smooth in $M$ and thus an immediate
consequence of analyticity, which follows from $\abs{\left(\Lambda{}_{I}^{\trans}\otimes I_{J}\right)-\left(I_{I}\otimes\Lambda_{J}\right)}\neq0$
by \prettyref{assu:general-ass}\prettyref{enu:no-evs} so that the
condition in the second display at the top of \citet[p. 90]{Sun1991}
is satisfied. For the case of real matrices, we appeal to \citet[Thm.~XIV.2.1]{Lang93}.
Therefore, $\hat{M}_{n}\inprob M$ implies that $\hat{R}_{n}\inprob R$
because the implicit map for $R$ defined by $MR=R\Lambda$ is analytic.
\prettyref{thm:consistency}\prettyref{enu:left-cons} follows from
the continuous mapping theorem and $\hat{L}_{n}^{\trans}=\hat{R}_{n}^{-1}$
where we note that matrix inversion is continuous. For \prettyref{thm:consistency}\prettyref{enu:root-cons},
$\hat{\Lambda}_{n}\inprob\Lambda$, we appeal to \citet[Thm.~3.9.1]{tyrtyshnikov1997brief}.
\prettyref{thm:consistency}\prettyref{enu:B-cons}-\prettyref{thm:consistency}\prettyref{enu:.cov-mat}
follow from the continuous mapping theorem. Finally, the map underlying
the estimator $\hat{D}_{n}$ defined in \prettyref{eq:ahat-1} inherits
the analyticity property by \prettyref{cor:derivative-estimator}.
Observing that $\psi_{D}\left(M\right)$ is smooth in $M$, we can
apply the continuous mapping theorem whence \prettyref{thm:consistency}\prettyref{enu:estimator}
follows.
\end{proof}

\subsection{Distribution of singular subspace estimators, \prettyref{cor:consistency-svd}}
\begin{proof}[Proof of \prettyref{cor:consistency-svd}.]
Note that right-singular vectors of $M$ are eigenvectors of $MM^{\trans}$.
For \prettyref{cor:consistency-svd}\prettyref{enu:B-svd-cons}, a
similar argument to the proof of \prettyref{thm:consistency}\prettyref{enu:B-cons}
given in above establishes the result. Similarly, \prettyref{cor:consistency-svd}\prettyref{enu:right-sg-cons}
follows from \prettyref{thm:consistency}\prettyref{enu:.right-cons}.
\end{proof}
Recall \prettyref{cor:folded-normal}, which states that
\[
\sqrt{n}\left(\hat{D}_{n}^{\trans}-D^{\trans}\right)=\sqrt{n}\upsilon_{\perp}^{\trans}\hat{R}_{n,J}\boldsymbol{S}^{-1}\left(\Delta_{n,\text{IJ}}\right)\hat{R}_{n,2,I}^{-1}+\sqrt{n}r\left(\hat{M}_{n}-M\right).
\]

\begin{proof}[Proof of \prettyref{cor:expansion-and-convergence-estimator}.]
 Apply \prettyref{cor:derivative-estimator} to the expression in
\prettyref{lem:ls-estimator-rewritten} and the result follows.
\end{proof}
Recall that \prettyref{assu:general-ass}\prettyref{enu:no-evs} was
sufficient for \prettyref{lem:open-and-smooth}\prettyref{enu:smooth-maps}.
To furnish necessity we have to show that in the following can assume
that invariant vectors are biorthogonal, i.e. that $R_{I}^{\trans}R_{I}=I_{q}$
without loss of generality. A corollary of this result is that the
eigenvectors do not need to be normalized for \prettyref{lem:open-and-smooth}\prettyref{enu:smooth-maps}
and \prettyref{lem:svd-smooth} to apply.
\begin{lem}
\label{lem:ls-estimator-rewritten}Let the estimator $\hat{D}_{n}^{\trans}=\hat{R}_{n,1,I}\hat{R}_{n,2,I}^{-1}$.
Then, the centered estimator
\begin{equation}
\left(\hat{D}_{n}-D\right)^{\trans}=\upsilon_{\perp}^{\trans}\hat{R}_{n,I}\hat{R}_{n,2,I}^{-1}+o_{p}\left(1\right).\label{eq:rewrittenestimator-1}
\end{equation}
\end{lem}
\begin{proof}[Proof of \prettyref{lem:ls-estimator-rewritten}.]
Add and subtract $\upsilon_{\perp}^{\trans}$ to obtain $\left(\left[\begin{smallmatrix}\hat{D}_{n}^{\trans} & I_{r}\end{smallmatrix}\right]\pm\left[\begin{smallmatrix}D^{\trans} & I_{r}\end{smallmatrix}\right]\right)\hat{R}_{n,I}=0$.
Therefore, 
\[
\left(\hat{D}_{n}^{\trans}-D^{\trans}\right)\hat{R}_{n,2,I}=\begin{bmatrix}\hat{D}_{n}^{\trans} & I_{r}\end{bmatrix}\left(\hat{R}_{n,I}-R_{I}R_{2,I}^{-1}\hat{R}_{n,2,I}\right),
\]
whence the result follows.
\end{proof}
To study the standardized and centered estimator, we can work with
the right-hand side of \prettyref{eq:rewrittenestimator-1}. For the
distributional limit, we have in analogy to \prettyref{lem:expansion-and-convergence},
\begin{cor}
\label{cor:expansion-and-convergence-estimator}We have the expansion
for the centered estimator\textup{
\[
\sqrt{n}\left(\hat{D}_{n}^{\trans}-D^{\trans}\right)=\sqrt{n}\upsilon_{\perp}^{\trans}\hat{R}_{n,J}\boldsymbol{S}^{-1}\left(\Delta_{n,\text{IJ}}\right)\hat{R}_{n,2,I}^{-1}+\sqrt{n}r\left(\hat{M}_{n}-M\right).
\]
}
\end{cor}
We note that \prettyref{lem:expansion-and-convergence}\prettyref{enu:remainder-term}
and \ref{enu:order} remain unchanged in the case of expansions for
$\hat{D}_{n}$ and the associated $t$-test statistic.

For the map induced by the normalized estimator $\hat{D}_{n}^{\trans}$,
we have in analogy to $\psi\left(M\right)$ defined in \prettyref{eq:inv-subspace-map},
the map $\psi_{D}\,:\,\set Q\mapsto\reals^{r\times q}$ with $\psi_{D}\left(M\right)\defeq\upsilon_{\perp}^{\trans}R_{I}\left(M\right)R_{2,I}^{-1}$.
Then, we have
\begin{cor}
\label{cor:derivative-estimator} The Fréchet derivative of $\psi_{D}$,
\[
\dot{\psi}_{D}\left(\hat{M}_{n}-M\right)=\upsilon_{\perp}^{\trans}R_{n,J}\boldsymbol{S}^{-1}\left(\hat{\Delta}_{n,\text{SL}}\right)L_{n,I}^{\trans}R_{n,2,I}^{-1}.
\]
\end{cor}
\begin{proof}[Proof of \prettyref{cor:derivative-estimator}.]
 In \prettyref{lem:derivative-1}, replace $L_{I}^{\trans}$ by $R_{2,I}^{-1}$
and the result follows.
\end{proof}
\begin{proof}[Proof of \prettyref{thm:waldstatdistrib}.]
Recall the map $\psi\left(M_{n}\right)$
\begin{equation}
\upsilon_{\perp}^{\trans}\hat{P}_{n,I}=\upsilon_{\perp}^{\trans}\left(\hat{P}_{n,I}-P_{I}\right)\label{eq:lefthalfwald}
\end{equation}
because $\upsilon_{\perp}^{\trans}P_{I}=0$ under the maintained hypothesis.
To expand the right-hand side of \prettyref{eq:lefthalfwald}, we
introduce the derivatives in Similarly, \prettyref{thm:waldstatdistrib}
follows from \prettyref{lem:deriv-conv}\prettyref{enu:proje} where
the degrees of freedom follow from \prettyref{lem:deriv-conv}\prettyref{enu:deg}.
The consistency of the covariance matrix and right and left eigenvector
column space estimates that is necessary for the convergences to hold
jointly follows from \prettyref{thm:consistency}.
\end{proof}
\begin{proof}[Proof of \prettyref{thm:tdist}.]
We now turn to proving \prettyref{thm:tdist}. By \prettyref{thm:estdist}
and \prettyref{assu:general-ass}, 
\[
\sqrt{n}\left(\hat{B}_{D}\hat{\Omega}_{M}\hat{B}_{D}^{\trans}\right)^{-1/2}\left(\vek\hat{D}_{n}^{\trans}-\vek D^{\trans}\right)\wkc N(0,I_{rq\times rq}).
\]
Let $e_{s,i}\in\reals^{s}$ denote a vector with zero everywhere except
for a $1$ in the $s$th position so that $d_{ij}\defeq e_{i,q}^{\trans}D^{\trans}e_{j,r}$
is the entry of $D^{\trans}$ in row $i$ and column $j$, which corresponds
to the $ij$th element of $\vek\hat{D}_{n}^{\trans}-\vek D_{0}^{\trans}$.
To find its associated standard error, we let $\sigma_{ij}^{2}$ denote
the $ij$th diagonal entry of $\hat{B}_{D}\hat{\Omega}_{M}\hat{B}_{D}^{\trans}$.
Then, we define as the standard error for $d_{ij}$, $\hat{\sigma}_{ij}:=\hat{\sigma}_{ij}^{2}/\sqrt{n}$
and find that $\frac{\hat{d}_{ij}-d_{ij}}{\hat{\sigma}_{ij}}\wkc N\left(0,1\right)$.
Note that $\hat{B}_{D}\hat{\Omega}_{M}\hat{B}_{D}^{\trans}$ is always
real by construction in analogy to the covariance matrix of the Wald
test statistic.

\citet{Tsagris_2014} give a proof of \prettyref{cor:folded-normal}.
This result allows us to construct one-sided confidence intervals
by folding the normal distribution over along the $y$-axis. An application
making use of these intervals appears in .
\end{proof}

\section{Jacobians and results for \prettyref{sec:jacobians-ho-perturbations}:
\nameref{sec:jacobians-ho-perturbations}\label{app:Jacobians}}

This appendix develops perturbation theory to obtain Jacobians of
invariant and singular subspace maps to obtain the distributions of
test statistics and the identification of associated basis vectors.
The aim is to outline the details of how we use perturbative arguments
to apply the delta method to invariant and singular subspace maps.

\subsection{\label{app:Additional-details-on-invar}Jacobian of invariant subspaces}

We introduce the off-diagonal element of $\Lambda$ corresponding
to $\Lambda_{JI}$ as 
\begin{equation}
\Delta_{\text{JI}}\defeq L_{J}^{\trans}\left(\tilde{M}-M\right)R_{I},\label{eq:argument-1}
\end{equation}
where $\tilde{M}=M$ implies $\Lambda_{JI}=0$. An invariant subspace
decomposition obtains from

\begin{equation}
M=R_{I}\left(\Lambda_{I}L_{I}^{\trans}+\Lambda_{IJ}L_{J}^{\trans}\right)+R_{J}\Lambda_{J}L_{J}^{\trans},\label{eq:gen-inv-subsp-decomp}
\end{equation}
where we use conjugate transposition on the matrices $L_{I}$ and
$L_{J}$ to allow for the fact that some eigenvectors may be complex.
Setting $\Lambda_{IJ}=0$ implements the familiar eigendecomposition.
The map that delivers the column span of $R_{I}$, which In this case,
based on \prettyref{eq:gen-inv-subsp-decomp}, we have the familiar
relation 
\begin{equation}
MR_{I}=R_{I}\Lambda_{I}.\label{eq:inv-subspace-def}
\end{equation}
Perturbations to $M$ affect the column span of $R_{I}$ via $\Delta_{\text{JI}}$,
where \prettyref{app:Additional-details-on-invar} contains additional
details.

\noindent To build a test statistic, we introduce 
\begin{equation}
P_{I}\defeq R_{I}L_{I}^{\trans},\label{eq:definingplu-1}
\end{equation}
which is a skew projection for eigenspaces ($\Lambda_{IJ}=0$) and
a regular projection for symmetric $M$. If the roots in the set $\eigs_{I}$
are closed\emph{ }under conjugation, then we may take $R_{I}$ and
$L_{I}$ to be real without loss of generality. However, even if $M$
itself is a real matrix, we may encounter estimation error that may
induce complex invariant vectors and, in the case of an eigendecomposition,
complex eigenvalues. Naturally, if invariant vectors of interest are
real, the imaginary parts of their estimators converge to zero in
probability. Under the hypothesis $H_{0}:\upsilon\in\risof$, $\upsilon_{\perp}^{\trans}P_{I}=0$.
To write out the first-order derivative of $\psi$, denoted by $\dot{\psi}$,
we introduce the operator $\boldsymbol{S}\,:\,\mathbb{R}^{r\times q}\rightarrow\mathbb{R}^{r\times q}$
such that 
\begin{equation}
\boldsymbol{S}\left(Q\right)\defeq Q\Lambda_{I}-\Lambda_{J}Q,\label{eq:lin-op}
\end{equation}
for a full-rank $Q\in\mathbb{R}^{r\times q}$. 

\noindent Then, the first-order term reads as
\begin{equation}
\dot{\psi}\left(M;\Delta_{\text{JI}}\right)=\upsilon_{\perp}^{\trans}R_{J}\boldsymbol{S}^{-1}\left(\Delta_{\text{JI}}\right)L_{I}^{\trans}.\label{eq:frechet}
\end{equation}

The subspace $\col R_{I}$ is invariant iff $\Lambda$ in
\begin{align}
R=\begin{bmatrix}R_{I}\,R_{J}\end{bmatrix}\,L=\begin{bmatrix}L_{I}\,L_{J}\end{bmatrix}\,\Lambda & =\begin{bmatrix}\Lambda_{I} & \Lambda_{IJ}\\
0 & \Lambda_{J}
\end{bmatrix},\label{eq:BlocksOfR}
\end{align}
is block upper-triangular and we have used \prettyref{assu:invar-sub-ass}
to partition the matrices according to their column spans. If $\Lambda_{JI}\neq0$
we would have
\begin{equation}
MR_{I}=R_{I}\Lambda_{I}+R_{J}\Lambda_{JI}.\label{eq:right-invar}
\end{equation}
Therefore, we can interpret $\norm{\Lambda_{\text{JI}}}$ as a measure
of how far $\col R_{I}$ is from being an invariant subspace. Regarding
the normalization \eqref{eq:eigenvector-normalization}, it is equivalent
to requiring $\upsilon_{\perp}^{\trans}N\upsilon_{\perp}=I_{q}$ for
normalized coefficients
\[
N=\left[\begin{smallmatrix}\left(\upsilon_{1}^{-\trans}\upsilon_{2}^{\trans}\upsilon_{2}\upsilon_{1}^{-1}\right)^{-1} & 0\\
0 & \left(\upsilon_{2}\left(\upsilon_{2}^{\trans}\upsilon_{2}\right)^{-1}\upsilon_{1}^{\trans}\upsilon_{1}\left(\upsilon_{2}^{\trans}\upsilon_{2}\right)^{-1}\upsilon_{2}^{\trans}-\upsilon_{2}\left(\upsilon_{2}^{\trans}\upsilon_{2}\right)^{-3}\upsilon_{2}^{\trans}\right)\upsilon_{2}\upsilon_{1}^{-1}
\end{smallmatrix}\right].
\]
\begin{align*}
 & \upsilon_{\perp}^{\trans}N\upsilon_{\perp}\\
 & =\upsilon_{1}^{\trans}\left(\upsilon_{1}^{-\trans}\upsilon_{2}^{\trans}\upsilon_{2}\upsilon_{1}^{-1}\right)^{-1}\upsilon_{1}\\
 & +\upsilon_{1}^{-\trans}\upsilon_{2}^{\trans}\left(\left(\upsilon_{2}\left(\upsilon_{2}^{\trans}\upsilon_{2}\right)^{-1}\upsilon_{1}^{\trans}\upsilon_{1}\left(\upsilon_{2}^{\trans}\upsilon_{2}\right)^{-1}\upsilon_{2}^{\trans}-\upsilon_{2}\left(\upsilon_{2}^{\trans}\upsilon_{2}\right)^{-3}\upsilon_{2}^{\trans}\right)\right)\upsilon_{2}\upsilon_{1}^{-1}\\
 & =I_{q}
\end{align*}
when we make the replacement $\upsilon_{1}\leftarrow I_{q}$ and $\upsilon_{2}\leftarrow\upsilon_{2}\upsilon_{1}^{-1}$.
Letting $\Omega_{M}$ be as in \prettyref{eq:homoskedastic-cols},
we have in the case of homoskedastic columns of $\hat{M}$ the alternative
normalization $\upsilon_{\perp}^{\trans}\Omega_{M}\upsilon_{\perp}=I_{q}$
where \citet{Anderson2010359} discusses potential drawbacks of this
normalization. For an approach that avoids normalizations altogether,
see \citet{silin2020hypothesis} for symmetric matrices.

\subsection{\label{app:perturb-to-deriv}Relationship between perturbation expansions
and derivatives}

To identify the correspondence between a perturbation expansion and
the appearance of derivatives, it is instructive to consider the following:
\begin{example}
\label{exa:perturb-deriv}Consider a continuous function $f\left(x\right)$
and a small $\epsilon$, that perturbs $x$, such that $\tilde{x}\defeq x+\epsilon$.
Recall the difference quotient 
\[
d\left(\tilde{x};x\right)\defeq\frac{f\left(\tilde{x}\right)-f\left(x\right)}{\tilde{x}-x}
\]
and observe that $\underset{x\goesto\tilde{x}}{\lim}\,d\left(\tilde{x};x\right)=f^{\trans}\left(x\right)$.
Now, let $f\left(x\right)=x^{2}$. Define a perturbation expansion
for $f$ with first-order perturbation term $\dot{f}$ like so
\begin{align*}
f\left(\tilde{x}\right) & =f\left(x\right)+\dot{f}\left(x\right)+O\left(\epsilon\right)\\
x^{2}+2\epsilon x+\epsilon^{2} & =x^{2}+2x\epsilon+O\left(\epsilon\right)
\end{align*}
and recognize that the relationship, as $\tilde{x}\goesto x$ and
$\epsilon\goesto0$, between the first-order perturbation expansion
term $\dot{f}$ and the derivative $f^{\trans}$ satisfies 
\begin{equation}
\dot{f}=f^{\trans}\epsilon,\label{eq:perturb-deriv}
\end{equation}
where we can identify the derivative from a perturbation expansion
without taking the limit. \citet[Thm. 11]{magnus2019matrix} extends
this relationship to matrix-valued functions. In non-scalar cases,
it is important that the perturbation $\epsilon$ in \prettyref{eq:perturb-deriv}
appears at the end of the expression, which we achieve in practice
by vectorizing matrix expressions and applying the identity $\vek YXZ=\left(Z^{\trans}\otimes Y\right)\vek X$.
\end{example}
We recall the maps delivering invariant and singular subspaces, $\psi\left(M;\upsilon\right)$
and $\Psi\left(M;\upsilon\right)$, respectively, which carry information
about invariant and singular subspaces. Further recall that, if $\upsilon$
belongs to a subspace of interest, i.e. $\upsilon\in\risof$ iff $\psi\left(M;\upsilon\right)=0$
and $\Psi\left(M;\upsilon\right)=0$. Exploiting the information thus
expressed, we can build inference methods for subspaces based on these
maps if we can derive the distribution of the sample analogues of
$\psi$ and $\Psi$ from that of $\sqrt{n}\vek\hat{M}_{n}$. The critical
ingredient for the extraction of the distribution are the Jacobians
of these maps and appear in our central results in \prettyref{lem:jacobian}.
We shall apply the delta method to transform the covariance matrix
of $\sqrt{n}\,\vek\hat{M}_{n}$ into those of $\psi$ and $\Psi$.
To derive the Jacobians, we use expansions of $\psi$ and $\Psi$
that originate from a perturbed version of $M$ analogously to \prettyref{exa:perturb-deriv},
denoted by $\tilde{M}$. In the spirit of \prettyref{exa:perturb-deriv},
we consider the perturbation explicitly, $E\defeq\tilde{M}-M$, as
the deviation from the true $M$. A typical perturbation expansion
of this type appears as $\psi\left(\tilde{M}\right)$ about $\psi\left(M\right)$,
so that 
\begin{equation}
\psi\left(\tilde{M}\right)=\psi\left(M\right)+\dot{\psi}\left(M-\tilde{M}\right)+R_{M}\left(E\right),\label{eq:perturb-exp}
\end{equation}
where we could always (naively) write $R_{M}\left(E\right)=O\left(\norm{M-\tilde{M}}\right)$.
An exact analogue with $\dot{\Psi}$ in lieu of $\dot{\psi}$ holds
for $\Psi$. Importantly, we shall argue that the first-order term
indeed corresponds to the Gateaux derivative of $\psi$ ($\Psi$).

Following \citet{magnus2019matrix}, we define the Jacobian matrix
in the direction of $E$ of a vectorized, matrix-valued function $F\left(X\right):\reals^{m\times l}\goesto\reals^{r\times q}$
via first differentials as the matrix $J\in\reals^{ml\times rq}$
satisfying
\[
\text{d}\vek F\left(X;E\right)=J\vek E,
\]
whence \citet[Thm. 11]{magnus2019matrix} establishes that $J$ is
indeed unique. In such case, 
\[
\psi\left(\tilde{M}\right)-\psi\left(M\right)\goesto\diff\psi
\]
as $\norm E\goesto0$ so that, by \citep[Remark 4.2]{Sun1991}, \prettyref{eq:perturb-exp}
turns into the definition of the Gateaux derivative. Formally, $R_{M}\left(E\right)=o\left(\norm E\right)$
is required to establish differentiability. Therefore, $\vek\dot{\psi}\left(M-\tilde{M}\right)$
will allow us to identify the Jacobian $J$ of the transformation
$\psi$. \prettyref{lem:perturb-diff} in \prettyref{app:Jacobians}
makes the preceding discussion formal and finds a tighter bound on
the remainder term $R_{M}\left(E\right)$ in \prettyref{eq:perturb-exp},
while \prettyref{subsec:expansion-psi} sketch the steps required
to derive the Jacobians for invariant and singular subspace maps.

By \prettyref{lem:open-and-smooth}, the map $\psi$ inherits the
infinite differentiability of the maps $R_{I}\left(M\right)$ and
$\Lambda_{I}\left(M\right)$, which enables application of the Taylor
expansion for $M\in\mathscr{M}$ and $E$ as in \prettyref{enu:invariant}
\begin{equation}
\psi\left(M+E\right)=\psi\left(M\right)+\diff\psi\left(M,E\right)+r\left(E\right),\label{eq:Taylor}
\end{equation}
where $r\left(.\right)$ is a remainder term. The fundamental result
is that we can identify the Jacobian $B$ of the map $\psi$ from
the differential in \prettyref{eq:Taylor} based on 
\begin{equation}
\diff\vek\psi\left(M,E\right)=B\vek E.\label{eq:Jacobian}
\end{equation}
\prettyref{lem:jacobian}\prettyref{enu:invariant} provides us with
an expression of the Jacobian of $\psi$. Formally, we have
\begin{lem}
\label{lem:perturb-diff}Let the Taylor expansion of $\psi$ be as
in \prettyref{eq:Taylor} and the Jacobian $B$ as in \prettyref{eq:jacobian}.
Then, the first-order perturbation term defined in \prettyref{eq:frechet},
\[
\dot{\psi}\left(M,E\right)=\diff\psi\left(M,E\right),
\]
which in vectorized form $\vek\dot{\psi}\left(M,E\right)=B\vek E$
identifies the Jacobian matrix $B$.
\end{lem}
\begin{proof}[Proof of \prettyref{lem:perturb-diff}.]
Denote by $\norm .$ Frobenius norms. Recognize that (2) in \citet[Ch. 5, Sec. 15, Def. 3]{magnus2019matrix}
is a perturbation expansion in the form of \prettyref{eq:perturb-exp}
and is equal to the definition of the Gateaux derivative before taking
limits. By (3) in \citet[Ch. 5, Sec. 15, Def. 3]{magnus2019matrix},
we need to have that as $\norm E\goesto0$, $\frac{R_{M}\left(E\right)}{\norm E}\goesto0$,
where $R_{M}\left(E\right)$ is defined in \prettyref{eq:perturb-exp}.
To establish this result, we must estimate (not in the statistical
sense) the orders of magnitude of $R_{M}\left(E\right)$ relative
to that of $\norm E$. First, observe that $\norm{R_{M}\left(E\right)}\goesto0$
iff $R_{M}\left(E\right)\goesto0$. The Gateaux argument is that $\smlnorm E$
is fixed while $t\goesto0$ so that we establish differentiability
from all directions. By \citet[Thm. 3.1]{Sun1991}, $\norm{R_{M}\left(E\right)}\leq\frac{\norm E}{\delta}\abs t+O\left(\abs t^{2}\right)$
for some parameter $t\in\complex$ which, using \citet[Remark 4.3]{Sun1991}
allows choosing $E$ such that $E\left(t\right)=tE$ so that $\tilde{M}\goesto M$
as $E\goesto0$ or $t\goesto0$, which permits us to write $\abs t\propto\norm E$
and showing Fréchet differentiability $\norm{R_{M}\left(E\right)}/\norm E\goesto0$
suffices. In other words, we have $\frac{1}{\delta}=\norm{\left[\left(\Lambda{}_{I}^{\trans}\otimes I_{J}\right)-\left(I_{I}\otimes\Lambda_{J}\right)\right]^{-1}}<\infty$
by \prettyref{assu:invar-sub-ass}\prettyref{enu:no-evs}. Thus, $\norm{R_{M}\left(E\right)}/\norm E\leq\frac{\norm E^{2}}{\delta}+O\left(\norm E\right)$,
whence we can see $\frac{\norm E^{2}}{\delta}+O\left(\norm E\right)\goesto0$
as $\norm E\goesto0$.
\end{proof}
\prettyref{lem:perturb-diff} allows us to treat perturbation terms
synonymously with derivatives and the critical ingredient in the proof
was that $\frac{1}{\delta}=\norm{\left[\left(\Lambda{}_{I}^{\trans}\otimes I_{J}\right)-\left(I_{I}\otimes\Lambda_{J}\right)\right]^{-1}}<\infty$
by \prettyref{assu:invar-sub-ass}\prettyref{enu:no-evs}. These derivatives
apply equally to singular subspaces as we can always represent them
as invariant subspaces of symmetric matrices. Therefore, the next
result's perturbation terms lead directly to derivatives of invariant
subspace maps. Recall the block-diagonalization of the perturbing
matrix $E=\tilde{M}-M$, \eqref{eq:argument-1-1-1}, which we restate
as 

\begin{equation}
\begin{bmatrix}\Delta_{I}\\
\Delta_{\text{JI}} & \Delta_{J}
\end{bmatrix}=\begin{bmatrix}L_{I}^{\trans}ER_{I}\\
L_{J}^{\trans}ER_{I} & L_{J}^{\trans}ER_{J}
\end{bmatrix}.\label{eq:argument-1-1}
\end{equation}

\begin{lem}
\label{lem:derivative-1}Suppose \prettyref{assu:general-ass} holds.
Then, the first and second order Gateaux derivatives of $\psi\left(M,E\right)$
at $M$ are equal to their perturbation expansion terms, i.e.
\begin{align}
\dot{\psi}\left(X,E\right)\vert_{X=M}= & \upsilon_{\perp}^{\trans}R_{J}\boldsymbol{S}^{-1}\left(\Delta_{\text{JI}}\right)L_{I}^{\trans}.\label{eq:firstderivative-1}\\
\ddot{\psi}\left(X,E\right)\vert_{X=M}= & 2\upsilon_{\perp}^{\trans}R_{J}\boldsymbol{S}^{-1}\left(\left(\Delta_{J}\boldsymbol{S}^{-1}\left(\Delta_{\text{JI}}\right)-\boldsymbol{S}^{-1}\left(\Delta_{\text{JI}}\right)\Delta_{I}\right)\right)L_{I}^{\trans}.\label{eq:secondderivative-1}
\end{align}
\end{lem}
Indeed, we now apply \prettyref{lem:perturb-diff} to \prettyref{lem:derivative-1}
to obtain the desired Jacobians, which the next result establishes:
\begin{lem}
\label{lem:neudeck-deriv}Suppose \prettyref{assu:general-ass} holds.
Then, the first order differential form of $\psi\left(M,E\right)$
at $M$ is
\begin{align}
\vek\diff\psi\left(M\right) & =\vek\upsilon_{\perp}^{\trans}R_{J}\boldsymbol{S}^{-1}\left(\Delta_{\text{JI}}\right)L_{I}^{\trans},\label{eq:diff-eq-perturb}
\end{align}
which lets us identify the Jacobian matrix $B$ via 
\begin{equation}
\vek\diff\psi\left(M;E\right)=\underset{B}{\underbrace{\left(L_{I}\otimes\upsilon_{\perp}^{\trans}R_{J}\right)\left\{ \left(\Lambda_{I}^{\trans}\otimes I_{J}\right)-\left(I_{I}\otimes\Lambda_{J}\right)\right\} ^{-1}\left(R_{I}^{\trans}\otimes L_{J}^{\trans}\right)}}\vek E\label{eq:jacob-expression}
\end{equation}
\end{lem}
\begin{proof}[Proof of \prettyref{lem:neudeck-deriv}.]
Expression \prettyref{eq:diff-eq-perturb} follows from applying
\prettyref{lem:perturb-diff} to \prettyref{lem:derivative-1} while
\prettyref{eq:jacob-expression} follows from \prettyref{lem:jacobian}\prettyref{enu:invariant}.
\end{proof}
The Gateaux derivatives in \prettyref{lem:derivative-1} depend on
the perturbation $E$, whose limit is the differential. An alternative
way to obtain such derivatives is to apply the rules of calculus and
to invoke \citet[Thm. 11]{magnus2019matrix}, as \citet{duffy2020cointegrated}
do. 

For easier reference, we state \citet[Thm. 2.1]{Sun1991}. For $X_{1}$,
such that $\spn X_{1}=\spn R_{I}$, we state.: 
\begin{thm}
\label{thm:expthm-1}There exists a unique $q$-dimensional invariant
subspace $\spn X_{1}\left(t\right)$ of $M\left(t\right)$ $\left(t\in\complex\right)$
such that $\spn X_{1}\left(0\right)=\spn X_{1}$ and the basis vectors
(columns of $X_{1}$) may be defined to be analytic functions of $t$
in some neighborhood of the origin of $\complex$. Further, the analytic
matrix-valued function $X_{1}\left(t\right)$ has the second order
perturbation expansion 
\begin{align}
X_{1}\left(t\right) & =X_{1}\label{eq:perturb-1}\\
 & +X_{2}\boldsymbol{S}^{-1}\left(\Delta_{\text{JI}}\right)t\nonumber \\
 & +X_{2}\boldsymbol{S}^{-1}\left(\Delta_{J}\boldsymbol{S}^{-1}\left(\Delta_{\text{JI}}\right)-\boldsymbol{S}^{-1}\left(\Delta_{\text{JI}}\right)\Delta_{I}-\boldsymbol{S}^{-1}\left(\Delta_{\text{JI}}\right)\Lambda_{IJ}\boldsymbol{S}^{-1}\left(\Delta_{\text{JI}}\right)\right)t^{2}\nonumber \\
 & +O\left(t^{3}\right).\nonumber 
\end{align}
\end{thm}

The perturbation expansion in \prettyref{eq:perturb-1} is valid for
basis vectors $X_{1}$ which may be unnormalized and are not orthogonal
in general. Without altering $\col R_{I}$, we can find orthonormal
basis vectors by
\begin{lem}
\label{lem:geometry-matching}Let the columns of $R$ span an invariant
subspace of $M$ such that $MR=R\Lambda$ for $\Lambda$ upper triangular
and $MR_{I}=R_{I}\Lambda_{I}$ as in \prettyref{assu:invar-sub-ass}.
Then, the following are equivalent:
\begin{enumerate}
\item There exists an orthonormal basis of vectors $s_{1},\dots,s_{q}$
for $\col R_{I}$ such that $\col R_{I}=\spn\left\{ s_{1},\dots,s_{q}\right\} $.
\item There exists a bijection to recover right-invariant vectors from the
orthonormal basis vectors $\left\{ s_{1},\dots,s_{q}\right\} \rightleftarrows\left\{ r_{1,}\dots,r_{q}\right\} $.
\end{enumerate}
\end{lem}
\prettyref{lem:geometry-matching} allows us to treat $R_{I}$ as
though they are orthogonal, the basis vectors spanning our invariant
subspace of interest, $\col R_{I}$, are not in general biorthogonal,
so a relation that left and right eigenvectors fulfill, $l_{i}^{\trans}r_{i}=\delta_{ij}$
does not in general hold for right or left eigenvectors alone, i.e.
$r_{i}^{\trans}r_{j}=\delta_{ij}$ applies only to symmetric matrices
when left and right eigenvectors coincide.
\begin{proof}[Proof of \prettyref{lem:geometry-matching}.]
Employ the QR factorization which exists for every real and complex
matrix whereby $R_{I}=SU$ with $S^{\trans}S=I_{q}$ and $U$ upper
triangular. Because $\rank U=q$, $\col R_{I}=\col S$ so that the
invariant subspace of interest that we wish to conduct inference on
can always be identified. Finally, $MR=R\Lambda$ implies 
\[
MS=S\underbrace{U\Lambda U^{-1}}_{\text{upper triangular}},
\]
so that $\col S=\col R_{I}$ is indeed an invariant subspace of $M$.
To see why $U\Lambda U^{-1}$ is upper triangular, we need only show
that $U^{-1}$ is upper triangular whenever $U$ is. Represent $U=D\left(I_{q}+K\right)$
where $D$ is diagonal and $K$ is strictly upper triangular. Then,
$\left(I_{q}+K\right)^{-1}=I_{q}-K+K^{2}-\dots+\left(-1\right)^{q-1}K^{q-1}$
which has only strictly upper triangular summands. Moreover, $K^{q}$
and higher powers are zero by the Cayley-Hamilton theorem. Then, the
columns of $S$ span the invariant subspace of $M$ associated with
roots in $R_{I}$. Finally, we have the relation for the columns of
$R_{I}$, $r_{i}=\sum_{k=1}^{i}u_{ki}s_{k}$, where $s_{k}$ are the
columns of $S$ and $u_{ki}$ is the $\left(k,i\right)$th entry of
$U$. This mapping is bijective because $U$ is of full rank so that
a reverse mapping can be achieved by swapping $r$ and $s$ and replacing
$u_{ki}$ with the entries of $U^{-1}$.
\end{proof}
\begin{proof}[Proof of \prettyref{lem:derivative-1}]

We wish to construct the first and second derivatives of the map 
\[
\psi\left(M\left(t\right)\right)=\upsilon_{\perp}^{\trans}R_{I}\left(t\right)L_{I}^{\trans}\left(t\right),
\]
where we have made the dependence on $t$ explicit and $M\left(0\right)=M$.
Then, by \prettyref{lem:perturb-diff} \citet[Remark 4.2]{Sun1991},
the Fréchet and Gateaux derivatives coincide and we can write 
\[
\dot{\psi}_{R}\left(\hat{M}_{n}-M\right)=\upsilon_{\perp}^{\trans}\dot{P}_{I}\left(t\right)|_{t=0},
\]
which by the product rule implies 
\[
\upsilon_{\perp}^{\trans}\dot{P}\left(0\right)=\upsilon_{\perp}^{\trans}\dot{R}_{I}\left(0\right)L_{I}^{\trans}+\upsilon_{\perp}^{\trans}R_{I}\dot{L}_{I}^{\trans}\left(0\right)=\upsilon_{\perp}^{\trans}\dot{R}_{I}\left(0\right)L_{I}^{\trans}
\]
because $\upsilon_{\perp}^{\trans}R_{I}$$\left(0\right)=\upsilon_{\perp}^{\trans}R_{I}=0.$
Therefore, we need to evaluate the derivative of basis vectors collected
in $R_{I}$. Note that normally, right-invariant vectors do not satisfy
$R_{I}^{\trans}R_{I}=I_{q}$ and $R_{J}^{\trans}R_{J}=I_{p-q}$ although
the perturbation theory requires this relationship. By \prettyref{lem:geometry-matching},
there always exist bijections $g_{1}$ and $g_{2}$ between $R_{I}$
and $X_{1}$ as well as $R_{J}$ and $X_{2}$ so that we can apply
Theorem 2.1 in \citet{Sun1991}. Therefore, let $X_{1}=g_{1}\left(R_{I}\right)\in\reals^{p\times q}$
with $\rank X_{1}=q$, $X_{1}^{\trans}X_{1}=I_{q}$ and $X_{2}=g_{2}\left(R_{J}\right)\in\reals^{p\times p-q}$
for $MX_{1}=X_{1}M_{1}$ for some matrix $M_{1}\in\reals^{q\times q}$.
Then $\spn X_{1}$ is an invariant subspace of $M\in\reals^{p\times p}$
if and only if there exists a non-singular matrix $X=\left[\begin{smallmatrix}X_{1} & X_{2}\end{smallmatrix}\right]\in\reals^{p\times p}$
with $X_{2}^{\trans}X_{2}=I_{p-q}$ such that
\[
X^{-1}MX=\begin{bmatrix}M_{11} & M_{12}\\
0 & M_{22}
\end{bmatrix},\,M_{11}\in\reals^{q\times q}.
\]
For the purposes of this study, we are interested in the span of the
vectors that are orthogonal to the columns of $X_{1}.$ Recall the
operator \prettyref{eq:lin-op}, 
\[
\boldsymbol{S}\left(Q\right)=QM_{11}-M_{22}Q
\]
and $Q\in\reals^{\left(p-q\right)\times q}.$ Then, a necessary condition
for Theorem 2.1 in \citet{Sun1991} to apply is that $M_{11}$ and
$M_{22}$ do not share any eigenvalues as stipulated in \prettyref{assu:invar-sub-ass}\prettyref{enu:no-evs}.
\end{proof}
\begin{proof}[Proof of \prettyref{lem:expansion-and-convergence}.]
~
\begin{enumerate}
\item We Taylor-expand 
\begin{align*}
\sqrt{n}\upsilon_{\perp}^{\trans}\hat{P}_{n,I} & =\sqrt{n}\psi\left(\hat{M}_{n}\right)\\
 & =_{\left(2\right)}\sqrt{n}\left(\psi\left(\hat{M}_{n}\right)-\psi\left(M\right)\right)\\
 & =\sqrt{n}\dot{\psi}\left(\hat{M}_{n}-M\right)+\sqrt{n}r\left(\hat{M}_{n}-M\right)\\
 & =_{\left(3\right)}\sqrt{n}\upsilon_{\perp}^{\trans}R_{n,J}\boldsymbol{S}^{-1}\left(\hat{\Delta}_{n,\text{IJ}}\right)L_{n,I}^{\trans}+\sqrt{n}r\left(\hat{M}_{n}-M\right)
\end{align*}
where $=_{\left(2\right)}$ follows from application of the null hypothesis
and $=_{\left(3\right)}$ follows from \prettyref{lem:derivative-1}. 
\item We obtain a remainder term 
\[
r\left(\hat{M}-M\right)\defeq\psi\left(\hat{M}\right)-\psi\left(M\right)-\dot{\psi}\left(\hat{M}-M\right).
\]
To ensure that the first order term dominates the remainder, we need
to verify that standardization by $\sqrt{n}$ does not cause dominant
second order terms. By \prettyref{assu:general-ass}\prettyref{enu:asymptotically-normal-estimator},
$\smlnorm{\hat{M}-M}=O_{p}(n^{-1/2})$.  For simplicity, we can treat
$(\hat{M}_{n}-M)$ as a deterministic sequence and then apply \citet[Lemma 2.12, Ch. 2]{book:770810}
to obtain the relevant statistical result. First, recognize that $\boldsymbol{S}^{-1}\left(a_{n}M\right)=a_{n}\boldsymbol{S}^{-1}\left(M\right)$
for any scalar sequence $a_{n}$ and admissible argument $M$. Then,
rewrite $r$ as the sum of the next higher order term and another
unspecified remainder term of known order via application of \prettyref{lem:derivative-1}.
Let $V$ be such that 
\[
\boldsymbol{S}^{-1}\left(L_{J}^{\trans}\left(M-\hat{M}_{n}\right)R_{I}\right)=V
\]
or
\[
VM_{11}-M_{22}V=L_{J}^{\trans}\left(M-\hat{M}_{n}\right)R_{I}.
\]
Then, $\smlnorm V=O(n^{-1/2})$ and $\smlnorm{\hat{\Delta}_{n,I}}=O(n^{-1/2})$
so that
\[
r\left(\hat{M}_{n}-M\right)=2\upsilon_{\perp}^{\trans}R_{J}\boldsymbol{S}^{-1}\left(\left(\hat{\Delta}_{n,J}V-V\hat{\Delta}_{n,I}\right)-V\Lambda_{IJ}V\right)L_{I}^{\trans}+r_{2}\left(\hat{M}_{n}-M\right)
\]
implies $n^{1/2}r(\hat{M}_{n}-M)=O(n^{-1/2})$ by \prettyref{lem:expansion-and-convergence}\prettyref{enu:order}
and \prettyref{lem:derivative-1}. By the expansion in \citet{Sun1991},
each higher order increases a power of the norm of $V$ so that $r_{2}$
is likewise negligible. Hence, a normalization by $\sqrt{n}$ is innocuous
and does not affect convergence.
\item For a consistent and asymptotically normal least-squares estimator
\prettyref{assu:general-ass}\prettyref{enu:asymptotically-normal-estimator}
and the definition in \prettyref{eq:argument-1-1}, the result follows
by block-diagonalizing $\hat{M}_{n}-M$ via 
\[
\hat{L}_{I}^{\trans}\left(\hat{M}_{n}-M\right)\hat{R}_{I}.
\]
\end{enumerate}
\end{proof}
We have seen how perturbation expansion terms correspond to derivatives.
The following proof shows how we bound the perturbation expansions
to higher orders.

\subsection{Additional results for higher-order Davis-Kahan bounds}

\label{app:inequality-block-diagonal}

Let the norm of the Jacobian $B$ be $\beta_{2}\defeq\norm B_{\text{F}}$.
First, note that the objects of interest are the spans of $\upsilon_{\perp}$,
$R_{I}$, $L_{I}$ so that they can always be orthogonalised using
Gram-Schmidt for example. Hence, we apply \citet[Lemma A1]{demetrius}
to find that $\smlnorm{\upsilon_{\perp}^{\trans}\hat{R}_{I}L_{I}^{\trans}}\leq\smlnorm{\hat{R}_{I}}$.
Then, using \eqref{eq:jacobian} and the preceding, we define 
\begin{align*}
\beta_{2} & \defeq\norm{L_{J}^{\trans}\left(A-\hat{A}\right)R_{I}}_{\text{F}}\leq\norm{A-\hat{A}}_{\text{F}}
\end{align*}
which implies that $\norm B_{\text{F}}\leq\norm{\left\{ \left(\Lambda_{I}^{\trans}\otimes I_{J}\right)-\left(I_{I}\otimes\Lambda_{J}\right)\right\} ^{-1}}_{\text{F}}$,
which for singleton sets $I=\left\{ 1\right\} $ and $J=\left\{ 2\right\} $
implies that $\norm B_{\text{F}}\leq\left(\lambda_{1}-\lambda_{2}\right)^{-1}$
so that \eqref{eq:davis-kahan-disguise} obtains.
\begin{proof}[Proof of \prettyref{thm:higher-order-dk}.]
Recursive substitution in \citet[Thm. 3.1]{Sun1991} yields the formula
for $m\geq3$,
\begin{align*}
\gamma_{m} & =\left(s\beta\right)^{m-2}\gamma_{2}+\sum_{k=0}^{m-3}\iota_{m,k}\\
\iota_{m,k} & =s^{k+1}\beta^{k}\beta_{1}\sum_{j=1}^{m-2-k}\gamma_{m-\left(k+1\right)-j}\gamma_{j}+s^{k+1}\beta^{k}\alpha_{12}\sum_{j=1}^{m-1-k}\gamma_{m-j-k}\gamma_{j}.
\end{align*}
Substituting $a\defeq\max\left\{ \beta,\beta_{1}\right\} $ for all
perturbation bounds yields a sequence that majorizes $\iota_{m,k}$,
i.e. we obtain 
\[
\iota_{m,k}\leq\left(\left(as\right)^{k+1}\sum_{j=1}^{m-2-k}\gamma_{m-\left(k+1\right)-j}\gamma_{j}+a^{k}s^{k+1}\alpha_{12}\sum_{j=1}^{m-1-k}\gamma_{m-k-j}\gamma_{j}\right).
\]
Now, the highest orders of the summands are $m-\left(k+1\right)$
and $m-k$ for $k=0,1,\dots,m-3$, which we prove by induction. Recall
$L(\gamma_{k+1})=\gamma_{1}\gamma_{k}$ and that all equalities hier
are modulo lower-order terms. We wish to show that $\gamma_{j}\gamma_{s}=L\left(\gamma_{j+s}\right)$
by repeatedly factoring out $\gamma_{1}$. A direct computation yields
$\gamma_{1}\gamma_{2}=L\left(\gamma_{3}\right)$. Assume $\gamma_{n}\gamma_{n+1}=L\left(\gamma_{2n+1}\right)$
and want to show that $\gamma_{n+1}\gamma_{n+2}=L\left(\gamma_{2n+3}\right)$.
By the induction hypothesis, $L\left(\gamma_{n+1}\gamma_{n+2}\right)=L\left(\gamma_{1}^{2}\gamma_{n}\gamma_{n+1}\right)=L\left(\gamma_{1}^{2}\gamma_{2n+1}\right)$.
The leading order of $\gamma_{2n+1}$ is $\left(sa\right)^{2n-1}\gamma_{2}$
so that $L\left(\gamma_{1}^{2}\left(sa\right)^{2n-1}\gamma_{2}\right)=L\left(\gamma_{2n+3}\right)$.
Then, because $L(\gamma_{n+1}\gamma_{n+2})=L(\gamma_{1}\gamma_{n}\gamma_{n+2})$,
the claim follows for arbitrary subscripts. If we deal with eigenvectors,
$\alpha_{12}=0$, so that the highest-order term is $\gamma_{1}\gamma_{m-\left(k+2\right)}$,
which by the preceding argument has maximal leading order $\left(as\right)^{k+1}L\left(\gamma_{m-k-1}\right)=L\left(\gamma_{m}\right)$,
which amounts to $L\left(s^{m},a^{m}\right)$ and thus establishes
the claim. To obtain the power of $\alpha_{12}$, consider the terms
$\gamma_{2}$ and $\gamma_{3}$, 
\begin{align*}
\gamma_{2} & =s\left(\beta\gamma_{1}+\alpha_{12}\gamma_{1}^{2}\right)\\
 & =s\beta\gamma_{1}+s\alpha_{12}\gamma_{1}^{2}\\
 & \leq s^{2}a^{2}+s^{3}a^{2}\alpha_{12}\\
\gamma_{3} & =s\left(\beta\gamma_{2}+\beta_{1}\gamma_{1}^{2}+\alpha_{12}\gamma_{1}\gamma_{2}\right)\\
 & =s\left(\beta s\left(\beta\gamma_{1}+\alpha_{12}\gamma_{1}^{2}\right)+\beta_{1}\gamma_{1}^{2}+\alpha_{12}\gamma_{1}s\left(\beta\gamma_{1}+\alpha_{12}\gamma_{1}^{2}\right)\right)\\
 & =s^{2}\beta^{2}\gamma_{1}+s^{2}\beta\alpha_{12}\gamma_{1}^{2}+s\beta_{1}\gamma_{1}^{2}+s^{2}\beta\alpha_{12}\gamma_{1}^{2}+\alpha_{12}^{2}s^{2}\gamma_{1}^{3}\\
 & =s^{2}a^{2}\gamma_{1}+s^{4}a^{3}\alpha_{12}+s^{3}a^{3}+s^{4}a^{3}\alpha_{12}+\alpha_{12}^{2}a^{3}s^{5}
\end{align*}
to deduce that $\gamma_{m}=L\left(\alpha_{12}^{m-1}\right)$ because
going from $\gamma_{m}$ to $\gamma_{m+1}$ increases the power of
$\alpha_{12}$ by one. Then, we obtain for the powers of $a$ and
$s$ a contribution of $\left(as\right)^{m}$ from $\gamma_{m}$.
In addition, however, we pick up a factor $s^{m-1}$ which can be
seen from expanding $\gamma_{m}$, which results in the highest power
of $\alpha_{12}$ being multiplied by $s^{2m-1}$ based on \citet[Eq. 3.4]{Sun1991}.
\[
\underbrace{s^{m-1}}_{\text{from mult. in front.}}\underbrace{a^{m}s^{m}}_{\text{from leading order of \ensuremath{\gamma_{m}}}}
\]
If $\alpha_{12}>0$, $\left(as\right)^{m-1}L\left(\gamma_{2}\right)=L\left(a^{m+1}s^{m+2}\right)$
so that $a^{k}s^{k+1}\alpha_{12}L\left(\gamma_{m-k}\right)=$ $\iota_{m,k}$
is $\left(as\right)^{k+1}$.
\end{proof}

\subsection{\label{subsec:svd-perturbation}Results for singular subspace maps}

To study the derivative of the implicit map that delivers $U_{I}$,
we consider a perturbation to $M$, $\tilde{M}$.

The aim is to construct an expansion for some $\tilde{M}$ such that
$\smlnorm{M-\tilde{M}}$ is small,
\begin{equation}
\Psi\left(\tilde{M}\right)=\dot{\Psi}\left(M-\tilde{M}\right)+O\left(\norm{M-\tilde{M}}\right).\label{eq:svd-perturb}
\end{equation}
Let $D_{\text{d}}$ be diagonal such that $D_{\text{d,}ij}\defeq\left(\iota_{i}^{2}-\iota_{j}^{2}\right)^{-1}$.
Following \citet{4301315}, we define 
\[
K\defeq D_{\text{d}}\cdot\left(U_{I}^{\trans}\left(M-\tilde{M}\right)V_{I}\Sigma_{I}+\Sigma_{I}V_{I}^{\trans}\left(M-\tilde{M}\right)^{\trans}U_{I}\right)
\]
and where $i\neq j,$ $D_{\text{d,}ii}=0$, and $i,j=1,\dots,F$ and
where $\cdot$ denotes the Hadamard product. We obtain the first-order
perturbation term
\begin{equation}
\dot{\Psi}\left(M-\tilde{M}\right)=\schwein\left(U_{I}K-U_{J}U_{J}^{\trans}\left(M-\tilde{M}\right)V_{I}\Sigma_{I}^{-1}\right)\label{eq:frechet-svd}
\end{equation}

\begin{proof}[Proof of \prettyref{lem:jacobian}\prettyref{enu:singular}.]
To identify the Jacobian from \prettyref{eq:frechet-svd}, we need
to apply \citet[Thm. 11]{magnus2019matrix} and \citet[Remark 4.2]{Sun1991}
to \prettyref{eq:svd-perturb}, which we achieve via vectorization
and studying the remainder term. We write for the mapping that represents
the test statistic 
\begin{align*}
 & \Psi\left(M\left(t\right)\right)\\
= & \upsilon_{\perp}^{\trans}U_{I}\left(t\right)\\
= & \upsilon_{\perp}^{\trans}\left(\tilde{U}_{I}-U_{I}\right)\left(t\right).
\end{align*}
To first-order, we expand $\Psi(M(t))$, using the expansion in \citet{4301315}:
\begin{align*}
 & \Psi\left(M\left(t\right)\right)\\
= & \schwein U_{I}\left(D_{\text{d}}\cdot\left(U_{I}^{\trans}E\left(t\right)V_{I}\Sigma_{I}+\Sigma_{I}V_{I}^{\trans}E^{\trans}\left(t\right)U_{I}\right)\right)\\
 & +\schwein U_{J}U_{J}^{\trans}E\left(t\right)V_{I}\Sigma_{I}^{-1}.
\end{align*}
In order to identify the Jacobian of the perturbation expansion of
$\Psi(M(t))$, we need to find $\vek\Psi(M(t))$, for which we write
\begin{align*}
 & \vek\Psi\left(M\left(t\right)\right)\\
= & \vek\schwein U_{I}\left(D_{\text{d}}\cdot\left(U_{I}^{\trans}E\left(t\right)V_{I}\Sigma_{I}+\Sigma_{I}V_{I}^{\trans}E^{\trans}\left(t\right)U_{I}\right)\right)\\
 & +\vek\schwein U_{J}U_{J}^{\trans}E\left(t\right)V_{I}\Sigma_{I}^{-1}.
\end{align*}
We focus on the second term first: 
\begin{align*}
 & \vek\schwein U_{J}U_{J}^{\trans}E\left(t\right)V_{I}\Sigma_{I}^{-1}\\
= & \left(\Sigma_{I}^{-1\trans}V_{I}^{\trans}\otimes\schwein U_{J}U_{J}^{\trans}\right)\vek E\left(t\right)
\end{align*}
Regarding the first term,
\begin{align*}
 & \vek\schwein U_{I}\left(D_{\text{d}}\cdot\left(U_{I}^{\trans}E\left(t\right)V_{I}\Sigma_{I}+\Sigma_{I}V_{I}^{\trans}E^{\trans}\left(t\right)U_{I}\right)\right)I_{F}\\
= & \left(I_{F}\otimes\schwein U_{I}\right)\left[\left(\vek D_{\text{d}}\right)\cdot\vek\left(U_{I}^{\trans}E\left(t\right)V_{I}\Sigma_{I}+\Sigma_{I}V_{I}^{\trans}E^{\trans}\left(t\right)U_{I}\right)\right]\\
= & \left(I_{F}\otimes\schwein U_{I}\right)\left[\left(\vek D_{\text{d}}\right)\cdot\left(\left(\Sigma_{I}^{\trans}V_{I}^{\trans}\otimes U_{I}^{\trans}\right)\vek\left(E\left(t\right)\right)+\left(U_{I}^{\trans}\otimes\Sigma_{I}V_{I}^{\trans}\right)K_{p}\vek E\left(t\right)\right)\right]\\
= & \left(I_{F}\otimes\schwein U_{I}\right)\left[\left(\vek D_{\text{d}}\right)\cdot\left(\left(\left(\Sigma_{I}^{\trans}V_{I}^{\trans}\otimes U_{I}^{\trans}\right)+\left(U_{I}^{\trans}\otimes\Sigma_{I}V_{I}^{\trans}\right)K_{p}\right)\vek E\left(t\right)\right)\right]
\end{align*}
where $I_{F}$ is an $F\times F$ identity matrix and $F$ is the
true rank of $M$ and $K_{p}$ satisfies $\vek E^{\trans}=K_{p}\vek E$.
We now apply \prettyref{lem:hadamard} to $\left(\vek D_{\text{d}}\right)\cdot\left[\left(\Sigma_{I}^{\trans}V_{I}^{\trans}\otimes U_{I}^{\trans}\right)\vek E\left(t\right)\right]$
to obtain $\left[\left(\indic^{\trans}\otimes\vek D_{\text{d}}\right)\cdot\left(\Sigma_{I}^{\trans}V_{I}^{\trans}\otimes U_{I}^{\trans}\right)\right]\vek E\left(t\right)$
and similarly, 
\begin{align*}
 & \left(\vek D_{\text{d}}\right)\cdot\left(U_{I}^{\trans}\otimes\Sigma_{I}V_{I}^{\trans}\right)K_{p}\vek E\left(t\right)\\
 & =\left[\left(\indic^{\trans}\otimes\vek D_{\text{d}}\right)\cdot\left(U_{I}^{\trans}\otimes\Sigma_{I}V_{I}^{\trans}\right)K_{p}\right]\vek E\left(t\right).
\end{align*}
Therefore, we obtain
\begin{align*}
 & \vek\Psi\left(M\left(t\right)\right)\\
= & \left(\Sigma_{I}^{-1\trans}V_{I}^{\trans}\otimes\schwein U_{J}U_{J}^{\trans}\right)\vek E\left(t\right)\\
 & +\left(I_{F}\otimes\schwein U_{I}\right)\left[\left(\indic^{\trans}\otimes\vek D_{\text{d}}\right)\cdot\left(\left(\Sigma_{I}^{\trans}V_{I}^{\trans}\otimes U_{I}^{\trans}\right)+\left(U_{I}^{\trans}\otimes\Sigma_{I}V_{I}^{\trans}\right)K_{p}\right)\right]\vek E\left(t\right).
\end{align*}
Hence, we can write the map $\Psi\left(M\left(t\right)\right)$ in
differential form as $\diff\Psi(M;E)=:\jacsvd\vek E$ for a perturbation
$E$ and identify the Jacobian matrix invoking \citet[Thm. 11]{magnus2019matrix}
as 
\begin{align*}
\jacsvd= & \left(\Sigma_{I}^{-1\trans}V_{I}^{\trans}\otimes\schwein U_{J}U_{J}^{\trans}\right)\\
+ & \left(I_{F}\otimes\schwein U_{I}\right)\left[\left(\indic^{\trans}\otimes\vek D_{\text{d}}\right)\cdot\left(\left(\Sigma_{I}^{\trans}V_{I}^{\trans}\otimes U_{I}^{\trans}\right)+\left(U_{I}^{\trans}\otimes\Sigma_{I}V_{I}^{\trans}\right)K_{p}\right)\right].
\end{align*}
Differentiability of $\Psi$, which follows from the differentiability
of $\psi$, because singular vectors of $M$ are eigenvectors of $MM^{\trans}$.
\end{proof}
\begin{lem}
\label{lem:hadamard}Let $A$, $B$, and $C$ be matrices of dimensions
$a\times b$, $ab\times cd$, and $c\times d$, respectively. Then,
for a matrix of ones of dimension $cd\times1$, $\indic$, we have
the identities
\[
\vek A\cdot\left(B\vek C\right)=\left[\left(\indic^{\trans}\otimes\vek A\right)\cdot B\right]\vek C
\]
and 
\[
\vek\left(A\cdot B\right)=\vek A\cdot\vek B.
\]
\end{lem}
\begin{proof}[Proof of \prettyref{lem:hadamard}.]
Writing out individual entries establishes the result.
\end{proof}
For a single vector of interest, the Jacobian simplifies considerably.
Write for the statistic 
\begin{align*}
 & \Psi_{s}\left(M\left(t\right)\right)\\
= & \schwein u_{i}\left(t\right)\\
= & \upsilon_{\perp}^{\trans}\left(\tilde{u}_{i}-u_{i}\right).\\
= & \schwein U_{I}D_{i}U_{I}^{\trans}E\left(t\right)v_{i}s_{i}\\
+ & \schwein U_{I}D_{i}\Sigma_{I}V_{I}^{\trans}E^{\trans}\left(t\right)u_{i}\\
+ & \schwein U_{J}U_{J}^{\trans}E\left(t\right)v_{i}s_{i}^{-1}.
\end{align*}
Again, we vectorize term by term, to obtain 
\[
\vek\schwein U_{I}D_{i}U_{I}^{\trans}Ev_{i}s_{i}=\left(s_{i}v_{i}^{\trans}\otimes\schwein U_{I}D_{i}U_{I}^{\trans}\right)\vek E\left(t\right)
\]
so that 
\[
\vek\schwein U_{I}D_{i}\Sigma_{I}V_{I}^{\trans}E^{\trans}\left(t\right)u_{i}=\left(u_{i}^{\trans}\otimes\schwein U_{I}D_{i}\Sigma_{I}V_{I}^{\trans}\right)K\vek E\left(t\right)
\]
and $\left(s_{i}^{-1}v_{i}^{\trans}\otimes\schwein U_{J}U_{J}^{\trans}\right)\vek E\left(t\right).$
\begin{align*}
 & \vek\Psi_{s}\left(M\left(t\right)\right)\\
= & \left(s_{i}v_{i}^{\trans}\otimes\schwein U_{I}D_{i}U_{I}^{\trans}\right)\vek E\left(t\right)\\
+ & \left(u_{i}^{\trans}\otimes\schwein U_{I}D_{i}\Sigma_{I}V_{I}^{\trans}\right)K\vek E\left(t\right)\\
+ & \left(s_{i}^{-1}v_{i}^{\trans}\otimes\schwein U_{J}U_{J}^{\trans}\right)\vek E\left(t\right),
\end{align*}
so that we can write the map $\Psi_{s}\left(M\left(t\right)\right)$
in differential form as $\text{d}\Psi_{s}\left(M;E\right)=:B_{s}\vek E$,
again invoking \citet[Thm. 11]{magnus2019matrix} where 
\begin{align*}
 & B_{s}\\
= & \left(s_{i}v_{i}^{\trans}\otimes\schwein U_{I}D_{i}U_{I}^{\trans}\right)\\
+ & \left(u_{i}^{\trans}\otimes\schwein U_{I}D_{i}\Sigma_{I}V_{I}^{\trans}\right)K\\
+ & \left(s_{i}^{-1}v_{i}^{\trans}\otimes\schwein U_{J}U_{J}^{\trans}\right).
\end{align*}

\section{Smoothness of invariant and singular subspace maps\label{app:using-ift}}

Recall the definitions of $\mathscr{M}$ and $\mathscr{S}$ from
\prettyref{assu:invar-sub-ass} and \prettyref{assu:svd-ass}. To
economize on notation, we shall establish the result for invariant
subspace maps and then argue that the same holds for singular subspace
maps via a small extension. As an additional bit of notation, matrices
with a '$0$' subscript denote the truth to help distinguish them
from variable parameters. We have \citet[Lemma B.1]{duffy2020cointegrated},
with a slightly rewritten proof. 
\begin{lem}
\label{lem:open-and-smooth}~
\begin{enumerate}
\item \label{enu:p-is-open}The set $\mathscr{M}$ is open.
\item \label{enu:smooth-maps}The maps $R_{I}\left(M\right)$, $\Lambda_{I}\left(M\right)$
are smooth on $\mathscr{M}$.
\end{enumerate}
\end{lem}
\begin{proof}
\prettyref{lem:open-and-smooth}\prettyref{enu:p-is-open} follows
from the fact that the set of full-rank matrices that $R_{I}$ belongs
to by \prettyref{assu:invar-sub-ass}\prettyref{enu:rankRlu} is open.
Then, by \citet[Thm. IV.1.1 and V.2.8]{Stewart1991} we have that
eigenvalues and simple invariant subspaces are continuous. In other
words, we have that $\forall\epsilon>0$, $\exists\delta>0$ such
that $\norm{M_{0}-M}<\delta$ implies $\norm{R_{0,I}-R_{I}}<\epsilon$.
Let the eigenvector map be defined by $f\left(M\right)=R_{I}$. Then,
we can show that the set of matrices $\mathscr{M}$ is open. The set
of full rank matrices $\mathcal{R}$ to which $R_{I}$ belongs is
open, therefore, for every $R\in\mathcal{R}$, $B\left(R,\delta\right)\subset\mathcal{R}$
for some $\delta>0$. Let $f^{-1}$ denote the pre-image rather than
inverse. Then, we have a small $r$ such that $B\left(f\left(M_{0}\right),r\right)\subset\mathcal{R}$
because for some $M_{0}\in f^{-1}\left(\mathcal{R}\right)$ it must
be that $f\left(M_{0}\right)\in\mathcal{R}$.

Continuity implies that for all $\epsilon>0$, there is a $\delta>0$
such that $M\in B\left(M_{0},\delta\right)\implies R\in B\left(f\left(M_{0}\right),\epsilon\right)\subset\mathcal{R}$,
where the last subset relation follows from the fact that we can choose
$\epsilon\leq r$ so that by openness of $\mathcal{R}$, it must contain
all its open balls such that the $\epsilon$-ball is wholly contained
in the $r$-ball easily found in $\mathcal{R}$. Formally, $R\in B\left(R_{0},\epsilon\right)\subseteq B\left(R_{0},r\right)\subseteq\mathcal{R}$.
Then, carrying on with the continuity argument, the implied relation
$R\in B\left(f\left(M_{0}\right),\epsilon\right)\subset\mathcal{R}$
furthermore implies $R=f\left(M\right)\in B\left(f\left(M_{0}\right),\epsilon\right)\subset\mathcal{R}$
so that $f\left(M\right)\in\mathcal{R}$ or $M\in f^{-1}\left(\mathcal{R}\right)$.
From before, we also have that $M\in B\left(M_{0},\delta\right)$Rewriting
the continuity implication from $M\in B\left(M_{0},\delta\right)\implies R\in B\left(f\left(M_{0}\right),\epsilon\right)\subset\mathcal{R}$
to $M\in B\left(M_{0},\delta\right)\implies M\in f^{-1}\left(\mathcal{R}\right)$,
i.e. if $M$ is in an open $\delta$-ball, it will also be in the
pre-image of $\mathcal{R}$, then we must have, that (since $M_{0}$
was arbitrary) the pre-image contains all open balls and is hence
open so that $\mathscr{M}$ is open.

To establish smoothness, consider the maps whose zero level sets define
the invariant subspaces~of interest 
\begin{align}
H^{\ast}\left(R_{I},\Lambda_{I};M\right) & \defeq\begin{bmatrix}R_{I}\Lambda_{I}-MR_{I};L_{0,I}^{\trans}R_{I}-I_{q}\end{bmatrix}=\begin{bmatrix}H_{1}^{\ast};H_{2}^{\ast}\end{bmatrix}\label{eq:h-star}\\
H\left(R_{I},\Lambda_{I};M\right) & \defeq\begin{bmatrix}R_{I}\Lambda_{I}-MR_{I};G^{\trans}R_{I}-I_{q}\end{bmatrix},\label{eq:h}
\end{align}
for $H,H^{\ast}:\reals^{p\times q}\times\reals^{q\times q}\rightarrow\reals^{\left(p+q\right)\times q}$.
Note that $H^{\ast}\left(R_{I,0},\Lambda_{I,0};M_{0}\right)=H^{\ast}\left(R_{I,0},\Lambda_{I,0};M_{0}\right)=0$,
i.e. the two maps agree at the truth. We now need to compute the Jacobian
of $H^{\ast}$ and show that it is non-singular in which case \citet[Thm.~XIV.2.1]{Lang93}
establishes that there exists a neighborhood $\mathscr{N}$ of $\mathscr{M}$
and smooth maps $R_{I}^{\ast}:\mathscr{N}\rightarrow\reals^{p\times q}$
and $\Lambda_{I}^{\ast}:\mathscr{N}\rightarrow\reals^{q\times q}$
such that $H^{\ast}\left(R_{I}^{\ast}\left(M\right),\Lambda^{\ast}\left(M\right);M\right)=0$
for all $M\in\mathscr{M}$. By the continuity of $R_{I}^{\ast}\left(\cdot\right)$
we may choose $\mathscr{N}$ such that $\rank G^{\trans}R_{I}=q$.
A non-singular Jacobian obtains if $\diff H^{\ast}=0\implies\diff R_{I}^{\ast}=0\,\wedge\diff\Lambda_{I}^{\ast}=0$.

Using matrix differentiation while holding $M$ at $M_{0}$ and evaluating
the other terms at the truth, too, yields $\diff H_{1}\defeq\diff R_{I}\Lambda_{I,0}+R_{I,0}\diff\Lambda_{I}-M_{0}\diff R_{I}$.
Expand $M_{0}\diff R_{I}=R_{I,0}\Lambda_{I,0}L_{I,0}^{\trans}\diff R_{I}+R_{J,0}\Lambda_{J,0}L_{J,0}^{\trans}\diff R_{I}$
and observe that the normalization of the second component in $H_{1}^{\ast}$
implies that $L_{I,0}^{\trans}\diff R_{I}=0$. Let $\diff\Lambda_{I}\leftarrow\diff\Lambda_{I}-\Lambda_{IJ,0}L_{J,0}^{\trans}\diff R_{I}$.
Biorthogonality would only establish $L_{J,0}^{\trans}R_{J,0}=I_{p-q}$.
Because $L_{I,0}^{\trans}\diff R_{I}=0$, i.e. $\diff R_{I}$ is orthogonal
to the columns of $L_{I,0}$, we have $\diff R_{I}\in\col R_{J,0}$
so that the projector $R_{J,0}L_{J,0}^{\trans}$ as an identity on
$\diff R_{I}$. Substituting these results into the previous expression
yields 
\begin{align}
\diff H_{1} & =\diff R_{I}\Lambda_{I,0}+R_{I,0}\diff\Lambda_{I}-R_{J,0}\Lambda_{J,0}L_{J,0}^{\trans}\diff R_{I}\nonumber \\
 & =_{\left(1\right)}R_{I,0}\diff\Lambda_{I}+R_{J,0}L_{J,0}^{\trans}\diff R_{I}\Lambda_{I,0}-R_{J,0}\Lambda_{J,0}L_{J,0}^{\trans}\diff R_{I}\nonumber \\
 & =R_{I,0}\diff\Lambda_{I}+R_{J,0}\left(L_{J,0}^{\trans}\diff R_{I}\Lambda_{I,0}-\Lambda_{J,0}L_{J,0}^{\trans}\diff R_{I}\right)\nonumber \\
 & =_{\left(2\right)}R_{I,0}\diff\Lambda_{I}+R_{J,0}\boldsymbol{S}\left(L_{J,0}^{\trans}\diff R_{I}\right)\label{eq:dh1-nonsing}
\end{align}
where $=_{\left(1\right)}$ follows from $R_{J,0}L_{J,0}^{\trans}=I_{p-q}$
and $=_{\left(2\right)}$ from \prettyref{eq:lin-op}.

Then, \prettyref{assu:invar-sub-ass}\prettyref{enu:rankRlu} implies
that $\diff\Lambda_{I}=0$ for the first term in \prettyref{eq:dh1-nonsing}.
Moreover, \prettyref{assu:invar-sub-ass}\prettyref{enu:no-evs} implies
that $\boldsymbol{S}\left(.\right)$ is non-singular \citep[Thm V.1.3]{Stewart1991}
and hence that we must have $\left[\begin{smallmatrix}L_{0,I} & L_{0,J}\end{smallmatrix}\right]^{\trans}\diff R_{I}=0$
which implies that $\diff R_{I}=0$ because \prettyref{assu:invar-sub-ass}\prettyref{enu:semisimple}
implies that $R_{0}$ has full rank and so must $L_{0}$.

Having established the result for $H^{\ast}$, the corresponding result
for $H$ now follows from the fact that the maps $R_{I}\left(M\right)$
and $\Lambda_{I}\left(M\right)$ can be obtained via observing that
\[
R_{I}\left(M\right)=R_{I}^{\ast}\left(G^{\trans}R_{I}^{\ast}\right)^{-1}
\]
obeys the normalization $G^{\trans}R_{I}\left(M\right)=I_{q}$, which
leaves invariant subspaces of interest undisturbed according to 
\[
\underset{R_{I}}{\underbrace{R_{I}^{\ast}\left(G^{\trans}R_{I}^{\ast}\right)^{-1}}}\underset{\Lambda_{I}}{\underbrace{\left(G^{\trans}R_{I}^{\ast}\right)\Lambda_{I}^{\ast}\left(G^{\trans}R_{I}^{\ast}\right)^{-1}}}\underset{\left(R_{I}\right)^{-1}}{\underbrace{\left(G^{\trans}R_{I}^{\ast}\right)\left(R_{I}^{\ast}\right)^{-1}}}
\]
so that we arrive at the renormalized maps
\begin{align*}
R_{I}\left(M\right) & =R_{I}^{\ast}\left(G^{\trans}R_{I}^{\ast}\right)^{-1}\\
\Lambda_{I}\left(M\right) & =\left(G^{\trans}R_{I}^{\ast}\right)\Lambda_{I}^{\ast}\left(G^{\trans}R_{I}^{\ast}\right)^{-1}.
\end{align*}
\end{proof}
\begin{lem}
\label{lem:svd-smooth}~
\begin{enumerate}
\item \label{enu:open-svd}The set $\mathscr{S}$ is open.
\item \label{enu:smooth-svd}The maps $U_{I}\left(M\right)$, $V_{I}\left(M\right)$,
and $\Sigma_{I}\left(M\right)$ are smooth on $\mathscr{S}$.
\end{enumerate}
\end{lem}
\begin{proof}
For \prettyref{lem:svd-smooth}\prettyref{enu:open-svd}, we observe
that there exists a continuous inverse that maps a pair $U_{I}$ and
$\Sigma_{I}$ to $MM^{\trans}$ and that the set of full-rank matrices
$U_{I}$ belongs to is open. 

For \prettyref{lem:svd-smooth}\prettyref{enu:smooth-svd}the case
of singular subspaces, we observe that the maps $H$ and $H^{\ast}$
establish the result for general invariant subspace maps. In the following,
restrict attention to the set of non-redundant columns as allowed
by \prettyref{assu:svd-ass}\prettyref{enu:svd-split}, i.e. $M=U_{I}\Sigma_{I}V_{I}^{\trans}$
so that we can ignore information in $M$ associated with zero singular
values. For $M\in\mathscr{M}$, consider the map $h:M\mapsto MM^{\trans}$
and observe that the map $H\circ h\left(M\right)$ delivers the result
for right-singular subspaces ($\spn U_{I}$). Changing $h$ to $h^{\trans}:M\mapsto M^{\trans}M$
delivers the result for left-singular subspaces as per the proof of
\prettyref{lem:open-and-smooth}.
\end{proof}

\subsection{\label{app:smooth-maps-perturb-expansions}Relationship between smooth
maps and perturbation expansions}

We begin by proving \prettyref{thm:resolvent}.
\begin{proof}[Proof of \prettyref{thm:resolvent}.]
To apply \citet[Thm.~XIV.2.1]{Lang93}, we need to show that $\diff K=0\implies\diff X=0$,
i.e. that the associated Jacobian is non-singular. Note that
\[
\diff K=-\left(X-\alpha I\right)^{-1}\left(\diff X\right)\left(X-\alpha I\right)^{-1}.
\]
Setting $\diff K=0$ implies $-\left(X-\alpha I\right)^{-1}\left(\diff X\right)\left(X-\alpha I\right)^{-1}=0$.
Clearly, $\left(X-\alpha I\right)^{-1}=0$ is impossible, so $\diff X=0$.
Alternatively, the Jacobian $J=\left[\left(X^{\trans}-\alpha I\right)^{-1}\otimes\left(X-\alpha I\right)^{-1}\right]\in\reals^{p^{2}\times p^{2}}$
can be identified from $\vek\diff K=J\diff\vek X$ and has to be non-singular
by the fact that $\rank\left(A\otimes B\right)=\rank A\rank B$ and
$\alpha\notin\eigs_{I}\union\eigs_{J}$.
\end{proof}
\prettyref{lem:open-and-smooth} established the result of infinite
differentiability via implicit differentiation and the fact that the
Jacobians of $H$ and $H^{\ast}$ have to be non-singular. As part
of this argument, we showed that the operator in \prettyref{eq:lin-op}
is non-singular. That same operator appears in \citet[Eq. 2.2 and p. 90, second display]{Sun1991}. 

The derivation in \citet{Sun1991} uses the following $\tilde{M}\left(t\right)\defeq L^{\trans}M\left(t\right)R$,
where if $L$ and $R$ are left- and right-invariant matrices, the
lower left $q\times\left(p-q\right)$ block ought to be zero. The
proof of \citet[Thm. 2.1]{Sun1991}, (his main result), proceeds to
study the conditions necessary for the equality $L_{I}^{\trans}M\left(t\right)R_{J}=0$
to hold and concludes that the operator in \prettyref{eq:lin-op}
has to be non-singular, i.e. $\norm{\left(\Lambda{}_{I}^{\trans}\otimes I_{J}\right)-\left(I_{I}\otimes\Lambda_{J}\right)}\neq0$
to achieve analyticity. Our result \prettyref{lem:open-and-smooth}\prettyref{enu:smooth-maps}
implies the analyticity result of \citet{Sun1991} and extends it
towards normalized bases as per the second component of \prettyref{eq:h-star}.
Therefore, the approach of appealing to implicit maps generalizes
the result of \citet{Sun1991} to that of normalized maps.

\section{Results and proofs for \prettyref{sec:centralities}: \nameref{sec:centralities}\label{app:results-network-centralities}}

\subsection{Convergence rates of network statistics}
\begin{proof}[Proof of \prettyref{thm:stoch-bound-prop}]
We write the perturbation expansion as

\begin{align*}
\psi\left(\hat{M}\right) & =\psi\left(M+O_{p}\left(r_{n}\right)\right)\\
 & =\psi\left(M\right)+\dot{\psi}\left(X;\Delta_{\text{JI}}\right)\mid_{X=M}^{\trans}\left(\vek\left(\hat{M}-M\right)\right)+o_{p}\left(\hat{M}-M\right)\\
 & =\psi\left(M\right)+\dot{\psi}\left(X;\Delta_{\text{JI}}\right)\mid_{X=M}^{\trans}\left(\vek\left(O_{p}\left(r_{n}\right)\right)\right)+o_{p}\left(\hat{M}-M\right)
\end{align*}
From the differentiability of $\psi\left(M\right)$ we have 
\begin{align*}
 & \psi\left(\hat{M}\right)-\psi\left(M\right)-\dot{\psi}\left(X;\Delta_{\text{JI}}\right)\mid_{X=M}^{\trans}\left(\vek\left(\hat{M}-M\right)\right)\\
 & =R\left(\norm{\hat{M}-M}\right)\\
 & =o_{p}\left(\norm{\hat{M}-M}\right)
\end{align*}
Define
\begin{align*}
\frac{1}{\delta} & \defeq\norm{\left[\left(\Lambda{}_{I}^{\trans}\otimes I_{J}\right)-\left(I_{I}\otimes\Lambda_{J}\right)\right]^{-1}}_{\text{op}}\\
\gamma_{1} & \defeq\frac{\beta_{2}}{\delta}\\
\beta_{2} & \defeq\norm{B_{21}}_{F}\\
\alpha_{12} & \defeq\norm{A_{12}}_{2}\\
\beta & \defeq\norm{B_{11}}_{2}+\norm{B_{22}}_{2}.
\end{align*}
We have from \citet[Thm. 3.1]{Sun1991} to second order 
\begin{align*}
d_{F}\left(M,\hat{M}\right) & \leq\gamma_{1}\abs t+\gamma_{2}\abs t^{2}+\underbrace{O\left(\abs t^{3}\right)}_{\text{bounded above by \ensuremath{\abs t^{3}}}}.\\
 & =\gamma_{1}\abs t+\frac{1}{\delta}\left(\beta\gamma_{1}+\alpha_{12}\gamma_{1}^{2}\right)\abs t^{2}\\
 & =\gamma_{1}\abs t+\gamma_{1}\abs t\frac{1}{\delta}\left(\beta+\alpha_{12}\gamma_{1}\right)\abs t\\
 & =\gamma_{1}\abs t+\gamma_{1}\abs t\frac{1}{\delta}\left(\beta+\alpha_{12}\gamma_{1}\right)\abs t\\
 & =\gamma_{1}\abs t\left(1+\frac{1}{\delta}\left(\beta+\alpha_{12}\gamma_{1}\right)\abs t\right)
\end{align*}
Recall that we model $\hat{M}=M+tE$ for a perturbing matrix $E$
with unit norm. Furthermore, evaluate $\psi$ in the $(L,R)$ basis
so that $E=L^{\trans}(\hat{M}-M)R$ and $\smlnorm{B_{21}}_{\text{F}}=\smlnorm{E_{21}}_{\text{F}}$.
Alternatively, we could appeal to the argument in \prettyref{app:inequality-block-diagonal},
and invoke \citet[Lemma A1]{demetrius}. We thus obtain 
\begin{align*}
R\left(\norm{\hat{M}-M},t\right) & \leq\gamma_{1}\frac{1}{\delta}\left(\beta+\alpha_{12}\gamma_{1}\right)\abs t^{2}+O\left(\abs t^{3}\right)\\
 & =\frac{\beta_{2}}{\delta^{2}}\left(\norm{B_{11}}_{2}+\norm{B_{22}}_{2}+\norm{A_{12}}_{2}\frac{\beta_{2}}{\delta}\right)\abs t^{2}+O\left(\abs t^{3}\right)\\
 & =\frac{\norm{B_{21}}_{F}}{\delta^{2}}\left(\norm{B_{11}}_{2}+\norm{B_{22}}_{2}+\norm{A_{12}}_{2}\frac{\norm{B_{21}}_{F}}{\delta}\right)\abs t^{2}+O\left(\abs t^{3}\right)\\
 & =\frac{\norm{E_{21}}_{F}}{\delta^{2}}\left(\norm{E_{11}}_{2}+\norm{E_{22}}_{2}+\norm{A_{12}}_{2}\frac{\norm{E_{21}}_{F}}{\delta}\right)\abs t^{2}+O\left(\norm{E_{21}}_{F}^{3}\abs t^{3}\right)
\end{align*}
where $\abs t\goesto0$. The above derivation used Sun's perturbation
estimations, equivalent to multiplying $d_{F}\left(M,\hat{M}\right)$
by the second-order Fréchet derivative $\gamma_{1}/\delta$ . It is
now easy to see that, recalling $\smlnorm{\hat{M}-M}_{2}=\smlnorm{tE}_{2}$,
implies 
\[
R\left(\norm{\hat{M}-M}_{2}\right)=o\left(\norm{\hat{M}-M}_{2}\right)
\]
as $\norm{tE}_{2}\goesto0$ because
\begin{align*}
\frac{R\left(\norm{\hat{M}-M}_{2}\right)}{\norm{\hat{M}-M}_{2}} & \leq\frac{\norm{E_{21}}_{F}}{\norm{tE}_{2}\delta^{2}}\left(\norm{E_{11}}_{2}+\norm{E_{22}}_{2}+\norm{A_{12}}_{2}\frac{\norm{E_{21}}_{F}}{\delta}\right)\abs t^{2}+O\left(\abs t^{3}\right).\\
 & \leq\frac{\norm{E_{21}}_{F}}{\norm{tE}_{2}\delta^{2}}\left(\norm{E_{11}}_{2}+\norm{E_{22}}_{2}+\norm{A_{12}}_{2}\frac{\norm{E_{21}}_{F}}{\delta}\right)\abs t^{2}+O\left(\abs t^{3}\right)\\
 & \leq\frac{\norm{E_{21}}_{F}}{\norm E_{2}\delta^{2}}\left(\norm{E_{11}}_{2}+\norm{E_{22}}_{2}+\norm{A_{12}}_{2}\frac{\norm{E_{21}}_{F}}{\delta}\right)\abs t+O\left(\abs t^{2}\right)\\
 & \goesto0.
\end{align*}
as $t\goesto0$, provided all other quantities are finite, i.e. the
perturbing matrix $E$ is such that $\frac{\norm{E_{21}}_{F}}{\norm E_{2}}<\infty$
and $\delta>0$. Therefore, we have established Fréchet differentiability.
We can thus approximate 
\[
\vek\left(\hat{M}-M\right)\mapsto\psi_{1}\left(\hat{M}\right)-\psi_{1}\left(M\right)
\]
by 
\[
\vek\left(\hat{M}-M\right)\mapsto\dot{\psi}_{1}\left(X;\Delta_{\text{JI}}\right)\mid_{X=M}^{\trans}\left(\vek\left(\hat{M}-M\right)\right)+R\left(\norm{\hat{M}-M}\right)
\]
where we have taken the approximation error estimate from the perturbation
estimation. Note, that we are suppressing dependence on the order
parameter $t$. To calculate sharp bounds on the convergence of $\psi\left(\hat{M}\right)-\psi\left(M\right)$,
we need dependence on the (norm of) the estimation error $\hat{M}-M$
alone. To first-order, we have
\begin{align*}
\norm{\psi_{1}\left(\hat{M}\right)-\psi_{1}\left(M\right)}_{2} & =\norm{\dot{\psi}_{1}\left(M;\Delta_{\text{JI}}\right)\mid_{M=M}^{\trans}\left(\vek\left(\hat{M}-M\right)\right)}_{2}+R\left(\norm{\hat{M}-M}\right)\\
 & =\norm{R_{J}\left(\lambda_{1}I_{p-q}-\Lambda_{J}\right)^{-1}L_{J}^{\trans}\left(M-\hat{M}\right)R_{I}}_{2}+R\left(\norm{\hat{M}-M}\right)\\
 & \leq_{\left(1\right)}\norm{R_{J}\left(\lambda_{1}I_{p-q}-\Lambda_{J}\right)^{-1}}_{2}\norm{L_{J}^{\trans}\left(M-\hat{M}\right)R_{I}}_{2}+R\left(\norm{\hat{M}-M}\right)\\
 & \leq_{\left(2\right)}\norm{\left(\lambda_{1}I_{p-q}-\Lambda_{J}\right)^{-1}}_{2}\norm{\left(M-\hat{M}\right)}_{2}+R\left(\norm{\hat{M}-M}\right),
\end{align*}
where $\leq_{\left(1\right)}$ follows from the sub-multiplicativity
of the Euclidean norm and $\leq_{\left(2\right)}$ follows from \citet[Lemma A1]{demetrius}.
Note that $AR_{J}=R_{J}\Lambda_{J}$ and similarly for $R_{I}$ where
we can always choose basis vectors such that $R_{J},$ $L_{J}$, and
$R_{I}$ have orthonormal columns using the Gram-Schmidt procedure
without altering the span of the invariant subspaces. Doing so alters
the bound by a fixed constant $\smlnorm{L_{J}^{\trans}Z_{1}}_{2}\smlnorm{Z_{2}R_{I}}_{2}$
for orthogonalising matrices $Z_{1/2}$. Finally, the term 
\begin{align}
\norm{\left[\left(\Lambda{}_{I}^{\trans}\otimes I_{J}\right)-\left(I_{I}\otimes\Lambda_{J}\right)\right]^{-1}}_{2} & =_{\left(1\right)}\max\left\{ \left(\lambda_{r-1}-\lambda_{r}\right)^{-1},\left(\lambda_{s}-\lambda_{s+1}\right)^{-1}\right\} \nonumber \\
 & =\frac{1}{\min\left\{ \left(\lambda_{r-1}-\lambda_{r}\right),\left(\lambda_{s}-\lambda_{s+1}\right)\right\} }\label{eq:sylvester}
\end{align}
where $=_{\left(1\right)}$ follows from the fact that if our vectors
are measured with the spectral (2) norm. For the principal component
expansion, we have $\norm{\left(\lambda_{1}I_{p-q}-\Lambda_{J}\right)^{-1}}_{2}=\frac{1}{\lambda_{1}-\lambda_{2}}$.
\end{proof}
\begin{proof}[Proof of \prettyref{thm:convergence-rate-clustering-coeff}]
Observe that 
\begin{align*}
\norm{M-\hat{M}}_{\text{F}} & =\sqrt{\sum_{i=1}^{p}\sum_{j=1}^{p}\abs{m_{ij}-\hat{m}_{ij}}^{2}}\\
 & =O_{p}\left(r_{n}\right)
\end{align*}
implies that $\abs{m_{ij}-\hat{m}_{ij}}=O_{p}\left(r_{n}\right)$
a fortiori. To find the convergence rate of the clustering coefficient,
we will work out denominator and numerator separately. Consider the
identity 
\begin{align}
\left(a-\hat{a}\right)\left(b-\hat{b}\right) & =ab-a\hat{b}-\hat{a}b+\hat{a}\hat{b}\pm ab\nonumber \\
 & =\left(\hat{a}\hat{b}-ab\right)+2ab-a\hat{b}-\hat{a}b\nonumber \\
 & =\left(\hat{a}\hat{b}-ab\right)+a\left(b-\hat{b}\right)+ab-\hat{a}b\nonumber \\
 & =\left(\hat{a}\hat{b}-ab\right)+a\left(b-\hat{b}\right)+\left(a-\hat{a}\right)b,\label{eq:identity}
\end{align}
so that $\left(\hat{a}\hat{b}-ab\right)=\left(a-\hat{a}\right)\left(b-\hat{b}\right)-a\left(b-\hat{b}\right)-\left(a-\hat{a}\right)b$.
We apply \eqref{eq:identity} to the centered estimator $\hat{m}_{ij}\hat{m}_{ik}-m_{ij}m_{ik}$
to obtain
\begin{align*}
m_{ij}m_{ik}-\hat{m}_{ij}\hat{m}_{ik} & =m_{ij}\left(m_{ik}-\hat{m}_{ik}\right)+\left(m_{ij}-\hat{m}_{ij}\right)m_{ik}-\left(m_{ij}-\hat{m}_{ij}\right)\left(m_{ik}-\hat{m}_{ik}\right)\\
 & =O_{p}\left(r_{n}\right)
\end{align*}
based on the fact that $m_{ik}-\hat{m}_{ik}=O_{p}\left(r_{n}\right)$.
Now, put $\bar{m}_{ijk}\defeq m_{ij}m_{ik}$ and consider the centered
numerator $\hat{\bar{m}}_{ijk}\hat{m}_{jk}-\bar{m}_{ijk}m_{jk}$.
Another application of the preceding argument to the centered numerator
of \eqref{eq:clustering-coeff} yields that $\hat{\bar{m}}_{ijk}\hat{m}_{jk}-\bar{m}_{ijk}m_{jk}=O_{p}\left(r_{n}\right)$.
To find the stochastic order of $\hat{\text{\text{cl}}}_{i}-\text{cl}_{i}$,
put $\bar{m}_{ijk}m_{jk}=n_{i}$ and $m_{ij}m_{ik}=d_{i}$ so that
\begin{align*}
\frac{m_{ij}m_{ik}m_{jk}}{m_{ij}m_{ik}}-\frac{\hat{m}_{ij}\hat{m}_{ik}\hat{m}_{jk}}{\hat{m}_{ij}\hat{m}_{ik}} & =\frac{\hat{n}_{i}d_{i}-n_{i}\hat{d}_{i}\pm\hat{n}_{i}\hat{d}_{i}}{\hat{d}_{i}d_{i}}\\
 & =\frac{\hat{n}_{i}\left(d_{i}-\hat{d}_{i}\right)}{d_{i}\hat{d}_{i}}-\frac{\left(n_{i}-\hat{n}_{i}\right)}{d_{i}}\\
 & =\underbrace{\frac{\hat{n}_{i}}{\hat{d}_{i}}}_{\text{A}}\underbrace{\left(1-\frac{\hat{d}_{i}}{d_{i}}\right)}_{\text{B}}-\underbrace{\frac{1}{d_{i}}\left(n_{i}-\hat{n}_{i}\right)}_{\text{C}}
\end{align*}
For term $\text{A}$, we expand the denominator like so $\frac{1}{\hat{d}_{i}}=\frac{1}{d_{i}+O_{p}\left(r_{n}\right)}=\frac{1}{d_{i}}\frac{1}{1+O_{p}\left(\frac{r_{n}}{d_{i}}\right)}$,
where we find a first-order Maclaurin expansion to obtain $\frac{1}{\hat{d}_{i}}=\frac{1}{d_{i}}\left(1-O_{p}\left(\frac{r_{n}}{d_{i}}\right)\right)$.
Given that $n_{i}=n+O_{p}\left(r_{n}\right)$ we have that 
\begin{align*}
\text{A} & =\hat{n}_{i}d_{i}^{-1}\left(1-O_{p}\left(\frac{r_{n}}{d_{i}}\right)\right)\\
 & =\left(O_{p}\left(r_{n}\right)+n_{i}\right)d_{i}^{-1}\left(1-O_{p}\left(\frac{r_{n}}{d_{i}}\right)\right)\\
 & =O_{p}\left(r_{n}\right)-O_{p}\left(\frac{r_{n}^{2}}{d_{i}^{2}}\right)+\frac{n_{i}}{d_{i}}-O_{p}\left(r_{n}\right)\left(\frac{n_{i}}{d_{i}^{2}}\right)\\
 & =\frac{n_{i}}{d_{i}}+O_{p}\left(r_{n}\right).
\end{align*}
Then, term $\text{B}=O_{p}\left(r_{n}\right)$ so that $\text{AB}=O_{p}\left(r_{n}\right)$.
Similarly, $\text{C}=O_{p}\left(r_{n}\right)$. In light of the fact
that $a_{n},b_{n}=O_{p}\left(r_{n}\right)$ implies that $\frac{\sum^{n}a_{n}+k_{1}}{\sum^{n}b_{n}+k_{2}}=\frac{k_{1}}{k_{2}}+O_{p}\left(r_{n}\right)$
so that $\hat{\text{\text{cl}}}_{i}-\text{cl}_{i}=O_{p}\left(r_{n}\right)$.
A similar argument establishes the claim for $\hat{\text{\text{cl}}}-\text{cl}$.
\end{proof}

\subsection{Network centrality measures\label{subsec:network-centrality-proofs}}

\begin{proof}[Proof of \prettyref{lem:equivalence-centr}.]
By the von Neumann series, $\underset{S\goesto\infty}{\lim}\sum_{s=0}^{S}\alpha^{s}M^{s}=\left(I-\alpha M\right)^{-1}$
which establishes \prettyref{lem:equivalence-centr}\prettyref{enu:diff-katz-equi}.
For \prettyref{lem:equivalence-centr}\prettyref{enu:eigenvec-equiv},
the trick is to multiply by $\lambda_{1}-\alpha$ before taking the
limit and observing that this transformation only scales the leading
eigenvector and $\lambda_{1}\goesto\alpha$ will zero out all other
eigendirections.
\end{proof}
\begin{proof}[Proof of \prettyref{thm:pagerank-katz-smooth}.]
\prettyref{thm:pagerank-katz-smooth}\prettyref{enu:pr-katz} follows
directly from differentiating 
\[
c_{P}\left(M\right)=e_{i}^{\trans}\left(M-\alpha I\right)^{-1}\beta
\]
with respect to $M$ and \prettyref{thm:resolvent}. \prettyref{thm:pagerank-katz-smooth}\prettyref{enu:diffusion}
follows from observing that $\diff\left(M^{f}\right)=\sum_{j=1}^{f}M{}^{j-1}\left(\diff M\right)M^{f-j}$
and writing out the definition of the diffusion centrality.
\end{proof}
\begin{proof}[Proof of \prettyref{thm:deg-centrality}]
Vectorize $c_{N}=M\indic$ and the result follows.
\end{proof}

\begin{algorithm}
\caption{\label{alg:wald-and-t}DGP underlying results in \prettyref{fig:quality-approx-invar-sub}.}
\centering
\begin{algorithmic} [1]
\REQUIRE Dimension of the matrix $p$, sample size $n$, Monte Carlo repetitions $K$, Dimension of invariant subspace of interest $q$, a candidate matrix $M$, and a null hypothesis $\upsilon_{\perp} = \begin{psmallmatrix} I_{r} & -D^\trans \end{psmallmatrix}$
\ENSURE $K$ samples of $t$ and Wald statistic.

\STATE Generate $n$ samples of $M_t = M+E_t \in \reals^{p \times p}$, $K$ times, where the columns of $E$ are drawn from a mean zero multivariate normal distribution with covariance matrix $\Omega_W$. Store the result.
\FOR{$i=1$ \TO $K$}
	\STATE Estimate $\hat{M}_n$ via $n^{-1} \sum_{t=1}^{n} M_t$.
	\STATE Estimate $\hat{R}_n$ and $\hat{L}_n$ and partition them like in  \prettyref{eq:BlocksOfR}, $R_{I} = R\left[ :,1:q \right]$ and $R_{J} = R\left[ :,q+1: p  \right]$. Find $L = \begin{psmallmatrix} L_{I} & L_{J} \end{psmallmatrix}$ via (generalized) inversion.
	\STATE Construct $\hat{D}$ based on \prettyref{eq:ahat-1}, take the $\left( i,j \right)$ element for $d_{i,j}$.
	\STATE Construct Wald or $t$-statistic for null hypothesis $\upsilon$ based on \prettyref{eq:waldstat-1} or \prettyref{eq:tstat-1-1}, respectively.
	\STATE For the $t$-statistic, take the real part only. The Wald statistic will necessarily be real.
\ENDFOR

\end{algorithmic}

\end{algorithm}
\begin{algorithm}
\begin{algorithmic} [1]
\REQUIRE Dimension of the matrix $m \times l$, sample size $n$, Monte Carlo repetitions $K$, Dimension of singular subspace of interest $F$, (if applicable) a singular vector of interested indexed by $f$, a candidate matrix $M$, and a null hypothesis $\upsilon_{\perp} = \begin{psmallmatrix} I_{r} & -D^\trans \end{psmallmatrix}$
\ENSURE $K$ samples of $t$ and Wald statistic.

\STATE Generate $n$ samples of $M_t = M+E_t \in \reals^{m \times l}$, $K$ times, where the columns of $E$ are drawn from a mean zero multivariate normal distribution with covariance matrix $\Omega_W$. Store the result.
\FOR{$i=1$ \TO $K$}
	\STATE Estimate $\hat{M}_n$ via $n^{-1} \sum_{t=1}^{n} M_t$.
	\STATE Estimate $\hat{U}_n$ and $\hat{V}_n$ and partition them like in  \prettyref{eq:BlocksOfR}, $U_{I} = U\left[ :,1:q \right]$ and $U_{J} = U\left[ :,q+1: p  \right]$.
	\STATE Construct $\hat{D}$ based on \prettyref{eq:ahat-1}, take the $\left( i,j \right)$ element for $d_{i,j}$.
	\STATE Construct Wald or $t$-statistic for null hypothesis $\upsilon$ based on \prettyref{eq:waldstat-1} or \prettyref{eq:tstat-1-1}, respectively.
	\STATE For the $t$-statistic, take the real part only. The Wald statistic will necessarily be real.
\ENDFOR

\end{algorithmic}

\caption{\label{alg:wald-and-t-svd}DGP underlying results in \prettyref{fig:quality-approx-svd}.}
\end{algorithm}

\end{document}